\documentclass[11pt]{article}%
\usepackage{amsmath,amsthm,amssymb,amsfonts}
\usepackage{hyperref}
\usepackage[authoryear]{natbib}
\usepackage{multibib}
\usepackage[english]{babel}
\usepackage{xcolor}
\usepackage{geometry}
\usepackage{multirow}
\usepackage{authblk}
\usepackage{amsmath}
\usepackage{amsfonts}
\usepackage{amssymb}
\usepackage{graphicx}%
\providecommand{\U}[1]{\protect\rule{.1in}{.1in}}
\newcites{supp}{Supplemental Material References}
\hypersetup{
colorlinks,
linkcolor={red!50!black},
citecolor={blue!50!black},
urlcolor={blue!80!black}
}
\newtheorem{theorem}{Theorem}
\newtheorem{proposition}{Proposition}
\newtheorem{corollary}{Corollary}
\newtheorem{lemma}{Lemma}
\newtheorem{assumption}{Assumption}
\newtheoremstyle{style2}{4pt}{4pt}{}{}{\bfseries}{.}{ }{}
\theoremstyle{style2}
\newtheorem{remark}{Remark}
\newtheorem{example}{Example}
\begin{document}

\title{Maximum Likelihood Estimation in Markov Regime-Switching Models with
Covariate-Dependent Transition Probabilities\thanks{We would like to thank
Noureddine El Karoui, Adityanand Guntuboyina, Michael Jansson, Ulrich
M\"{u}ller, and three anonymous referees for useful comments and suggestions.
All errors are our own.}}
\author{Demian Pouzo,\thanks{ e-mail: dpouzo@berkeley.edu (corresponding author).
Address: UC Berkeley, Dept. of Economics; 530 Evans Hall, Berkeley CA 94619,
USA.} Zacharias Psaradakis\thanks{ e-mail: z.psaradakis@bbk.ac.uk. Address:
Birkbeck, University of London, Dept. of Economics, Mathematics \& Statistics;
Malet Street, London WC1E 7HX, UK} and Martin Sola\thanks{ e-mail:
msola@utdt.edu. Address: Universidad Torcuato Di Tella, Dept. of Economics;
Figueroa Alcorta 7350, C1428 Buenos Aires, Argentina.}}
\maketitle

\begin{abstract}
\noindent This paper considers maximum likelihood (ML) estimation in a large
class of models with hidden Markov regimes. We investigate consistency of the
ML estimator and local asymptotic normality for the models under general
conditions which allow for autoregressive dynamics in the observable process,
Markov regime sequences with covariate-dependent transition matrices, and
possible model misspecification. A Monte Carlo study examines the
finite-sample properties of the ML estimator in correctly specified and
misspecified models. An empirical application is also discussed.\bigskip

\bigskip\newpage

\end{abstract}

\section{Introduction}

Stochastic models with parameters that are subject to changes driven by an
unobservable Markov chain (the regime or state sequence) have attracted
considerable attention in many different areas, the influential work by
\cite{hamilton89} being a prominent example from the econometrics literature.
An important subclass of such models, the so-called hidden Markov models, in
which observations are conditionally independent given the regime sequence,
are also widely used in a variety of disciplines. A common assumption in these
models is that the unobservable Markov chain is temporally homogeneous.

In this paper, we focus on a larger class of models in which the hidden regime
process and the observation process (conditional on the regimes) are both
temporally inhomogeneous Markov chains. This is a useful generalization of
models with a time-invariant transition mechanism that has found numerous
applications, especially in economics and finance.\footnote{Examples include,
among many others, applications to the analysis of business-cycle fluctuations
(e.g., \cite{filardo94}, \cite{FilardoGordon98}, \cite{ravn99},
\cite{Simpson01}, \cite{rivas15}), interest rates and yields (e.g.,
\cite{Gray96}, \cite{BekaertH01}, \cite{ang02}, \cite{angbek02j},
\cite{pses21}), consumption growth (\cite{Whitelaw00}), currency crises (e.g.,
\cite{peria02}, \cite{mourat08}), and cryptocurrency returns (\cite{tan21}).}
Statistical inference in this class of models is predominantly likelihood
based, even though very little is known about the asymptotic properties of the
relevant inferential procedures. In a typical application, inference is
conducted on the implicit assumption that the maximum likelihood (ML)
estimator of unknown parameters has its familiar properties of consistency,
asymptotic normality, and asymptotic efficiency, and associated confidence
sets and hypotheses tests are constructed in the usual manner. It hardly needs
noting that, unless these properties of likelihood-based inferential
procedures known from regular parametric estimation problems hold for Markov
regime-switching models with covariate-dependent transition probabilities,
inferences drawn from them cannot be justified in any meaningful way and
should be interpreted very cautiously. Arguing, for instance, as is common in
applied work, that an economic variable is a useful leading indicator for
business-cycle phases because it appears to have a `statistically significant'
coefficient in the transition functions of a Markov-switching model for output
growth is problematic when little is known about the properties of the
relevant estimators and related tests.

The main contribution of this paper is to provide consistency and asymptotic
normality results for a large class of models that are relevant in
applications. Our approach allows for autoregressive dynamics in the
observable process, covariate-dependence in the transition functions of the
hidden regime process, and potential model misspecification. To the best of
our knowledge, the only asymptotic results available on ML estimation in
Markov regime-switching models with covariate-dependent transition
probabilities are those of \cite{aill13}, who investigate consistency of the
ML estimator in a correctly specified model (i.e., a model which contains the
data-generating process). Our results include both consistency of the ML
estimator and local asymptotic normality (LAN) for the model, from which
asymptotic normality of the ML estimator can be inferred. Unlike
\cite{aill13}, who allow for a general hidden state space, we require the
latter to be finite, but do not restrict the model to be correctly specified.
In doing so, we also extend some results of \cite{white82} for independent,
identically distributed (i.i.d.) data to the case of dependent observations
and for classes of parametric distributions associated with dynamic models
with hidden Markov regimes. As such stochastic specifications are typically
highly parametric, it is important to understand the properties of
likelihood-based inferential procedures in situations where the true
probability structure of the data does not necessarily lie within the
parametric family of distributions specified by the model. We show that the ML
estimator in our setting converges to the true parameter value if the model is
correctly specified and to a pseudo-true parameter set if the model is
misspecified.\footnote{We refer to estimators for the various models discussed
throughout as ML estimators even though they may be obtained from a
pseudo-likelihood based on a misspecified model.} We also show that the sample
log-likelihood satisfies the LAN property, establish an asymptotic linear
representation for the ML estimator, obtain the asymptotic distribution of the
estimator, and present results relating to consistent estimation of its
asymptotic covariance matrix. These are the most general results available for
Markov regime-switching models with autoregressive dynamics and
covariate-dependent transition probabilities.

In related earlier work, \cite{Mevel04} examine consistency and asymptotic
normality of the ML estimator in misspecified hidden Markov models with a
finite state space, while \cite{douc12} consider consistency under general
state spaces. \cite{bickel96}, \cite{bickel98}, \cite{Jensen99},
\cite{douc01}, and \cite{douc11} investigate consistency and/or asymptotic
normality in correctly specified hidden Markov models with regime sequences
defined on either a finite or a general state space. \cite{Francq98} and
\cite{Krishna98} consider consistency in correctly specified autoregressive
models with Markov regimes defined on a finite state space. \cite{douc04} and
\cite{Kasahara2019} investigate consistency and asymptotic normality in a
similar autoregressive setup, but allow the regime sequence to take values in
a space that is not necessarily countable. In all of these papers, the hidden
regime sequence is assumed to be a temporally homogeneous Markov chain.

In the sequel, we follow \cite{bickel98} and \cite{douc04} fairly closely in
terms of the technical tools and the arguments used to establish our results,
but our setup is more general in certain respects. Like \cite{bickel98}, we
consider models with a finite hidden state space, but allow for autoregressive
dynamics in the observation sequence, covariate-dependence in the transition
probabilities of the regime sequence, and potential model misspecification. In
\cite{douc04}, the hidden Markov chain is allowed to take values in a compact
topological space, but is restricted to be temporally homogeneous and the
model is assumed to be correctly specified. The cornerstone of the methods
used in these papers for establishing the asymptotic properties of the ML
estimator are mixing-type results for the unobservable regime sequence
conditional on the observation sequence (see also \cite{bickel96}). This is
also true for our approach, although the aforementioned results cannot be
invoked directly because they are established under the assumption of temporal
homogeneity of the hidden Markov chain. We, therefore, extend these results to
allow for more general Markov regime sequences; in particular, we establish
mixing-type results for the unobservable regime sequence given the observed
data, allowing for a particular form of covariate-dependence in the transition
kernels. This last result is, to our knowledge, novel and may be of interest
in its own right.

The remainder of the paper is organized as follows. Section~\ref{sec:models}
defines the class of models under consideration, gives sufficient conditions
for stationarity and ergodicity of the observation process, and describes the
estimation problem. Section~\ref{sec:consistent} investigates consistency of
the ML estimator in a general setting. Section~\ref{sec:anormal} contains
results on the LAN property of the model and the asymptotic normality of the
ML estimator. Section~\ref{sec:simulation} presents simulation results on the
finite-sample properties of estimators based on well-specified and
misspecified likelihoods. Section~\ref{sec:illustration} presents an
illustration using real-world data. Proofs of the main results are gathered in
an Appendix.
%Appendix \ref{sec:proofs} gathers the proofs.
%Section~\ref{sec:summary} summarizes and concludes. All proofs are gathered in an Appendix.

The following notational conventions are used throughout the paper. For an
infinite sequence $(V_{j})_{j}$, $V_{a}^{b}=(V_{a},\ldots,V_{b})$ for any
$a\leq b$; $\mathcal{P}(\mathbb{V})$ denotes the set of Borel probability
measures on a Polish space $\mathbb{V}$; for a probability measure $P$,
$E_{P}(\cdot)$ denotes expectation with respect to $P$, $o_{P}(\cdot)$ and
$O_{P}(\cdot)$ indicate order in probability under $P$, $\Rightarrow_{P}$
signifies weak convergence under $P$, and $L^{r}(P)$, $1\leq r<\infty$,
denotes the class of measurable functions integrable to order $r$ with respect
to $P$; $\nabla_{\vartheta}$ and $\nabla_{\vartheta}^{2}$ are the gradient and
Hessian operators, respectively, with respect to $\vartheta$; $\left\Vert
\cdot\right\Vert $ denotes the Euclidean norm of a vector or matrix;
$1\{\cdot\}$ denotes the indicator function; $\mathbb{N}$ denotes the set of
positive integers. Unless stated otherwise, limits are taken as the sample
size, $T$, diverges to infinity.

\section{Model and Estimation}

\label{sec:models}

\subsection{Statistical Model}

\label{sec:model}

Let $(X_{t},S_{t})_{t=0}^{\infty}$ be a discrete-time stochastic process such
that, for each $t\in\mathbb{N}$, $S_{t}\in\mathbb{S}\equiv\{s_{1}%
,\ldots,s_{|\mathbb{S}|}\}\subset\mathbb{R}$ is the unobservable state and
$X_{t}\in\mathbb{X}\subseteq\mathbb{R}^{h}$, for some $h\in\mathbb{N}$, is the
observable state. Moreover, for each $t\in\mathbb{N}$,~the conditional
distribution of $X_{t}$, given $X_{0}^{t-1}$ and $S_{0}^{t}$, depends only on
$X_{t-1}$ and $S_{t}$, and the conditional distribution of $S_{t}$, given
$X_{0}^{t-1}$ and $S_{0}^{t-1}$, depends only on $X_{t-1}$ and $S_{t-1}$, so
that
\begin{align*}
&  X_{t}\mid(X_{0}^{t-1},S_{0}^{t})\sim P_{\ast}(X_{t-1},S_{t},\cdot),\\
&  S_{t}\mid(X_{0}^{t-1},S_{0}^{t-1})\sim Q_{\ast}(X_{t-1},S_{t-1},\cdot),
\end{align*}
with $(x,s)\mapsto P_{\ast}(x,s,\cdot)\in\mathcal{P}(\mathbb{X})$ and
$(x,s)\mapsto Q_{\ast}(x,s,\cdot)\in\mathcal{P}(\mathbb{S})$ denoting the true
transition probabilities. It is further assumed that, for each $(x,s)\in
\mathbb{X}\times\mathbb{S}$, $P_{\ast}(x,s,\cdot)$ admits a density $p_{\ast
}(x,s,\cdot)$ with respect to some $\sigma$-finite measure on $\mathbb{X}$.
Our framework imposes no additional restrictions on this measure; for
instance, it can be the Lebesgue measure (i.e., allow for continuous $X_{t}$)
or a counting measure (i.e., allow for discrete $X_{t}$).

The researcher's model is given by a family of transition probabilities
$(x,s)\mapsto P_{\theta}(x,s,\cdot)\in\mathcal{P}(\mathbb{X})$ and
$(x,s)\mapsto Q_{\theta}(x,s,\cdot)\in\mathcal{P}(\mathbb{S})$ indexed by an
(unknown) parameter $\theta\in\Theta\subseteq\mathbb{R}^{q}$, for some
$q\in\mathbb{N}$, such that, for each $\theta\in\Theta$,
\[
X_{t}\mid(X_{0}^{t-1},S_{0}^{t})\sim P_{\theta}(X_{t-1},S_{t},\cdot),
\]%
\[
S_{t}\mid(X_{0}^{t-1},S_{0}^{t-1})\sim Q_{\theta}(X_{t-1},S_{t-1},\cdot),
\]
and, for each $(x,s)\in\mathbb{X}\times\mathbb{S}$, $P_{\theta}(x,s,\cdot)$
admits a density $p_{\theta}(x,s,\cdot)$ with respect to the same measure used
to define $p_{\ast}(x,s,\cdot)$.

A few remarks about this setup are worth making. First, and perhaps most
importantly, the unobservable states (regimes) are a Markov chain whose
transition kernel can depend on the lagged value of the observable state. This
is the main departure from prior literature which, with the exception of
\cite{aill13}, has focused on the case where the conditional distribution of
$S_{t}$, given $X_{0}^{t-1}$ and $S_{0}^{t-1}$, depends only on $S_{t-1}$.
Second, the model $\{(P_{\theta},Q_{\theta})\colon\theta\in\Theta\}$ is
allowed to be misspecified in the sense that $(P_{\ast},Q_{\ast}%
)\notin\{(P_{\theta},Q_{\theta})\colon\theta\in\Theta\}$. This setup
encompasses a rich family of models that arise in econometric and statistical
applications, some examples of which are given below. Note that, although the
family $\{Q_{\theta}\colon\theta\in\Theta\}$ may or may not contain $Q_{\ast}%
$, it is defined on the same finite state space $\mathbb{S}$ as $Q_{\ast}$, so
misspecification of the number of unobservable regimes is ruled
out.\footnote{Such misspecification would introduce additional complexities
(e.g., locally non-quadratic likelihood surfaces, non-identifiable parameters)
which are beyond the scope of this paper.} Finally, even though the
conditional distribution of $S_{t}$, given $X_{0}^{t-1}$ and $S_{0}^{t-1}$, is
assumed to depend only on $X_{t-1}$ and $S_{t-1}$, it is straightforward to
extend all our results to situations where this distribution is dependent on
additional higher-order lags of $X_{t}$ and/or $S_{t}$.

\begin{example}
[Hidden Markov Model with Covariate-Dependent Transition Probabilities]%
\label{exa:HMM} Let $x=(y,z)\in\mathbb{X}=\mathbb{R}^{2}$ and $\mathbb{S}%
=\{0,1\}$. Let $P_{\ast}$ be determined by the equations
\begin{align*}
Y_{t}  &  =\mu^{\ast}(S_{t})+\sigma^{\ast}(S_{t})U_{1,t},\\
Z_{t}  &  =\mu_{2}^{\ast}+\psi^{\ast}Z_{t-1}+\sigma_{2}^{\ast}U_{2,t},
\end{align*}
where $(U_{1,t},U_{2,t})_{t}$ are i.i.d. (independent of $(S_{t})_{t}$) with
zero mean and covariance matrix indexed by a parameter $\rho^{\ast}$, e.g.,
$\left[
\begin{array}
[c]{cc}%
1 & \rho^{\ast}\\
\rho^{\ast} & 1
\end{array}
\right]  $. The transition probabilities of $(S_{t})_{t}$ are allowed to
depend on $Z_{t-1}$; for instance, $(s,z)\mapsto Q_{\ast}(z,s,s)\equiv
\Pr(S_{t}=s\mid Z_{t-1}=z,S_{t-1}=s)=[1+\exp(-\alpha_{s}^{\ast}-\beta
_{s}^{\ast}z)]^{-1}$ for $s\in\mathbb{S}$. The homogeneous specification with
$\beta_{0}^{\ast}=\beta_{1}^{\ast}=0$ has been used to model regime shifts in
a variety of economic and financial time series, including output growth
(\cite{AlbertChib93}, \cite{rivas15}), foreign exchange rates
(\cite{engelhamilton90}, \cite{bollen0}), and equity returns (\cite{ryden98},
\cite{ang02}). The homogeneity restriction is relaxed in \cite{dieb94},
\cite{Engel96}, and \cite{ang02}, among others, to allow the transition
probabilities to depend on $Z_{t-1}$. The results obtained here establish the
asymptotic properties of the ML estimator in this class of models. $\triangle$
\end{example}

\begin{example}
[Markov-Switching Autoregressive Model with Covariate-Dependent Transition
Probabilities]\label{exa:MSR} A useful generalization of the previous example
is one where the outcome equation is extended to
\[
Y_{t}=\mu^{\ast}(S_{t})+\phi^{\ast}Y_{t-1}+\sigma^{\ast}(S_{t})U_{1,t}.
\]
Variations of the model with $\beta_{0}^{\ast}=\beta_{1}^{\ast}=0$ have found
widespread application in economics (e.g., \cite{Hansen92},
\cite{McCulloghTsay94}, \cite{Murcia95}, \cite{angbekwei08}) and beyond.
Generalizations of the model without the restriction of time-invariant
transition probabilities are also very popular and have been used, for
example, in the modeling of output growth (\cite{rivas15}), interest rates
(\cite{angbek02j}), consumption growth (\cite{Whitelaw00}), and bond spreads
(\cite{pses21}). The results obtained here establish the asymptotic properties
of the ML estimator in this class of models. $\triangle$
\end{example}

\begin{example}
[Mixture Autoregressive Model]\label{exa:MAM}Let $x\in\mathbb{X}=\mathbb{R}$,
$\mathbb{S}=\{0,1\}$, and for each $t\in\mathbb{N}$, $\Pr(S_{t}=0\mid
X_{t-1})=G_{\ast}(X_{t-1})$ for some $x\mapsto G_{\ast}(x)\in\lbrack0,1]$ and
$X_{t}\sim P_{\ast}(X_{t-1},S_{t},\cdot)$. This specification implies a
conditional density for $X_{t}$, given $X_{t-1}$, that is a mixture of the
type
\[
x\mapsto p_{\ast}(x\mid x_{t-1})=G_{\ast}(x_{t-1})p_{\ast}(x_{t-1}%
,0,x)+(1-G_{\ast}(x_{t-1}))p_{\ast}(x_{t-1},1,x).
\]
The covariate-dependence of the transition functions is reflected by the fact
that $G_{\ast}$ depends on $X_{t-1}$. Models which belong to the general class
of mixture autoregressive models (e.g., \cite{dss07}, \cite{Tadjuidje09},
\cite{dsps11}, \cite{saik15}) are covered by this framework. $\triangle$
\end{example}

The examples above, as well as those that follow, illustrate that in many
areas of application the stochastic process $(X_{t},S_{t})_{t=0}^{\infty}$ is
typically highly complex and it is natural/desirable to allow for feedback
from past realizations of the observable process $(X_{t})_{t=0}^{\infty}$ to
the law of the unobservable regime sequence $(S_{t})_{t=0}^{\infty}$; a
tractable way for modeling such feedback is to allow the transition kernel
$Q_{\theta}$ to depend on $X_{t-1}$. This feature adds an additional level of
complexity and with it sources of potential misspecification. For instance, a
common assumption in most applications is that the transition probabilities of
hidden regimes are time-invariant, an assumption which may result in
misspecification of the transitions functions. In applications that allow for
covariate-dependent transition functions, a somewhat more subtle and often
overlooked source of misspecification, associated with endogeneity of the
transition-driving covariates, may come into play. Specifically, in the
context of models such as those in Examples~\ref{exa:HMM} and \ref{exa:MSR},
contemporaneous correlation between $U_{1,t}$ and $U_{2,t}$ is typically
ignored and inference is based on the likelihood implied by the outcome
equation alone; we consider this case in more detail in Sections
\ref{sec:simulation2} and \ref{sec:illustration}.

Before discussing estimation of $\theta$, we give a result regarding the
mixing and ergodicity properties of $(X_{t})_{t=0}^{\infty}$. To do so, let
$\bar{P}_{\ast}^{\kappa}$ denote the true distribution over $(X_{t}%
)_{t=0}^{\infty}$ when the distribution of $(X_{0},S_{0})$ is $\kappa$. Under
the following assumptions, Lemma~\ref{lem:sta-ergo} below ensures that there
exists a Borel probability measure on $\mathbb{X}\times\mathbb{S}$, denoted
henceforth by $\nu$, for which $(X_{t})_{t=0}^{\infty}$ is stationary and ergodic.

\begin{assumption}
\label{ass:BDD_Q} There exists a continuous function $\underline{q}%
:\mathbb{X}\rightarrow\mathbb{R}_{+}\setminus\{0\}$ such that, for all
$Q\in\{Q_{\theta}\colon\theta\in\Theta\}\cup Q_{\ast}$, $Q(x,s,s^{\prime}%
)\geq\underline{q}(x)$ for all $(s^{\prime},s,x)\in\mathbb{S}^{2}%
\times\mathbb{X}$.
\end{assumption}

\begin{assumption}
\label{ass:fx-ergo} There exist constants $\lambda^{\prime}\in(0,1)$,
$\gamma\in(0,1)$, $b^{\prime}>0$ and $R>2b^{\prime}/(1-\gamma)$, a lower
semi-continuous function $\mathcal{U}:\mathbb{X}\rightarrow\lbrack1,\infty)$,
and a measure $\varpi\in\mathcal{P}(\mathbb{X})$ such that, for all
$s\in\mathbb{S}$: (i)$~\int_{\mathbb{X}}\mathcal{U}(x^{\prime})P_{\ast
}(x,s,dx^{\prime})\leq\gamma\mathcal{U}(x)+b^{\prime}1\{x\in A\}$, with
$A\equiv\{x\in\mathbb{X}\colon\mathcal{U}(x)\leq R\}$; (ii)$~A$ is bounded and
$\varpi(A)>0$; (iii)$~\inf_{x\in A}P_{\ast}(x,s,C)\geq\lambda^{\prime}%
\varpi(C)$ for any Borel set $C\subseteq\mathbb{X}$.
%	XXXXXX
%		There exists a lower semi-continuous function $\mathcal{U} : \mathbb{X} \rightarrow \mathbb{R}_{+}$  such that: (i) $\int_{\mathbb{X}}
%	\mathcal{U}(x^{\prime }) P_{\theta^{\ast}} (dx^{\prime }|x,s) \leq \gamma
%	\mathcal{U}(x) + K$ for some $\gamma \in [0,1)$ and $K \geq 0$; (ii)  $\{ x \in \mathbb{X}
%	: \mathcal{U}(x) \leq R\}$ is bounded for all $R>0$; (iii) there exists a $%
%	\lambda >0$ and a  measure $\varpi$ such that $P_{\theta^{\ast}} (C \mid x,s) \geq  \lambda \varpi(C)$ for any $C \subseteq \mathbb{X}$ Borel and any $x
%	\in \{ x \in \mathbb{X} :  \mathcal{U}(x) \leq R \}$ some $R > 2K/(1-\gamma)$, and any $s \in \mathbb{S}$.

\end{assumption}

The following lemma establishes stationarity, ergodicity, and $\beta$-mixing
of $(X_{t})_{t=0}^{\infty}$.

\begin{lemma}
\label{lem:sta-ergo} Suppose Assumptions \ref{ass:BDD_Q} and \ref{ass:fx-ergo}
hold. Then, there exists a $\nu\in\mathcal{P}(\mathbb{X}\times\mathbb{S})$
such that, under $\bar{P}_{\ast}^{\nu}$, $(X_{t})_{t=0}^{\infty}$ is
stationary, ergodic, and $\beta$-mixing with mixing coefficients $\beta
_{n}=O(\gamma^{n})$, $n\in\mathbb{N}$.
\end{lemma}

\begin{proof}
See Supplemental Material \ref{SM:Ergo}.
\end{proof}

The result follows in a standard manner by using Assumptions \ref{ass:BDD_Q}
and \ref{ass:fx-ergo} to establish that the implied transition kernel of the
joint process $(X_{t},S_{t})_{t=0}^{\infty}$ has a unique invariant
distribution and also that it is Harris recurrent and aperiodic. This fact, in
turn, is used to show that $(X_{t})_{t=0}^{\infty}$ is stationary, ergodic,
and $\beta$-mixing at a geometric rate.

\medspace

\begin{remark}
[Discussion of Assumptions \ref{ass:BDD_Q} and \ref{ass:fx-ergo}%
]\label{rem:BDD-Q} Assumption \ref{ass:BDD_Q} is an extension of a common
assumption in the literature (cf. \cite{douc04}, \cite{aill13}) to the case
where the transition kernel of $(S_{t})_{t=0}^{\infty}$ depends on $X_{t-1}$.
Allowing the lower bound $\underline{q}$ to depend on $x$ is especially
relevant when the support of $X_{t}$ is unbounded because, while
$\underline{q}(x)>0$, it is allowed to converge to zero as $\left\Vert
x\right\Vert \rightarrow\infty$. Although this assumption is not innocuous, we
view it as mild because it accommodates the typical specifications used in the
literature, where $Q$ is parameterized by a standard Gaussian or logistic
cumulative distribution function and a single index $x^{\intercal}\beta$, with
$\beta$ restricted to take values in a bounded subset of a finite-dimensional
Euclidean space.

Assumption \ref{ass:fx-ergo}(iii) is an analogous condition for the transition
kernel $P_{\ast}$. By inspection of the proof of Lemma~\ref{lem:sta-ergo}, it
is easy to see that it suffices to obtain a minorization condition for the
\textquotedblleft joint\textquotedblright\ kernel, i.e., $\inf_{x\in A}%
P_{\ast}(x,s^{\prime},C)Q(x,s,s^{\prime})\geq\lambda\widetilde{\varpi
}(C,s^{\prime})$ for any Borel set $C\subseteq\mathbb{X}$ and for some
$\widetilde{\varpi}\in\mathcal{P}(\mathbb{X}\times\mathbb{S})$ and $\lambda
\in(0,1)$. Thus, Assumptions \ref{ass:BDD_Q}(i) and \ref{ass:fx-ergo}(i) could
be relaxed; e.g., the former could be relaxed to $Q(x,s,s^{\prime}%
)\geq\underline{q}(x)\varrho(s^{\prime})$, where $\varrho\in\mathcal{P}%
(\mathbb{S})$, or the latter could be relaxed to $\inf_{x\in A}P_{\ast
}(x,s^{\prime},C)\geq\lambda^{\prime}\widetilde{\varpi}(C,s^{\prime})$, where
$\widetilde{\varpi}\in\mathcal{P}(\mathbb{X}\times\mathbb{S})$.

Assumption \ref{ass:fx-ergo}(i),(ii) is a so-called Foster--Lyapunov drift
condition; see \cite{meyn2012markov} and references therein for a discussion
of the assumption. $\triangle$
\end{remark}

\medspace

In view of Lemma~\ref{lem:sta-ergo}, under $\nu$, the process $(X_{t}%
)_{t=0}^{\infty}$ can be extended to a two-sided sequence $(X_{t})_{t=-\infty
}^{\infty}$. With a slight abuse of notation, we still use $\bar{P}_{\ast
}^{\nu}$ to denote the true probability distribution over $(X_{t}%
,S_{t})_{t=-\infty}^{\infty}$; $\bar{P}_{\theta}^{\nu}$ is defined analogously
for the model $(Q_{\theta},p_{\theta},\nu)$.\footnote{Throughout the text, we
use $\bar{P}_{\theta}^{\nu}$ to denote any marginal or conditional
probabilities associated with $\bar{P}_{\theta}^{\nu}$; the same holds for
$\bar{P}_{\ast}^{\nu}$.}

\subsection{Parameter Estimation}

\label{sec:estimation}

For any $T\in\mathbb{N}$, let $\ell_{T}^{\nu}:\mathbb{X}^{T+1}\times
\Theta\rightarrow\mathbb{R}$ be the sample criterion function given by
%\begin{equation*}
%\ell _{T}^{\nu }(X_{0}^{T},\theta )=T^{-1}\log p_{T}^{\nu }(X_{0}^{T},\theta
%)=T^{-1}\sum_{t=1}^{T}\log p_{t}^{\nu }(X_{t}\mid X_{0}^{t-1},\theta ),
%\end{equation*}%
\begin{equation}
\ell_{T}^{\nu}(X_{0}^{T},\theta)=T^{-1}\sum_{t=1}^{T}\log p_{t}^{\nu}%
(X_{t}\mid X_{0}^{t-1},\theta), \label{logl}%
\end{equation}
where $p_{t}^{\nu}(X_{t}\mid X_{0}^{t-1},\theta)$ denotes the conditional
density of $X_{t}$ given $X_{0}^{t-1}$ for any $\theta\in\Theta$; the latter
is defined recursively as follows: for any $t\geq1,$
\[
p_{t}^{\nu}(X_{t}\mid X_{0}^{t-1},\theta)=\sum_{s^{\prime}\in\mathbb{S}}%
\sum_{s\in\mathbb{S}}p_{\theta}(X_{t-1},s^{\prime},X_{t})Q_{\theta}%
(X_{t-1},s,s^{\prime})\delta_{t}^{\theta,\nu}(s),
\]
and $s\mapsto\delta_{t}^{\theta,\nu}(s)\equiv\bar{P}_{\theta}^{\nu}%
(S_{t-1}=s\mid X_{0}^{t-1})$. For each $t\geq2$ and any $s\in\mathbb{S}$,
$s\mapsto\delta_{t}^{\theta,\nu}(s)$ satisfies the recursion
%\begin{align}
%\delta _{t}^{\theta ,\nu }(s)=& \sum_{\tilde{s}\in \mathbb{S}}\Pr (s\mid
%\tilde{s},X_{0}^{t-1})\Pr (\tilde{s}\mid X_{0}^{t-1})  \notag \\
%=& \sum_{\tilde{s}\in \mathbb{S}}\frac{\Pr (s\mid \tilde{s},X_{0}^{t-1})\Pr (%
%\tilde{s},X_{t-1},X_{0}^{t-2})}{\sum_{s^{\prime }\in \mathbb{S}}p_{\theta
%}(X_{t-1}\mid X_{t-2},s^{\prime })\delta _{t-1}^{\theta ,\nu }(s^{\prime })}
%\notag \\
%=& \sum_{\tilde{s}\in \mathbb{S}}\frac{Q_{\theta }(s\mid \tilde{s}%
%,X_{t-1})p_{\theta}(X_{t-1}\mid X_{t-2},\tilde{s})\delta _{t-1}^{\theta
%,\nu }(\tilde{s})}{\sum_{s^{\prime }\in \mathbb{S}}p_{\theta}(X_{t-1}\mid
%X_{t-2},s^{\prime })\delta _{t-1}^{\theta ,\nu }(s^{\prime })},
%\label{eqn:funcf}
%\end{align}%
\[
\delta_{t}^{\theta,\nu}(s)=\sum_{\tilde{s}\in\mathbb{S}}\frac{Q_{\theta
}(X_{t-1},\tilde{s},s)p_{\theta}(X_{t-2},\tilde{s},X_{t-1})\delta
_{t-1}^{\theta,\nu}(\tilde{s})}{\sum_{s^{\prime}\in\mathbb{S}}p_{\theta
}(X_{t-2},s^{\prime},X_{t-1})\delta_{t-1}^{\theta,\nu}(s^{\prime})},
\]
with $s\mapsto\delta_{1}^{\theta,\nu}(s)=\sum_{\tilde{s}\in\mathbb{S}%
}Q_{\theta}(X_{0},\tilde{s},s)\nu(\tilde{s}|X_{0})$, where $\nu(\cdot|\cdot)$
is the conditional density corresponding to $\nu$.

For a given initial distribution $\kappa\in\mathcal{P}(\mathbb{X}%
\times\mathbb{S})$ over $(X_{0},S_{0})$, we define our estimator as
$\hat{\theta}_{\kappa,T}$, where
\begin{equation}
\ell_{T}^{\kappa}(X_{0}^{T},\hat{\theta}_{\kappa,T})\geq\sup_{\theta\in\Theta
}\ell_{T}^{\kappa}(X_{0}^{T},\theta)-\eta_{T}, \label{mle}%
\end{equation}
for some $\eta_{T}\geq0$ and $\eta_{T}=o(1)$.

\section{Consistency}

\label{sec:consistent}

The main result of this section establishes convergence of the estimator
$\hat{\theta}_{\nu,T}$ to the set of points in $\Theta$ that are closest to
the true model under the Kullback--Leibler information criterion.

Let $H^{\ast}:\Theta\rightarrow\mathbb{R}_{+}\cup\{\infty\}$ be the
Kullback--Leibler information criterion $\theta\mapsto H^{\ast}(\theta)$,
which is given by
\[
H^{\ast}(\theta)=E_{\bar{P}_{\ast}^{\nu}}\left[  \log\frac{p_{\ast}^{\nu
}(X_{0}\mid X_{-\infty}^{-1})}{p^{\nu}(X_{0}\mid X_{-\infty}^{-1},\theta
)}\right]  ,
\]
where, for any $\theta\in\Theta$, $p^{\nu}(X_{0}\mid X_{-\infty}^{-1},\theta)$
denotes the conditional density of $X_{0}$ given $X_{-\infty}^{-1}$ induced by
$(P_{\theta},Q_{\theta},\nu)$, and $p_{\ast}^{\nu}(X_{0}\mid X_{-\infty}%
^{-1})$ is its counterpart induced by the true transition kernels $(P_{\ast
},Q_{\ast},\nu)$; we refer the reader to the Supplemental Material
\ref{SM:PDFs} for details about the construction of these objects and their properties.

Our results allow for misspecified models, and thus $p_{\ast}^{\nu}%
\notin\{p^{\nu}(\cdot\mid\cdot,\theta):\theta\in\Theta\}$. Hence, as in
\cite{white82}, the relevant limiting set for our estimator is
\[
\Theta_{\ast}=\underset{\theta\in\Theta}{\arg\min}H^{\ast}(\theta),
\]
which is the \emph{pseudo-true parameter (set)} that minimizes the
Kullback--Leibler information criterion. Under Assumption \ref{ass:exist}
below, $\Theta_{\ast}$ is non-empty and compact by the Weierstrass Theorem.

\begin{assumption}
\label{ass:exist} (i) $\Theta$ is compact; (ii) $H^{\ast}$ exists and is lower semi-continuous.
\end{assumption}

The following additional assumptions are used to establish the main result. To
state these, for any $\delta>0$ and any $\dot{\theta}\in\Theta$, let
$B(\delta,\dot{\theta})\equiv\{\theta\in\Theta\colon||\dot{\theta}%
-\theta||<\delta\}$.

\begin{assumption}
\label{ass:bdd} (i) For any $\epsilon>0$, there exists some $\delta>0$ such
that
\[
\max_{\dot{\theta}\in\Theta}E_{\bar{P}_{\ast}^{\nu}}\left[  \sup_{\theta\in
B(\delta,\dot{\theta})}\frac{p^{\nu}(X_{0}\mid X_{-\infty}^{-1},\theta
)}{p^{\nu}(X_{0}\mid X_{-\infty}^{-1},\dot{\theta})}\right]  \leq1+\epsilon;
\]
(ii) there exists a function $(x,x^{\prime})\mapsto C(x,x^{\prime}%
)\in\mathbb{R}_{+}$ such that $\sup_{\theta\in\Theta}\frac{\max_{s\in
\mathbb{S}}p_{\theta}(X,s,X^{\prime})}{\min_{s\in\mathbb{S}}p_{\theta
}(X,s,X^{\prime})}\leq C(X,X^{\prime})$ and $\frac{\max_{s\in\mathbb{S}%
}p_{\ast}(X,s,X^{\prime})}{\min_{s\in\mathbb{S}}p_{\ast}(X,s,X^{\prime})}\leq
C(X,X^{\prime})$ a.s.-$\bar{P}_{\ast}^{\nu}$.
\end{assumption}

\begin{assumption}
\label{ass:qbar-sum} $T^{-1}\sum_{t=1}^{T}\max\{1,C(X_{t-1},X_{t}%
)\}\prod_{i=0}^{t-1}(1-\underline{q}(X_{i}))=o_{\bar{P}_{\nu}^{\ast}}(1).$
\end{assumption}

\begin{remark}
[Discussion of Assumptions \ref{ass:exist}, \ref{ass:bdd} and
\ref{ass:qbar-sum}]Assumption~\ref{ass:exist}(i) is standard.
Assumption~\ref{ass:exist}(ii) is high-level but can be obtained from
lower-level conditions (e.g., by following the reasoning in Proposition~11 of
\cite{douc12}).

Assumption \ref{ass:bdd}(i) is a high-level condition used for establishing
uniform law of large numbers results (see Lemma~\ref{lem:lik-conv} in Section
\ref{app:consistent}). Assumption \ref{ass:bdd}(ii) is akin to Assumption~A4
in \cite{bickel98}; it essentially restricts the support of $p_{\theta}$ and
$p_{\ast}$ for different values of the state variable.\footnote{The second
part of Assumption \ref{ass:bdd}(ii) is used to show that $p_{\ast}^{\nu
}(\cdot|X_{-\infty}^{-1})$ integrates to one.}

Assumption \ref{ass:qbar-sum} essentially requires that $\underline{q}$ is not
\textquotedblleft too close\textquotedblright\ to zero on average and that
$C(X_{t-1},X_{t})$ is finite a.s.-$\bar{P}_{\nu}^{\ast}$. For instance, if
$\underline{q}(x)\geq c$ for some $c>0$, then Assumption \ref{ass:qbar-sum} is
automatically satisfied provided $E_{\bar{P}_{\nu}^{\ast}}[C(X_{0}%
,X_{1})]<\infty$. Moreover, by exploiting the fact that $(X_{t})_{t}$ is
$\beta$-mixing (see Lemma \ref{lem:sta-ergo}), Lemma \ref{lem:suff-A5andA8} in
the Supplemental Material \ref{SM:SuffAssumptions} provides sufficient
conditions of the form $E_{\bar{P}_{\nu}^{\ast}}\left[  \underline{q}%
(X_{1})\right]  >0$ and $E_{\bar{P}_{\nu}^{\ast}}\left[  C(X_{1},X_{0}%
)^{l}\right]  <\infty$ for some $l>1$.\footnote{We thank a referee and the
editor for suggestions on how to weaken Assumption \ref{ass:bdd}, which lead
to these sufficient conditions for Assumption \ref{ass:qbar-sum}.}
%
%XXXXXXXXXXXXX
%A low level sufficient condition is given by: For any $\epsilon>0$, there exists a $\delta>0$ such that
%\begin{align*}
%	\sup_{||\theta_{1} - \theta_{2}|| \leq \delta}  \max_{s \in\mathbb{S}} |f_{\theta_{1}}(X_{1} \mid X_{0} ,s )/f_{\theta_{2}}(X_{1} \mid X_{0} ,s ) - 1 | \leq C(X_{1},X_{0})\epsilon,~with~E[C(X_{1},X_{0})] < \infty,
%\end{align*}
%a.s.-$\bar{P}^{\nu}_{\theta^{\ast}}$. Because, since $p^{\nu}(X_{0} \mid X^{-1}_{-\infty}, \theta) = \sum_{s \in \mathbb{S}} f_{\theta}(X_{1} \mid X_{0} ,s ) Pr_{\theta}(s \mid X^{-1}_{-\infty})$, under this condition and Assumption \ref{ass:bdd}(ii), it follows that for any $\epsilon^{\prime }>0$ and choosing $\delta>0$ as in the condition,
%\begin{align*}
%	\sup_{\theta \in ||\theta ^{\prime }-\theta
%		||<\delta }\frac{p^{\nu }(X_{0}\mid X_{-\infty }^{-1},\theta^{\prime} )}{p^{\nu
%		}(X_{0}\mid X_{-\infty }^{-1},\theta)} \leq & \sup_{\theta \in ||\theta ^{\prime }-\theta
%		||<\delta }\frac{\sum_{s \in \mathbb{S}} f_{\theta}(X_{1} \mid X_{0} ,s )(1 + C(X_{1},X_{0})\epsilon^{\prime }) Pr_{\theta^{\prime}}(s \mid X^{-1}_{-\infty})   }{\sum_{s \in \mathbb{S}} f_{\theta}(X_{1} \mid X_{0} ,s ) Pr_{\theta}(s \mid X^{-1}_{-\infty}) } \\
%	\leq&(1 + C(X_{1},X_{0})\epsilon^{\prime }) C  \frac{\sum_{s \in \mathbb{S}} Pr_{\theta^{\prime}}(s \mid X^{-1}_{-\infty})   }{\sum_{s \in \mathbb{S}} Pr_{\theta}(s \mid X^{-1}_{-\infty}) }\\
%	   = & C(1 + C(X_{1},X_{0})\epsilon^{\prime }).
%\end{align*}
%So be choosing $\delta$ so that $\epsilon = E[]C(X_{1},X_{0})] \epsilon'$, Assumption \ref{ass:bdd}(i) follows.
$\triangle$
\end{remark}

We now establish consistency of the estimator defined by (\ref{mle}).

\begin{theorem}
\label{thm:consistent} Suppose Assumptions \ref{ass:BDD_Q}--\ref{ass:bdd}
hold. Then, $d_{\Theta}(\hat{\theta}_{\nu,T},\Theta_{\ast})=o_{\bar{P}_{\ast
}^{\nu}}(1).$\footnote{For any set $A\subseteq\Theta$, $d_{\Theta}%
(\theta,A)\equiv\inf_{\dot{\theta}\in A}||\theta-\dot{\theta}||$.}
\end{theorem}

\begin{proof}
See Appendix \ref{app:consistent}.
\end{proof}

This result is analogous to Theorem~2 in \cite{douc12} but for a somewhat
different setup; specifically, we allow for autoregressive dynamics and
covariate-dependent transition probabilities, but restrict $\mathbb{S}$ to be
finite.\footnote{Finiteness of $\mathbb{S}$ simplifies the proofs as we do not
have to be concerned with uniformity issues of certain quantities as functions
of the hidden state. By combining techniques available in the literature
(e.g., \cite{douc12}) -- which essentially amount to imposing requirements
like compactness of $\mathbb{S}$ and continuity -- with ours, we conjecture
that our results can be extended to the case where $\mathbb{S}$ is a compact
set.}

Clearly, if the model is correctly specified and point identified, i.e., there
exists a $\theta\in\Theta$ such that $(P_{\ast},Q_{\ast})=(P_{\theta
},Q_{\theta})$, then $\Theta_{\ast}=\{\theta\}$ and our estimator converges in
$\bar{P}_{\ast}^{\nu}$-probability to this point. If, however, the model is
misspecified, our estimator converges to
%a pseudo-true parameter (set), which is
the set of parameters that is closest to the true set, when closeness is
measured by means of the Kullback--Leibler information criterion (cf.
\cite{white82}, \cite{douc12}).

To prove Theorem \ref{thm:consistent}, we first show that $T^{-1}\sum
_{t=1}^{T}\log p_{t}^{\nu}(X_{t}\mid X_{0}^{t-1},\theta)$ is well-approximated
by $T^{-1}\sum_{t=1}^{T}\log p^{\nu}(X_{t}\mid X_{-\infty}^{t-1},\theta)$ (see
Lemma~\ref{lem:lik-approx} in Appendix~\ref{app:consistent}). Second, relying
on ergodicity (Lemma~\ref{lem:sta-ergo}) and Assumption \ref{ass:bdd}, we
establish a \emph{uniform} law of large numbers for the latter quantity (see
Lemma~\ref{lem:lik-conv} in Appendix~\ref{app:consistent}). The proof of
consistency then follows the standard Wald approach.

The approximation result in the first step relies on \textquotedblleft
mixing\textquotedblright\ results for the process $(S_{t})_{t=-\infty}%
^{\infty}$, given $(X_{t})_{t=-\infty}^{\infty}$. The following theorem, which
might be of independent interest, establishes such a \textquotedblleft
mixing\textquotedblright\ result in our setting.

\begin{theorem}
\label{thm:Q-ergo} Take any $(j,m)\in\mathbb{N}^{2}$. Suppose that, for any
$\theta\in\Theta$, there exist mappings $x\mapsto\varrho(x,\cdot
)\in\mathcal{P}(\mathbb{S})$ and $\underline{q}:\mathbb{X}\rightarrow
\mathbb{R}_{+}$ such that, for all $(s,s^{\prime})\in\mathbb{S}^{2}$,
\begin{equation}
Q_{\theta}(X,s,s^{\prime})\geq\underline{q}(X)\varrho(X,s^{\prime}%
)\quad~a.s.\text{\textrm{-}}\bar{P}_{\ast}^{\nu}. \label{eqn:Q-ergo}%
\end{equation}
Then,\footnote{For any $P,Q\in$ $\mathcal{P}(\mathbb{S})$, $||P-Q||_{1}%
\equiv\sum_{s\in\mathbb{S}}|P(s)-Q(s)|$.}
\[
\max_{(b,c)\in\mathbb{S}^{2}}\left\Vert \bar{P}_{\theta}^{\nu}(S_{j+1}%
=\cdot|S_{-m}=b,X_{-m}^{j})-\bar{P}_{\theta}^{\nu}(S_{j+1}=\cdot
|S_{-m}=c,X_{-m}^{j})\right\Vert _{1}\leq\prod_{n=-m}^{j}(1-\underline{q}%
(X_{n})).
\]

\end{theorem}

\begin{proof}
See Appendix \ref{app:mixing}.
\end{proof}

This result is analogous to results in \cite{douc04} (e.g., their Lemma 1 and
Corollary 1) but for a more general transition function $Q_{\theta}$ (albeit
under the requirement that $\mathbb{S}$ is finite). Specifically, we allow
$Q_{\theta}$ to depend on $X$, and do not restrict the lower bound in
(\ref{eqn:Q-ergo}) to be uniform in $s^{\prime}$; these features are, to our
knowledge, novel.

%The proof relies on bounding the Dobrushin coefficient of the transition
%kernel $\bar{P}_{\theta }^{\nu }(S_{l+1}=\cdot |S_{l}=\cdot ,X_{-m}^{j})$ by
%$1-\underline{q}(X_{l})$ for each $l\in \{-m,...,j\}$. If $\varrho (x,\cdot
%) $, in condition (\ref{eqn:Q-ergo}), is uniformly bounded from below (e.g.,
%this will be the case under Assumption \ref{ass:BDD_Q}), then such bound is
%obtained by elementary calculations. For the general case, however, we had
%to use a different approach based on coupling techniques; see Section \ref%
%{app:mixing} for details.
The proof of Theorem~\ref{thm:Q-ergo} relies on bounding the Dobrushin
coefficient of the transition kernel $\bar{P}_{\theta}^{\nu}(S_{l+1}%
=\cdot|S_{l}=\cdot,X_{-m}^{j})$ by $1-\underline{q}(X_{l})$, for each
$l\in\{-m,...,j\}$. The case where $\varrho(x,\cdot)$, in condition
(\ref{eqn:Q-ergo}), is uniformly bounded from below has been studied in the
literature, and such a bound can be obtained by elementary calculations. For
the general case, however, the standard technique cannot be used, and we
develop a novel, to our knowledge, approach based on coupling techniques; see
Appendix~\ref{app:mixing} for details.

%\medspace

\begin{remark}
Remark \ref{rem:BDD-Q} and the fact that Theorem \ref{thm:Q-ergo} is
established under condition (\ref{eqn:Q-ergo}), imply that, in regards to
consistency, Assumption \ref{ass:BDD_Q} could be replaced by the weaker
condition (\ref{eqn:Q-ergo}). This remark, however, does not extend to the LAN
results of Section~\ref{sec:anormal}, since we do not know whether Assumption
\ref{ass:BDD_Q} could be weakened to condition (\ref{eqn:Q-ergo}) in this
case. $\triangle$
\end{remark}

%XXX PUT CONT OF EXAMPLES HERE. EXAMPLES (1) CURRENT EXAMPLE OF MISSPECIFIED
%JOINT LIKELIHOOD. NOT MUCH TO SAY. (2) EXAMPLE OF MISSPEFICIEF MIXTURE.
%XXX GOAL OF THESE EXAMPLES: (1) VERIFY ASSUMPTIONS. (2) PRESENT
%MISSPEFICIFIED CASES XXX

We discuss next a canonical example that encompasses many commonly used
specifications, such as those in Examples \ref{exa:HMM} and \ref{exa:MSR}. The
purpose of this example is to verify our assumptions under primitive and
low-level conditions.

\begin{example}
[Canonical Example]\label{exa:Canon} We verify the regularity conditions for a
model with $\mathbb{S}=\{0,1\}$ and
\begin{align*}
&  X_{t}=\mu(S_{t})+\Phi^{\intercal}X_{t-1}+\Sigma^{1/2}(S_{t})\varepsilon
_{t},\\
&  S_{t}\sim Q_{\bar{\vartheta}}(X_{t-1},S_{t-1},\cdot),
\end{align*}
where $(\varepsilon_{t})_{t}\sim\mathrm{i.i.d.}~\mathcal{N}(0,I)$, $I$ being
the identity matrix, and $x\mapsto Q_{\bar{\vartheta}}(x,s,s)\equiv\Pr
(S_{t}=s\mid X_{t-1}=x,S_{t-1}=s)=\Psi(x^{\intercal}\bar{\vartheta}_{s})$ for
$s\in\mathbb{S}$, where $\Psi$ is a full-support, continuous cumulative
distribution function (e.g., logistic or normal). It is assumed that the
parameter set $\Theta$ is compact and that any $\theta\equiv((\mu
(s),\Sigma(s))_{s\in\mathbb{S}},\Phi,\bar{\vartheta})$ in it is such that
$\Sigma(\cdot)$ and $\Phi\Phi^{\intercal}$ have eigenvalues uniformly bounded
away from zero and infinity. For simplicity, we assume that the true model is
indexed by $((\mu_{\ast}(s),\Sigma_{\ast}(s))_{s\in\mathbb{S}},\Phi_{\ast
},\bar{\vartheta}_{\ast})$, for which analogous conditions hold.

%the same assumptions hold for the parameter $((\mu_{\ast}(s),\Sigma_{\ast}(s))_{s\in\mathbb{S}},\Phi_{\ast},\beta_{\ast})$ that indexes the true model. %The noise $(\varepsilon_{t})_{t}$ is assumed to be independent of $(S_{t})_{t}$.

Let $q(x)\equiv\inf_{b}\Psi(x^{\intercal}b)$ and note that $q(x)>0$ for each
finite $x$ as $\Psi$ has full support.
%We also assume that $q$ is a lower bound for the primitive $Q_{\ast}$.
Thus, Assumption \ref{ass:BDD_Q} holds and $E_{\bar{P}_{\nu}^{\ast}}[q(X)]>0$.
The conditions of Assumption \ref{ass:fx-ergo} follow by the results in
\cite{douc2004practical}; see Lemma \ref{lem:fx-ergo-holds} in the
Supplemental Material \ref{SM:examples} for the formal argument. In Lemma
\ref{lem:exa-CBound} in the Supplemental Material \ref{SM:examples}, it is
shown that
\[
C^{-1}\underline{p}(x,y)\leq f_{\mathcal{N}}(\{y-\Phi^{\intercal}%
x-\mu(s)\}\Sigma^{-1/2}(s))\leq C\overline{p}(x,y),
\]
for some $C\geq1$, where $f_{\mathcal{N}}$ is the $\mathcal{N}(0,I)$
probability density function,
\begin{align*}
\underline{p}(x,y)\equiv &  \exp\{-0.5(a_{1}\left\Vert y\right\Vert ^{2}%
+a_{2}\left\Vert x\right\Vert ^{2}+2a_{3}\left\Vert x\right\Vert \left\Vert
y\right\Vert )-a_{4}\left\Vert y\right\Vert -a_{5}\left\Vert x\right\Vert
\},\\
\overline{p}(x,y)\equiv &  \exp\{-0.5(b_{1}\left\Vert y\right\Vert ^{2}%
+b_{2}\left\Vert x\right\Vert ^{2}-2b_{3}\left\Vert x\right\Vert \left\Vert
y\right\Vert )+b_{4}\left\Vert y\right\Vert +b_{5}\left\Vert x\right\Vert \},
\end{align*}
and $a_{i},b_{i}$ $(i=1,\ldots,5)$ are positive constants that are functions
of uniform bounds of $\Sigma(\cdot)$, $\mu(\cdot)$ and $\Phi$ over $\Theta$;
they are defined in Lemmas \ref{lem:exa-CBound} and \ref{lem:exa-CBound2} in
the Supplemental Material \ref{SM:SuffAssumptions}. Therefore, by
Lemma~\ref{lem:prop-ergoPDF} in the Supplemental Material \ref{SM:PDFs}, it
follows that Assumptions \ref{ass:exist} and \ref{ass:bdd} hold with
$(x,x^{\prime})\mapsto C(x,x^{\prime})\equiv\overline{p}(x,x^{\prime
})/\underline{p}(x,x^{\prime})$, provided that $E_{\bar{P}_{\nu}^{\ast}}%
[\exp\{-0.5l((b_{1}-a_{1})\left\Vert Y\right\Vert ^{2}+(b_{2}-a_{2})\left\Vert
X\right\Vert ^{2}-2(b_{3}+a_{3})\left\Vert X\right\Vert \left\Vert
Y\right\Vert )+l(b_{4}+a_{4})\left\Vert y\right\Vert +l(b_{5}+a_{5})\left\Vert
X\right\Vert \}]<\infty$ for some $l\geq1$. Lemma \ref{lem:exa-CBound2} in the
Supplemental Material \ref{SM:examples} shows that the latter condition holds
under some restrictions on the eigenvalues of $\Sigma(\cdot)$, $\Sigma_{\ast
}(\cdot)$, $\Phi\Phi^{\intercal}$, and $\Phi_{\ast}\Phi_{\ast}^{\intercal}$
(see the aforementioned lemma and its subsequent remark for the exact
formulation and more details). Finally, under these conditions, Lemma
\ref{lem:suff-A5andA8} in the Supplemental Material \ref{SM:SuffAssumptions}
shows that Assumption \ref{ass:qbar-sum} holds. $\triangle$
\end{example}

Lastly, we note that further characterization and interpretation of the
pseudo-true parameter set $\Theta_{\ast}$ is not straightforward in our
general setup, and we believe that progress in this direction should be made
on a case-by-case basis. While pursuing this program is outside the scope of
the present paper, in Section \ref{sec:simulation} we present numerical
results that attempt to shed some light on the characterization of this set.
In addition, Example~\ref{exa:HMM-miss} in the Supplemental Material
\ref{SM:examplemixture} provides analytical results on the characterization of
$\Theta_{\ast}$ for a subclass of the models in Example \ref{exa:MSR}, namely
hidden Markov models with covariate-dependent transition probabilities, when
the (misspecified) model is a simple mixture.

\section{Asymptotic Distribution Theory}

\label{sec:anormal}

In this section, we establish a LAN property (\cite[Ch. II]{IbragHas81},
\cite{leCam2012}) for our model and an asymptotic linear representation for
our estimator, from which asymptotic normality of the estimator is inferred.
For this, we require the following assumptions.

%Our main result extends results in \cite%
%{bickel96} and \cite{bickel98} to the case where the regime process is a
%time-inhomogeneous Markov chain and the observation process exhibits
%autoregressive dynamics (conditionally on the regimes). XXX ALSO misspecified no? XXX

\begin{assumption}
\label{ass:Theta-int} (i) $\Theta_{\ast}=\{\theta_{\ast}\}\subset
\mathrm{int}(\Theta)$; (ii) $\theta\mapsto p_{\theta}(X,s,X^{\prime})$ and
$\theta\mapsto Q_{\theta}(X,s,s^{\prime})$ are twice continuously
differentiable $a.s.$-$\bar{P}_{\ast}^{\nu}$ for all $(s,s^{\prime}%
)\in\mathbb{S}^{2}$.
\end{assumption}

\begin{assumption}
\label{ass:deriva-bdd} For some $\delta>0$ and $a\geq1$, and for all
$(s^{\prime},s)\in\mathbb{S}^{2}$: (i)%
\[
E_{\bar{P}_{\ast}^{\nu}}\left[  \sup_{\theta\in B(\delta,\theta_{\ast}%
)}\left\Vert \nabla_{\theta}\log p_{\theta}(X,s,X^{\prime})\right\Vert
^{2a}\right]  <\infty\text{ and}~E_{\bar{P}_{\ast}^{\nu}}\left[  \sup
_{\theta\in B(\delta,\theta_{\ast})}\left\Vert \nabla_{\theta}\log Q_{\theta
}(X,s,s^{\prime})\right\Vert ^{2a}\right]  <\infty;
\]
(ii)
\[
E_{\bar{P}_{\ast}^{\nu}}\left[  \sup_{\theta\in B(\delta,\theta_{\ast}%
)}\left\Vert \nabla_{\theta}^{2}\log p_{\theta}(X,s,X^{\prime})\right\Vert
^{2a}\right]  <\infty\text{ and}~E_{\bar{P}_{\ast}^{\nu}}\left[  \sup
_{\theta\in B(\delta,\theta_{\ast})}\left\Vert \nabla_{\theta}^{2}\log
Q_{\theta}(X,s,s^{\prime})\right\Vert ^{2a}\right]  <\infty.
\]

\end{assumption}

\begin{assumption}
\label{ass:q-sum-p} $\sum_{j=0}^{\infty}\left(  E_{\bar{P}_{\ast}^{\nu}%
}\left[  \prod_{i=0}^{j}(1-\underline{q}(X_{i}))^{\frac{2a}{1-a}}\right]
\right)  ^{p\left(  \frac{1-a}{2a}\right)  }<\infty$ for some $p\in(0,2/3)$
(and the same $a\geq1$ that appears in Assumption \ref{ass:deriva-bdd}).
\end{assumption}

\begin{remark}
[Discussion of Assumptions \ref{ass:Theta-int}, \ref{ass:deriva-bdd} and
\ref{ass:q-sum-p}]Part (i) of Assumption~\ref{ass:Theta-int} is standard in
the literature. The restriction that $\Theta_{\ast}$ is a singleton could be
relaxed using the ideas of \cite{liu2003asymptotics} for non-identified ML
estimators. This extension, albeit interesting, would present nuances that are
beyond the scope of the present paper. Part~(ii) of
Assumption~\ref{ass:Theta-int} is also standard, and so is
Assumption~\ref{ass:deriva-bdd} (see \cite{bickel98} for a discussion).
Finally, Assumption~\ref{ass:q-sum-p} is a strengthening of Assumption
\ref{ass:qbar-sum}, and is required in order to establish the existence of a
random sequence $(\Delta_{t}(\theta_{\ast}))_{t}$ which approximates the
\textquotedblleft score\textquotedblright\ function well (in the sense of
Lemma~\ref{lem:score_approx} in Appendix \ref{app:LAR}). As was the case with
Assumption \ref{ass:qbar-sum}, Lemma \ref{lem:suff-A5andA8} in the
Supplemental Material \ref{SM:SuffAssumptions} provides sufficient conditions
of the form $E_{\bar{P}_{\nu}^{\ast}}\left[  (1-\underline{q}(X_{1}%
))^{2a/(1-a)}\right]  <1$. $\triangle$
\end{remark}

The next theorem establishes a LAN-type property for the log-likelihood
criterion function defined in (\ref{logl}).

\begin{theorem}
\label{thm:LAN} Suppose Assumptions \ref{ass:BDD_Q}, \ref{ass:fx-ergo},
\ref{ass:Theta-int}, \ref{ass:deriva-bdd} and \ref{ass:q-sum-p} hold. Then,
there exists a stationary and ergodic process $(\Delta_{t}(\theta_{\ast}%
))_{t}$ in $L^{2}(\bar{P}_{\ast}^{\nu})$, a sequence of negative definite
matrices $(\xi_{t}(\theta_{\ast}))_{t}$, and a compact set $K\subseteq\Theta$
that includes $0$, such that, for any $v\in K$,
\begin{align*}
\ell_{T}^{\nu}(X_{0}^{T},\theta_{\ast}+v)-\ell_{T}^{\nu}(X_{0}^{T}%
,\theta_{\ast})=  &  v^{\intercal}\left(  T^{-1}\sum_{t=0}^{T}\Delta
_{t}(\theta_{\ast})+o_{\bar{P}_{\ast}^{\nu}}(T^{-1/2})\right) \\
&  +(1/2)v^{\intercal}\left(  T^{-1}\sum_{t=0}^{T}\xi_{t}(\theta_{\ast
})+o_{\bar{P}_{\ast}^{\nu}}(1)\right)  v+R_{T}(v),
\end{align*}
where $v\mapsto R_{T}(v)\in\mathbb{R}$ is such that $\lim_{\delta\rightarrow
0}\bar{P}_{\ast}^{\nu}\left(  \sup_{v\in B(\delta,0)}||v||^{-2}R_{T}%
(v)\geq\epsilon\right)  =0$ for any $\epsilon>0$ and any $T\in\mathbb{N}$.
\end{theorem}

\begin{proof}
See Appendix \ref{app:LAR}.
\end{proof}

Theorem~\ref{thm:LAN} extends the results in \cite{bickel96} and
\cite{bickel98} (see their remark on p.~1620) to a more general setup which
allows for covariate-dependent transition probabilities, autoregressive
dynamics, and misspecified models. The proof develops along the same lines as
theirs. The main difference relates to the way in which one establishes that
the \textquotedblleft score\textquotedblright\ $\nabla_{\theta}\ell_{T}^{\nu
}(\cdot,\theta_{\ast})$ and the Hessian $\nabla_{\theta}^{2}\ell_{T}^{\nu
}(\cdot,\theta_{\ast})$ can be approximated by $T^{-1}\sum_{t=0}^{T}\Delta
_{t}(\theta_{\ast})$ and $T^{-1}\sum_{t=0}^{T}\zeta_{t}(\theta_{\ast})$,
respectively (see Lemmas \ref{lem:score_approxV2} and \ref{lem:Hess-approx} in
Appendix \ref{app:LAR}). As mentioned earlier, these approximations rely on
\textquotedblleft mixing\textquotedblright\ properties of the \emph{temporally
inhomogeneous} hidden Markov chain; see Lemma~\ref{lem:Ker-rev} in the
Supplemental Material \ref{SM:LAR}.

Theorem \ref{thm:LAN} may be used to establish the following asymptotic linear
representation for our estimator in terms of $(\Delta_{t}(\theta_{\ast
}))_{t=0}^{\infty}$ and $(\xi_{t}(\theta_{\ast}))_{t=0}^{\infty}$.

\begin{theorem}
\label{thm:anormal} Suppose Assumptions \ref{ass:BDD_Q}--\ref{ass:q-sum-p}
hold and $\eta_{T}=o(T^{-1})$. Then,
\[
\frac{\sqrt{T}(\hat{\theta}_{\nu,T}-\theta_{\ast})}{\sqrt{\mathrm{tr}%
\{\Sigma_{T}(\theta_{\ast})\}}}=-\{E_{\bar{P}_{\ast}^{\nu}}[\xi_{1}%
(\theta_{\ast})]+o_{\bar{P}_{\ast}^{\nu}}(1)\}^{-1}T^{-1/2}\sum_{t=0}^{T}%
\frac{\Delta_{t}(\theta_{\ast})}{\sqrt{\mathrm{tr}\{\Sigma_{T}(\theta_{\ast
})\}}}+o_{\bar{P}_{\ast}^{\nu}}\left(  1\right)  ,
\]
where $\Sigma_{T}(\theta_{\ast})\equiv T^{-1}E_{\bar{P}_{\ast}^{\nu}}%
[\{\sum_{t=0}^{T}\Delta_{t}(\theta_{\ast})\}\{\sum_{t=0}^{T}\Delta_{t}%
(\theta_{\ast})^{\intercal}\}]$.
\end{theorem}

\begin{proof}
See Appendix \ref{app:LAR}.
\end{proof}

Theorem \ref{thm:anormal} readily implies that, if $T^{-1/2}\Sigma_{T}%
(\theta_{\ast})^{-1/2}\sum_{t=0}^{T}\Delta_{t}(\theta_{\ast})\Rightarrow
_{\bar{P}_{\ast}^{\nu}}\mathcal{N}(0,I)$, then
\begin{equation}
\sqrt{T}\Sigma_{T}(\theta_{\ast})^{-1/2}E_{\bar{P}_{\ast}^{\nu}}[\xi
_{1}(\theta_{\ast})](\hat{\theta}_{\nu,T}-\theta_{\ast})\Rightarrow_{\bar
{P}_{\ast}^{\nu}}\mathcal{N}(0,I). \label{eqn:anormtheta}%
\end{equation}
(If $\Sigma_{T}(\theta_{\ast})\rightarrow\Sigma(\theta_{\ast})$, $\Sigma
_{T}(\theta_{\ast})$ may be replaced by $\Sigma(\theta_{\ast})$ in these
statements). This result is akin to results in \cite{white82} and shares the
same features, i.e., the asymptotic covariance matrix $\Omega_{T}(\theta
_{\ast})\equiv(E_{\bar{P}_{\ast}^{\nu}}[\xi_{1}(\theta_{\ast})])^{-1}%
\Sigma_{T}(\theta_{\ast})(E_{\bar{P}_{\ast}^{\nu}}[\xi_{1}(\theta_{\ast
})])^{-1}$ has a \textquotedblleft sandwich\textquotedblright\ form and the
familiar Fisher information equality does not necessarily hold (see also
\cite[Ch. 6]{white94}).

The next theorem complements Theorem \ref{thm:anormal} by presenting results
on consistent estimation of the asymptotic covariance matrix of $\hat{\theta
}_{\nu,T}$, vis., $\Omega_{T}(\theta_{\ast})$, in both correctly specified and
potentially misspecified models. In the latter case, the proposed estimator is
of the \textquotedblleft heteroskedasticity and autocorrelation
consistent\textquotedblright\ type. As in many other results for such
estimators (e.g., \cite{nw87}, \cite[Ch. 8]{white94}), consistency requires
some conditions over the score covariance process $(\Delta_{t+\tau}%
(\cdot)\Delta_{t}(\cdot)^{\intercal})_{t=0}^{T}$, $\tau\in\mathbb{N}$, namely
continuity and law of large numbers results (see Lemma
\ref{lem:Charac.ScoreProcess} in the Supplemental Material \ref{sm:StdErrors}%
). In what follows, the modulus of continuity is defined as
\[
\ddot{\varpi}(\delta^{\prime})\equiv\max_{t}\left\Vert \sup_{||\theta
-\theta_{\ast}||<\delta^{\prime}}\left\Vert \Delta_{t}(\theta)\Delta
_{0}(\theta)^{\intercal}-\Delta_{t}(\theta_{\ast})\Delta_{0}(\theta_{\ast
})^{\intercal}\right\Vert \right\Vert _{L^{1}(\bar{P}_{\ast}^{\nu})},
\]
for any $\delta^{\prime}\in(0,\delta]$, where $\delta$ is as in Assumption
\ref{ass:deriva-bdd}.\footnote{Lemma \ref{lem:Charac.ScoreProcess} in the
Supplementary Material \ref{sm:StdErrors} shows that the modulus of continuity
converges to zero as $\delta$ vanishes.}

\begin{theorem}
\label{thm:StdErrors} Suppose the assumptions of Theorem \ref{thm:anormal}
hold.
%and $E_{\bar{P}_{\ast}^{\nu}}[ \xi_{1}(\theta_{\ast} ) ]$ is non-singular.

\begin{enumerate}
\item[(a)] If $\theta_{\ast}$ is such that, for any $t\geq0$ and $T\geq1$,
$p_{t}^{\nu}(\cdot\mid X_{t-T}^{t-1};\theta^{\ast})=p_{\ast}^{\nu}(\cdot\mid
X_{t-T}^{t-1})$, then
\[
||\Omega_{T}(\theta_{\ast})-\{-H_{T}^{-1}(\hat{\theta}_{\nu,T})\}||=o_{\bar
{P}_{\ast}^{\nu}}(1),
\]
where
\[
H_{T}(\theta)\equiv T^{-1}\sum_{t=1}^{T}\nabla_{\theta}^{2}\log p_{t}^{\nu
}(X_{t}|X_{0}^{t-1},\theta).
\]

\item[(b)] If, for any $l\geq0$,
%$||E_{\bar{P}_{\ast}^{\nu}}[\Delta_{l}(\theta_{\ast})\Delta_{0}(\theta_{\ast})^{\intercal}]||\leq\bar{\upsilon}(l)$ for some integrable function
$|| E_{\bar{P}_{\ast}^{\nu}}[ \Delta_{l}(\theta_{\ast}) \Delta_{0}%
(\theta_{\ast})^{\intercal} ] || \leq\bar{\upsilon}(l)$
%$||\Delta_{l}(\theta_{\ast})\Delta_{0}(\theta_{\ast})^{\intercal}%
%||_{L^{1}(\bar{P}_{\ast}^{\nu})} \leq\bar{\upsilon}(l)$
for some integrable function $l\mapsto\bar{\upsilon}(l)\in\mathbb{R}_{+}$,
then
\[
||\Omega_{T}(\theta_{\ast})-H_{T}^{-1}(\hat{\theta}_{\nu,T})J_{T}(\hat{\theta
}_{\nu,T})H_{T}^{-1}(\hat{\theta}_{\nu,T})||=o_{\bar{P}_{\ast}^{\nu}}(1),
\]

\end{enumerate}

where
\begin{align*}
J_{T}(\theta)  &  \equiv T^{-1}\sum_{t=1}^{T}\nabla_{\theta}\log p_{t}^{\nu
}(X_{t}|X_{0}^{t-1},\theta)\nabla_{\theta}\log p_{t}^{\nu}(X_{t}|X_{0}%
^{t-1},\theta)^{\intercal}\\
&  +\sum_{\tau=1}^{L_{T}}\frac{\omega(\tau,L_{T})}{T-\tau}\sum_{t=1}^{T-\tau
}\nabla_{\theta}\log p_{t+\tau}^{\nu}(X_{t+\tau}|X_{0}^{t+\tau-1}%
,\theta)\nabla_{\theta}\log p_{t}^{\nu}(X_{t}|X_{0}^{t-1},\theta)^{\intercal
}\\
&  +\sum_{\tau=1}^{L_{T}}\frac{\omega(\tau,L_{T})}{T-\tau}\sum_{t=1}^{T-\tau
}\nabla_{\theta}\log p_{t}^{\nu}(X_{t}|X_{0}^{t-1},\theta)\nabla_{\theta}\log
p_{t+\tau}^{\nu}(X_{t+\tau}|X_{0}^{t+\tau-1},\theta)^{\intercal},
\end{align*}
$\omega(\cdot,\cdot)$ are bounded real weights with $\lim_{T\rightarrow\infty
}\omega(\tau,L_{T})=1$ for all $\tau\geq1$, and $(L_{T})_{T=1}^{\infty
}\subseteq\mathbb{N}$ is such that $L_{T}(\ddot{\varpi}(T^{-1/2}\log\log
T)\log\log T+r_{T}+T^{-1/2})=o(1)$, with $(r_{T})_{T}$ being a positive
sequence converging to zero.
\end{theorem}

\begin{proof}
	See Supplemental Material \ref{sm:StdErrors}.
\end{proof}

Part (a) of Theorem~\ref{thm:StdErrors} deals with correctly specified models
and makes use of the Fisher information equality. Part (b) provides a
consistent covariance estimator in the general case of potentially
misspecified models. The conditions for the weights $\omega(\tau,L_{T})$ and
the tuning parameters $(L_{T})_{T}$ are standard for estimators of this
type.\footnote{The additional condition that $\omega(\tau,L_{T})=\sum
_{j=1+\tau}^{L_{T}}c(j,L_{T})c(j-\tau,L_{T})$ for some constants
$c(1,L_{T}),\ldots,c(L_{T},L_{T})$ guarantees that $J_{T}(\hat{\theta}_{\nu
,T})$ is positive semidefinite. Weights obtained from the commonly used
Bartlett, Parzen and quadratic-spectral kernels satisfy this condition.} The
terms $\ddot{\varpi}(T^{-1/2}\log\log T)\log\log T$ and $r_{T}$ in the growth
condition for $L_{T}$ are analogous to those appearing in \cite{nw87} and
arise as \textquotedblleft costs\textquotedblright\ of working with
$\nabla_{\theta}\log p_{t}^{\nu}(X_{t}|X_{0}^{t-1},\hat{\theta}_{\nu,T})$ as
opposed to $\nabla_{\theta}\log p_{t}^{\nu}(X_{t}|X_{0}^{t-1},\theta_{\ast})$,
and of working with sample averages as opposed to population means,
respectively.\footnote{The $\log\log T$ factor is used in order to avoid
working with constants.}
%For instance, in the case of Corollary \ref{cor:anormal}(b), the form of the modulus of continuity can be derived from the primitive smoothness conditions on $p_{\theta}$ and $Q_{\theta}$ and the rate $r_{T}$ can be derived using exponential inequalities for $\beta$-mixing process (e.g. XXX) since the score is $\beta$-mixing (because $(X_{t-1},...,X_{t-L})$ is $\beta$-mixing by Lemma \ref{lem:sta-ergo}).
On the other hand, the $T^{-1/2}$ term does not appear in \cite{nw87} and
arises because we need to approximate $\Delta_{t}(\theta_{\ast})$, that
depends on $X_{-\infty}^{t}$, with the score, that depends on $X_{0}^{t}$.

Theorems~\ref{thm:anormal} and \ref{thm:StdErrors}, together with an
asymptotic normality result such as (\ref{eqn:anormtheta}), provide the means
for constructing asymptotically correct confidence sets and hypotheses tests
for $\theta_{\ast}$. In correctly specified models, $(\Delta_{t}(\theta_{\ast
}))_{t=0}^{\infty}$ is a martingale difference sequence, and thus result
(\ref{eqn:anormtheta}) can be obtained by invoking a martingale central limit
theorem. In potentially misspecified models, $(\Delta_{t}(\theta_{\ast
}))_{t=0}^{\infty}$ will not, in general, be a martingale difference sequence,
so one should use a different approach; in some situations, a central limit
theorem for $\beta$-mixing sequences can be used instead. The following
corollary formalizes this discussion.

\begin{corollary}
\label{cor:anormal} Suppose the assumptions of Theorem \ref{thm:anormal}
hold.
%and $E_{\bar{P}_{\ast}^{\nu}}[ \nabla^{2}_{\theta} \log p_{1}^{\nu}(X_{1} \mid X^{0}_{-\infty} ; \theta_{\ast} ) ]$ is non-singular.

\begin{enumerate}
\item[(a)] If $\theta_{\ast}$ is such that, for any $t\geq0$ and $T\geq1$,
$p_{t}^{\nu}(\cdot\mid X_{t-T}^{t-1};\theta^{\ast})=p_{\ast}^{\nu}(\cdot\mid
X_{t-T}^{t-1})$, then
\[
\sqrt{T}\{-H_{T}(\hat{\theta}_{\nu,T})\}^{1/2}(\hat{\theta}_{\nu,T}%
-\theta_{\ast})\Rightarrow_{\bar{P}_{\ast}^{\nu}}\mathcal{N}(0,I).
\]

\item[(b)] If there exists $\bar{L}>0$ such that, for any $k,T>\bar{L}$,
$p_{k}^{\nu}(X_{k}\mid X_{k-T}^{k-1};\theta^{\ast})=p_{k}^{\nu}(X_{k}\mid
X_{k-\bar{L}}^{k-1};\theta^{\ast})$, $\lim\inf_{T\rightarrow\infty}e_{\min
}(\Sigma_{T}(\theta_{\ast}))>0$,\footnote{For any real symmetric matrix $A$,
$e_{\min}(A)$ denotes its minimum eigenvalue.} and $E_{\bar{P}_{\ast}^{\nu}%
}[\left\Vert \Delta_{1}(\theta_{\ast})\right\Vert ^{4+4\delta}]<\infty$ for
some $\delta>0$, then
\[
\sqrt{T}\hat{\Omega}_{T}(\hat{\theta}_{\nu,T})^{-1/2}(\hat{\theta}_{\nu
,T}-\theta_{\ast})\Rightarrow_{\bar{P}_{\ast}^{\nu}}\mathcal{N}(0,I),
\]
where $\hat{\Omega}_{T}(\hat{\theta}_{\nu,T})\equiv H_{T}^{-1}(\hat{\theta
}_{\nu,T})J_{T}(\hat{\theta}_{\nu,T})H_{T}^{-1}(\hat{\theta}_{\nu,T})$.
\end{enumerate}
\end{corollary}

\begin{proof}
	See Appendix \ref{app:LAR}.
\end{proof}

Regarding the asymptotic normality result, both parts of the corollary rely on
the structure of the process defined by the random variables $\nabla_{\theta
}\log p_{k}^{\nu}(X_{k}\mid X_{-\infty}^{k-1};\theta^{\ast})$; this is a
martingale difference sequence in part (a) and a geometrically $\beta$-mixing
sequence in part (b). The latter result is a consequence of the $\beta$-mixing
structure of $(X_{t})_{t=-\infty}^{\infty}$ and of the fact that the
conditional density at $\theta_{\ast}$ depends on a finite number of lags.
Examples of models for which such properties hold true are mixture models (cf.
Example~\ref{exa:HMM-miss} in the Supplemental Material
\ref{SM:examplemixture}) and mixture autoregressive models (cf.
Example~\ref{exa:MAM}).\footnote{At this level of generality, we cannot
establish an asymptotic normality result in the general case where the model
is misspecified and $\nabla_{\theta}\log p_{k}^{\nu}(X_{k}\mid X_{-\infty
}^{k-1};\theta^{\ast})$ depends on the entire $X_{-\infty}^{k-1}$. The reason
is that, although the process $X_{-\infty}^{\infty}$ is $\beta$-mixing, there
is no guarantee that $\nabla_{\theta}\log p_{k}^{\nu}(X_{k}\mid X_{-\infty
}^{k-1};\theta^{\ast})$ inherits these mixing properties, or that it is a
martingale difference sequence.}

Regarding the covariance estimators used in parts (a) and (b) of the
corollary, these rely on the corresponding parts of Theorem
\ref{thm:StdErrors}. The result in part (b) is established by exploiting the
$\beta$-mixing structure of the score process $(\Delta_{k}(\theta_{\ast
}))_{k=-\infty}^{\infty}$ and the finiteness of its $4+4\delta$ moments (under
$\bar{P}_{\ast}^{\nu}$) for some $\delta>0$. Lemma \ref{lem:suff.HAC.Rate} in
Appendix \ref{app:LAR} shows that $r_{T}=L_{T}T^{-1/2}$, and, since
$\theta\mapsto\log p_{\theta}$ and $\theta\mapsto\log Q_{\theta}$ are smooth
(see Assumptions \ref{ass:Theta-int} and \ref{ass:deriva-bdd}), it follows
that $\delta^{\prime}\mapsto\ddot{\varpi}(\delta^{\prime})=C\delta^{\prime}$
for some finite constant $C$. Hence, the growth condition on $(L_{T})_{T}$
translates into $L_{T}T^{-1/4}\log\log T=o(1)$, which is analogous to that in
\cite{nw87} (apart from the $\log\log T$ factor, which was introduced in
Theorem \ref{thm:StdErrors} for convenience).

The next example verifies the assumptions in the context of the models
considered in Example \ref{exa:Canon}.

\begin{example}
\label{exa:Canon.Inference copy(1)} In view of the results in Example
\ref{exa:Canon}, we only need to verify Assumptions \ref{ass:Theta-int}%
--\ref{ass:q-sum-p}. Part (i) of Assumption \ref{ass:Theta-int} is standard
and is directly imposed, while part (ii) follows from the setup of the
example. Assumption \ref{ass:deriva-bdd} follows by the continuity of the
derivatives. Finally, Lemma \ref{lem:suff-A5andA8} in the Supplemental
Material \ref{SM:consistent} implies that Assumption \ref{ass:q-sum-p} holds.

Thus, Theorem \ref{thm:anormal} holds for the class of models considered in
Example \ref{exa:Canon}. In particular, in the correctly specified case,
Corollary \ref{cor:anormal}(a) guarantees asymptotic normality of the
studentized ML estimator of $\theta_{\ast}$, thereby providing the basis for
inference. These results are, to our knowledge, new in the context of
Markov-switching autoregressive models with covariate-dependent transition
probabilities. $\triangle$
\end{example}

\section{Monte Carlo Simulations}

\label{sec:simulation}

The objective in this section is twofold. First, to assess the quality of
approximations provided by our asymptotic results by examining the
finite-sample properties of the ML estimator and related statistics in a
correctly specified Markov-switching autoregressive model with
covariate-dependent transition probabilities. Second, to explore the effects
of a type of empirically relevant misspecification which involves the use of
an incomplete approximation to the likelihood function that ignores potential
contemporaneous correlation between the observation variable $(Y_{t})$ and the
variable $(Z_{t})$ upon the lagged value of which the transition probabilities depend.

Monte Carlo experiments are based on artificial data $(X_{t}=(Y_{t}%
,Z_{t}))_{t}$ generated according to the equations
\begin{align}
Y_{t}  &  =\mu_{0}(1-S_{t})+\mu_{1}S_{t}+\phi Y_{t-1}+[\sigma_{0}%
(1-S_{t})+\sigma_{1}S_{t}]U_{1,t},\label{mcy}\\
Z_{t}  &  =\mu_{2}+\psi Z_{t-1}+\sigma_{2}U_{2,t}, \label{mcz}%
\end{align}
for $t\in\mathbb{N}$, with $X_{0}=(0.5,0.4)$, $\mu_{0}=-\mu_{1}=1$, $\phi
=0.9$, $\sigma_{0}=\sigma_{1}=1$, $\mu_{2}=0.2$, $\psi=0.8$, $\sigma_{2}=0.6$,
and $(U_{1,t},U_{2,t})_{t}\thicksim$ $\mathrm{i.i.d.}~\mathcal{N}\left(
\left[
\begin{array}
[c]{c}%
0\\
0
\end{array}
\right]  ,\left[
\begin{array}
[c]{cc}%
1 & \rho\\
\rho & 1
\end{array}
\right]  \right)  $ with $\rho\in\{0,0.8\}$. The regimes $(S_{t})_{t}$ are a
Markov chain on $\{0,1\}$, independent of $(U_{1,t},U_{2,t})_{t}$, with
transition probabilities
\begin{equation}
Q_{\theta}(z,s,s)\equiv\Pr(S_{t}=s\mid Z_{t-1}=z,S_{t-1}=s)=[1+\exp
(-\alpha_{s}-\beta_{s}z)]^{-1}, \label{tvtpm}%
\end{equation}
for $s\in\{0,1\}$, where $\alpha_{0}=\alpha_{1}=2$ and $\beta_{0}=-\beta
_{1}=-0.5$. The model defined by (\ref{mcy})--(\ref{tvtpm}) is a prototypical
Markov-switching autoregressive model with covariate-dependent transition
probabilities (cf. Example \ref{exa:MSR}). In each of 1000 independent Monte
Carlo replications, $100+T$ data points for $(X_{t})_{t}$ are generated, with
$T\in\{200,800,1600,3200\}$, and the last $T$ points are used to compute
estimates of the parameters of interest. In order to conserve space, only a
selection of the results are reported (the full set of results is available
upon request).

\subsection{Correct Specification}

\label{sec:simulation1}

In the first set of experiments, we consider estimation of the parameters of
the model in (\ref{mcy})--(\ref{tvtpm}) using the likelihood function based on
the conditional distribution of $X_{t}$ given $X_{0}^{t-1}$. Table
\ref{tab:MC-1n} reports the deviation of the mean of the finite-sample
distributions of the ML estimators of the elements of $\vartheta=(\mu_{0}%
,\mu_{1},\phi,\sigma_{0},\sigma_{1},\alpha_{0},\beta_{0},\alpha_{1},\beta
_{1})$ from the corresponding true parameter values (bias) when $\rho=0.8$. We
also report the ratio of the sampling standard deviation of the estimators to
the estimated standard errors (averaged across replications); the latter are
computed using the Hessian estimator (cf. Theorem~\ref{thm:StdErrors}%
(a)).\footnote{Results for $\rho=0$ are not reported since they are very
similar to those for $\rho=0.8$.} Although the estimators of $\beta_{0}$ and
$\beta_{1}$ are somewhat biased when $T=200$, bias is insignificant in the
rest of the cases. Estimated standard errors are somewhat downwards biased in
most cases but, unless $T$ is small, the bias is not generally substantial and
decreases as $T$ increases.
%Conventional measures of skewness and kurtosis
%based on standardized third and fourth empirical cumulants (not shown) reveal
%that the distributions of the ML estimators of some parameters, especially
%those associated with the transition probabilities, tend to deviate from the
%symmetric and mesokurtic distributions predicted by large-sample theory when
%$T=200$. However, the quality of the Gaussian approximation quickly improves
%as the sample size increases.

{\footnotesize \begin{table}[ptb]
\caption{Bias and Standard Deviation of ML Estimators $(\rho=0.8)$}%
\label{tab:MC-1n}
\begin{center}
{\footnotesize
\begin{tabular}
[c]{cccccccccc}\hline\hline
$T$ & $\mu_{0}$ & $\mu_{1}$ & $\alpha_{1}$ & $\beta_{1}$ & $\alpha_{0}$ &
$\beta_{0}$ & $\sigma_{0}$ & $\sigma_{1}$ & $\phi$\\\hline
\multicolumn{1}{l}{} & \multicolumn{9}{l}{Bias}\\
\multicolumn{1}{l}{200} & -0.007 & 0.023 & 0.004 & 0.126 & 0.043 & 0.209 &
-0.011 & -0.005 & -0.005\\
\multicolumn{1}{l}{800} & -0.002 & 0.004 & -0.002 & 0.023 & 0.002 & 0.026 &
-0.001 & 0.000 & -0.001\\
\multicolumn{1}{l}{1600} & -0.001 & 0.002 & -0.002 & 0.019 & 0.007 & 0.019 &
0.001 & 0.000 & -0.001\\
\multicolumn{1}{l}{3200} & 0.001 & 0.002 & 0.002 & -0.005 & 0.002 & 0.011 &
0.001 & 0.000 & 0.000\\
\multicolumn{1}{l}{} & \multicolumn{9}{l}{Ratio of sampling standard deviation
to estimated standard error}\\
\multicolumn{1}{l}{200} & 1.085 & 1.016 & 1.084 & 1.078 & 1.132 & 1.178 &
1.066 & 1.043 & 1.072\\
\multicolumn{1}{l}{800} & 0.983 & 0.991 & 0.980 & 1.024 & 1.044 & 1.018 &
1.012 & 1.010 & 1.045\\
\multicolumn{1}{l}{1600} & 1.021 & 1.002 & 0.964 & 0.955 & 1.033 & 1.042 &
1.000 & 0.999 & 0.997\\
\multicolumn{1}{l}{3200} & 1.026 & 0.979 & 1.010 & 1.017 & 0.990 & 0.955 &
1.012 & 1.020 & 1.020\\\hline\hline
\end{tabular}
}
\end{center}
\end{table}}

We also examine conventional hypothesis tests for $\vartheta$. Table
\ref{tab:MC-2n} reports the rejection frequencies of: (i)~a $t$-type test of
$\mathcal{H}_{0}:\vartheta_{j}=\vartheta_{j}^{\ast}$ versus $\mathcal{H}%
_{1}:\vartheta_{j}\neq\vartheta_{j}^{\ast}$, where $\vartheta_{j}$ is the
$j$-th element of $\vartheta$ and $\vartheta_{j}^{\ast}$ is its true value;
(ii)$~$a $t$-type test of $\mathcal{H}_{0}:\vartheta_{j}=0$ versus
$\mathcal{H}_{1}:\vartheta_{j}\neq0$. These rejection frequencies are referred
to as \textquotedblleft size\textquotedblright\ and \textquotedblleft
power\textquotedblright, respectively, and are computed using the 0.975
standard-normal quantile as critical value.\footnote{Results should be
interpreted with caution in the case of $\mathcal{H}_{0}:\sigma_{i}=0$,
$i\in\{0,1\}$, because the null value of $\sigma_{i}$ is on the boundary of
the maintained hypothesis. Our asymptotic theory does not allow for parameters
that may lie on the boundary of the parameter space.} Tests tend to have
Type~I error probabilities which are generally close to the nominal 0.05
level, especially for $T>200$. Tests are also powerful enough to reject the
hypothesis of a zero parameter value, except in the case of $\beta_{0}$ and
$\beta_{1}$ with $T=200$. The distributions of studentized statistics
associated with the elements of the ML estimator of $\vartheta$ (ratio of
estimation error to corresponding estimated standard error) generally tend to
have mean and variance (not shown) that do not differ substantially from zero
and one, respectively, and Gaussianity is never rejected for $T>200$%
.\footnote{Statistics associated with $\phi$, $\sigma_{0}$ and $\sigma_{1}$
appear to fare somewhat worse than others when $T=200$, a finding similar to
that reported in \cite{ps98} for models with a time-invariant transition
mechanism.}

{\footnotesize \begin{table}[ptb]
\caption{Size and Power of $t$-Type Tests $(\rho=0.8)$}%
\label{tab:MC-2n}
\begin{center}
{\footnotesize
\begin{tabular}
[c]{cccccccccc}\hline\hline
$T$ & $\mu_{0}$ & $\mu_{1}$ & $\alpha_{1}$ & $\beta_{1}$ & $\alpha_{0}$ &
$\beta_{0}$ & $\sigma_{0}$ & $\sigma_{1}$ & $\phi$\\\hline
\multicolumn{1}{l}{} & \multicolumn{9}{l}{Size}\\
\multicolumn{1}{l}{200} & 0.063 & 0.076 & 0.069 & 0.069 & 0.081 & 0.043 &
0.084 & 0.088 & 0.108\\
\multicolumn{1}{l}{800} & 0.053 & 0.063 & 0.061 & 0.051 & 0.074 & 0.049 &
0.065 & 0.065 & 0.073\\
\multicolumn{1}{l}{1600} & 0.065 & 0.061 & 0.047 & 0.043 & 0.060 & 0.048 &
0.052 & 0.047 & 0.068\\
\multicolumn{1}{l}{3200} & 0.049 & 0.058 & 0.053 & 0.058 & 0.054 & 0.042 &
0.051 & 0.061 & 0.053\\
\multicolumn{1}{l}{} & \multicolumn{9}{l}{Power}\\
\multicolumn{1}{l}{200} & 0.999 & 1.000 & 0.990 & 0.261 & 0.934 & 0.207 &
1.000 & 1.000 & 1.000\\
\multicolumn{1}{l}{800} & 1.000 & 1.000 & 1.000 & 0.812 & 0.999 & 0.732 &
1.000 & 1.000 & 1.000\\
\multicolumn{1}{l}{1600} & 1.000 & 1.000 & 1.000 & 0.988 & 1.000 & 0.966 &
1.000 & 1.000 & 1.000\\
\multicolumn{1}{l}{3200} & 1.000 & 1.000 & 1.000 & 1.000 & 1.000 & 1.000 &
1.000 & 1.000 & 1.000\\\hline\hline
\end{tabular}
}
\end{center}
\end{table}}

\subsection{Misspecification}

\label{sec:simulation2}

In the second set of experiments, we consider estimation of the parameter
$\vartheta=(\mu_{0},\mu_{1},\phi,\sigma_{0},\sigma_{1},\alpha_{0},\beta
_{0},\alpha_{1},\beta_{1})$ using the partial likelihood function based on the
conditional distribution of $Y_{t}$ given $(Y_{t-1},S_{t})$. Inference in
models like (\ref{mcy})--(\ref{tvtpm}) is predominantly based on such a
partial likelihood that ignores the equation for $Z_{t}$ (see, e.g.,
\cite{dieb94}, \cite{filardo94}). Formally, the misspecified model may be
viewed as defined by equations (\ref{mcy})--(\ref{tvtpm}), with the additional
assumption that $\rho=0$. Under this assumption, estimation of $\vartheta$ is
based on (\ref{mcy}) alone, the (potentially incorrect) rational behind this
approach being that, since the conditional distribution of $Y_{t}$ given
$(Y_{t-1},S_{t})$ and the transition probabilities of $(S_{t})_{t}$ depend
only on $Z_{t-1}$, (\ref{mcy}) may be analyzed, without loss of relevant
information, independently of $(Z_{t})_{t}$. Even though such an approach may
be appealing because of its relative simplicity, it is far from obvious that
it provides a valid way for conducting inference on $\vartheta$; it is
unclear, for example, what the limit point of the ML estimator based on the
partial likelihood might be when $\rho\neq0$. In earlier sections, we
considered a theoretical framework that acknowledges this source of
misspecification (among others) and provided tools for asymptotically valid
inference. We now quantify the implications of this misspecification in finite
samples. For brevity, we refer to the maximizer of the partial likelihood
function associated with the conditional distribution of $Y_{t}$ given
$(Y_{t-1},S_{t})$ as the `partial ML' estimator to distinguish it from the
`joint ML' estimator based on the joint model for the conditional distribution
of $X_{t}$ given $X_{0}^{t-1}$ (cf. Section~\ref{sec:simulation1}).

%XXXXX
%In a second set of experiments, we consider estimation of the parameter
%$\vartheta=(\mu_{0},\mu_{1},\phi,\sigma_{0},\sigma_{1},\alpha_{0},\beta
%_{0},\alpha_{1},\beta_{1})$ using the partial likelihood function based on the
%conditional distribution of $Y_{t}$ given $X_{0}^{t-1}$. In empirical
%applications, inference in models like (\ref{mcy})--(\ref{tvtpm}) is
%predominantly based on such a partial likelihood that ignores the equation for
%$Z_{t}$ (cf. \cite{dieb94}, \cite{filardo94}, \cite{rivas15}), the implicit
%assumption being that potential endogeneity of $Z_{t}$ is of little
%consequence since it is only $Z_{t-1}$ that appears in the transition
%mechanism (\ref{tvtpm}). For brevity, we shall refer to this approach as
%limited-information \textquotedblleft partial\textquotedblright\ ML, to
%distinguish it from the full-information \textquotedblleft
%joint\textquotedblright\ ML approach considered in
%Section~\ref{sec:simulation1}.

Table \ref{tab:MC-3n} shows the estimated bias of the partial ML estimators of
the elements of $\vartheta$ and the ratio of the sampling standard deviation
of the estimators to the estimated standard errors (averaged across
replications) when $\rho=0.8$. To reflect what is common practice in applied
research, standard errors are computed using the Hessian estimator (which
relies on the assumption of a correctly specified likelihood) instead of a
\textquotedblleft sandwich\textquotedblright\ estimator (which allows for
mispecification). It is immediately apparent that the partial ML estimator of
most of the parameters is considerably more biased than the joint ML
estimator. The differences between the two estimators are more pronounced for
parameters associated with the transition probabilities ($\alpha_{0}$,
$\beta_{0}$, $\alpha_{1}$, $\beta_{1}$), the partial ML estimators of which
are significantly biased even for $T=3200$. This suggests that the bias of the
partial ML estimator when $\rho\neq0$ is not associated only with small
samples, a finding that is consistent with our asymptotic results. Regarding
the accuracy of estimated standard errors, the latter are downwards biased in
most cases, the bias being somewhat larger than it is for joint ML estimators.
However, unless $T$ is small, this bias is not generally substantial and
declines as $T$ increases, despite the fact that standard errors are obtained
from the Hessian.\footnote{Results for $\rho=0$ (not shown) are not
substantially different from those obtained from the joint ML\ procedure. This
is not surprising since the joint and partial ML estimators are both
consistent for the true parameter value when $\rho=0$.}
{\footnotesize \begin{table}[ptb]
\caption{Bias and Standard Deviation of Partial ML Estimators $(\rho=0.8)$}%
\label{tab:MC-3n}
\begin{center}
{\footnotesize
\begin{tabular}
[c]{cccccccccc}\hline\hline
$T$ & $\mu_{0}$ & $\mu_{1}$ & $\alpha_{1}$ & $\beta_{1}$ & $\alpha_{0}$ &
$\beta_{0}$ & $\sigma_{0}$ & $\sigma_{1}$ & $\phi$\\\hline
\multicolumn{1}{l}{} & \multicolumn{9}{l}{Bias}\\
\multicolumn{1}{l}{200} & 0.019 & 0.046 & -0.172 & 0.898 & 0.397 & 1.363 &
-0.056 & -0.039 & -0.001\\
\multicolumn{1}{l}{800} & 0.017 & 0.017 & -0.156 & 0.485 & 0.216 & 0.838 &
-0.025 & -0.024 & 0.002\\
\multicolumn{1}{l}{1600} & 0.008 & 0.011 & -0.150 & 0.446 & 0.211 & 0.824 &
-0.018 & -0.021 & 0.003\\
\multicolumn{1}{l}{3200} & 0.012 & 0.008 & -0.149 & 0.444 & 0.196 & 0.773 &
-0.020 & -0.019 & 0.003\\
\multicolumn{1}{l}{} & \multicolumn{9}{l}{Ratio of sampling standard deviation
to estimated standard error}\\
\multicolumn{1}{l}{200} & 1.244 & 1.135 & 1.395 & 1.585 & 1.349 & 1.406 &
1.161 & 1.062 & 1.250\\
\multicolumn{1}{l}{800} & 1.027 & 1.041 & 1.153 & 1.170 & 1.067 & 1.119 &
0.988 & 1.020 & 1.080\\
\multicolumn{1}{l}{1600} & 0.991 & 1.014 & 1.070 & 1.075 & 1.072 & 1.161 &
1.006 & 0.959 & 1.054\\
\multicolumn{1}{l}{3200} & 1.047 & 1.027 & 1.050 & 1.063 & 1.034 & 1.088 &
1.021 & 1.004 & 1.056\\\hline\hline
\end{tabular}
}
\end{center}
\end{table}}

We note that hypothesis tests based on studentized statistics analogous to
those considered in Section~\ref{sec:simulation1} (not shown) are unreliable
when partial ML estimates are used. This is especially true in the case of
parameters associated with the transition probabilities, the corresponding
tests being either excessively conservative or excessively liberal. Although
tests of this type are extensively used in applied work, they should be
interpreted with caution since the statistics on which they are based have an
asymptotically normal null distribution only when $\rho=0$. We also note that
using the \textquotedblleft sandwich\textquotedblright\ estimator of
Theorem~\ref{thm:StdErrors}(b) instead of the estimator based on the observed
information matrix is not without difficulty when $\rho\neq0$. As
\cite{Freedman06} points out, the use of such an estimator for inference is
unlikely to produce results that are any less misleading under
misspecification since the problem of bias/inconsistency of the ML estimator
for the true parameter value remains. It is indeed clear from the results in
Table~\ref{tab:MC-3n} that the bias of the partial ML estimator presents a
much more serious problem in our setting than the inaccuracy of conventionally
computed standard errors.

%XXXXX
%In summary, the experiments in Sections <ref>sec:simulation1</ref> and <ref>sec:simulation2</ref> reveal that the
%joint ML estimator has good statistical properties regardless of potential contemporaneous correlation between
%the outcome variable Y_{t} and the information variable Z_{t} driving the hidden Markov transition mechanism.
%By contrast, in the presence of substantial contemporaneous correlation between Y_{t} and Z_{t}, the popular
%partial ML estimator based on the outcome equation alone is severely biased even for what are very large
%sample sizes by the standards of empirical applications.

\section{Empirical Illustration}

\label{sec:illustration}

In this section, we present an empirical illustration based on a
regime-switching model of a type that is commonly used in
economics.\footnote{An additional empirical example, which examines the
predictive ability of an index of leading indicators for regime changes in
U.S. output growth, is discussed in the Supplemental
Material~\ref{SM:empiricalf}.} Specifically, we investigate the potential
contribution of the interest rate spread and the growth in tax revenues in
predicting regime changes in U.S. real output growth. The model is a variant
of the specification used in the simulations and is given by
\begin{align}
Y_{t}  &  =\mu_{0}(1-S_{t})+\mu_{1}S_{t}+\sum_{i=1}^{h_{1}}\phi_{i}%
Y_{t-i}+\sigma_{1}U_{1,t},\label{eqy}\\
Z_{t}  &  =\mu_{2}+\sum_{i=1}^{h_{2}}\psi_{i}Z_{t-i}+\sigma_{2}U_{2,t},
\label{eqz}%
\end{align}
for some $h_{1},h_{2}\in\mathbb{N}$, with the hidden, two-state Markov chain
$(S_{t})_{t}$ being governed by the transition probabilities
\begin{equation}
Q_{\theta}(z,s,s)\equiv\Pr(S_{t}=s\mid Z_{t-1}=z,S_{t-1}=s)=[1+\exp
(-\alpha_{s}-\beta_{s}z)]^{-1}, \label{eqtp}%
\end{equation}
with $s\in\{0,1\}$, and $(U_{1,t},U_{2,t})_{t}$ postulated to be i.i.d.
$\mathcal{N}\left(  \left[
\begin{array}
[c]{c}%
0\\
0
\end{array}
\right]  ,\left[
\begin{array}
[c]{cc}%
1 & \rho\\
\rho & 1
\end{array}
\right]  \right)  $ and independent of $(S_{t})_{t}$. In (\ref{eqy}%
)--(\ref{eqtp}), $Y_{t}$ stands for the growth rate of real gross domestic
product and $Z_{t}$ is either the spread between the 10-year Treasury note
rate and the 3-month Treasury bill rate or the growth rate of real government
receipts of direct and indirect taxes.\footnote{The model could be generalized
to allow for Markov changes in all the parameters. However, since $Z_{t}$ is
thought of here as a potential leading indicator for business-cycle phases, it
does not seem sensible to allow the parameters in both (\ref{eqy}) and
(\ref{eqz}) to be subject to changes driven by $(S_{t})_{t}$. Modeling regime
changes in $(Y_{t})_{t}$ and $(Z_{t})_{t}$ as being driven by two independent
Markov processes is more attractive, but we choose to abstract from this as it
is not directly related to the main problem under study.} The data are
quarterly and span the period 1954:3--2009:2.\footnote{Interest rate data are
taken from the FRED$^{\textregistered}$ database; output and tax data are
taken from \cite{auerbach12}. The likelihood ratio test of \cite{Hansen92}
rejects the hypothesis that $\mu_{0}=\mu_{1}$ in (\ref{eqy}).}

Since the aim is not only to assess the predictive ability of the interest
rate spread and tax revenues for regime changes in output growth but also to
examine whether treating these variables as exogenous yields results which are
different from those obtained from a joint model, we compute two sets of
estimates: partial ML estimates based on (\ref{eqy}) alone and joint ML
estimates based on the system (\ref{eqy})--(\ref{eqz}). We note that in
econometric models of the business cycle such as (\ref{eqy})--(\ref{eqtp}), it
is common to rely on partial ML estimation (see, e.g., \cite{filardo94}).
Parameter estimates are reported in Tables \ref{tab:app-1} and \ref{tab:app-2}%
, with estimated standard errors given in parentheses; the latter are obtained
from the \textquotedblleft sandwich\textquotedblright\ estimator of
Theorem~\ref{thm:StdErrors}(b).\footnote{The weights $\omega(\cdot,L_{T})$ are
obtained from the Parzen kernel and $L_{T}$ is determined using the automatic
plug-in method of \cite{andrews91}.
\par
{}} On the basis of $t$-type tests based on joint ML estimates, at least one
of the parameters $(\beta_{0},\beta_{1})$ is significantly different from
zero, indicating that the spread and tax revenues contain significant
information about the probability of switching between the two regimes.

{\footnotesize \ \begin{table}[ptb]
\caption{ML Estimates (Output Growth, Interest Rate Spread)}%
\label{tab:app-1}%
{\footnotesize  \ \  }
\par
\begin{center}
{\footnotesize \ \
\begin{tabular}
[c]{lllllllll}\hline\hline
\multicolumn{2}{c}{Partial ML} &  &  &  & \multicolumn{4}{c}{Joint ML}\\\hline
$\mu_{0}$ & 0.0091 (0.0012) &  &  &  & $\mu_{0}$ & 0.0092 (0.0012) &  & \\
$\mu_{1}$ & 0.0001 (0.0013) &  &  &  & $\mu_{1}$ & 0.0004 (0.0013) & $\mu_{2}$
& 0.2956 (0.0984)\\
$\phi_{1}$ & 0.1543 (0.0761) &  &  &  & $\phi_{1}$ & 0.1449 (0.0763) &
$\psi_{1}$ & 0.8098 (0.0012)\\
$\phi_{2}$ & 0.0620 (0.0855) &  &  &  & $\phi_{2}$ & 0.0595 (0.0852) &
$\psi_{2}$ & -0.1065 (0.1251)\\
$\phi_{3}$ & -0.0468 (0.0867) &  &  &  & $\phi_{3}$ & -0.0435 (0.0759) &
$\psi_{3}$ & 0.3368 (0.1747)\\
$\phi_{4}$ & -0.0136 (0.0929) &  &  &  & $\phi_{4}$ & -0.0183 (0.0914) &
$\psi_{4}$ & -0.2495 (0.0859)\\
$\alpha_{0}$ & -1.3367 (3.8866) &  &  &  & $\alpha_{0}$ & -1.5210 (2.6935) &
$\sigma_{2}$ & 0.7053 (0.1091)\\
$\beta_{0}$ & 8.8363 (9.4778) &  &  &  & $\beta_{0}$ & 9.1169 (6.4769) &
$\rho$ & -0.0849 (0.1198)\\
$\alpha_{1}$ & 3.1927 (0.7646) &  &  &  & $\alpha_{1}$ & 3.1218 (0.7688) &  &
\\
$\beta_{1}$ & -1.0779 (0.4304) &  &  &  & $\beta_{1}$ & -1.0338 (0.4293) &  &
\\
$\sigma_{1}$ & 0.0077 (0.0006) &  &  &  & $\sigma_{1}$ & 0.0078 (0.0006) &  &
\\\hline\hline
\end{tabular}
}
\end{center}
\end{table}}

Regarding the implications of treating $Z_{t}$ as exogenous, the differences
between partial and joint ML estimates are substantial in the model with tax
revenues (especially for autoregressive coefficients and the parameters
associated with the transition probabilities) but much less so in the model
with the interest rate spread. This is not entirely unexpected in view of the
fact that the estimated value of the conditional correlation $\rho$ is
relatively large (0.6034) in the former model but much smaller ($-$0.0849),
and insignificantly different from zero, in the latter. Such findings are in
line with the analytical and simulation results presented in previous
sections. The relatively large estimate of $\rho$ in Table~\ref{tab:app-2}
also suggests that inference based on the partial ML estimator is potentially
misleading because of the likely bias of the estimator.

{\footnotesize \ \begin{table}[pth]
\caption{ML Estimates (Output Growth, Growth in Taxes)}%
\label{tab:app-2}%
{\footnotesize  \  }
\par
\begin{center}
{\footnotesize \ \
\begin{tabular}
[c]{lllllllll}\hline\hline
\multicolumn{2}{c}{Partial ML} &  &  &  & \multicolumn{4}{c}{Joint ML}\\\hline
$\mu_{0}$ & 0.0071 (0.0014) &  &  &  & $\mu_{0}$ & 0.0081 (0.0014) &  & \\
$\mu_{1}$ & -0.0119 (0.0034) &  &  &  & $\mu_{1}$ & -0.0100 (0.0075) &
$\mu_{2}$ & 0.5455 (0.2633)\\
$\phi_{1}$ & 0.2076 (0.0679) &  &  &  & $\phi_{1}$ & 0.0987 (0.0665) &
$\psi_{1}$ & 0.1723 (0.1220)\\
$\phi_{2}$ & 0.0709 (0.0953) &  &  &  & $\phi_{2}$ & 0.0754 (0.0833) &
$\sigma_{2}$ & 3.0458 (0.3082)\\
$\phi_{3}$ & -0.0530 (0.0754) &  &  &  & $\phi_{3}$ & -0.1168 (0.0610) &
$\rho$ & 0.6034 (0.0673)\\
$\phi_{4}$ & -0.0291 (0.0885) &  &  &  & $\phi_{4}$ & -0.0345 (0.0721) &  & \\
$\alpha_{0}$ & 3.4588 (0.6379) &  &  &  & $\alpha_{0}$ & 3.8835 (1.3467) &  &
\\
$\beta_{0}$ & 0.2754 (0.1237) &  &  &  & $\beta_{0}$ & 0.4061 (0.1379) &  & \\
$\alpha_{1}$ & 0.3852 (0.8704) &  &  &  & $\alpha_{1}$ & -2.5349 (4.7906) &  &
\\
$\beta_{1}$ & 0.2558 (0.1053) &  &  &  & $\beta_{1}$ & 0.0579 (0.1713) &  & \\
$\sigma_{1}$ & 0.0076 (0.0006) &  &  &  & $\sigma_{1}$ & 0.0084 (0.0009) &  &
\\\hline\hline
\end{tabular}
}
\end{center}
\end{table}}

%\bibliographystyle{econometrica}
%{\small {\
%\bibliographystyle{econometrica}
%\bibliographystyle{econometrica}
%\bibliography{markov}
%}}

%\clearpage

\newpage\appendix

\section{Appendix: Proofs}

\label{sec:proofs}

\subsection{Consistency}

\label{app:consistent}

In order to prove Theorem \ref{thm:consistent}, we need two lemmas (the proofs
of which are relegated to the Supplemental Material \ref{SM:consistent}).
%\begin{lemma}
%\label{lem:lik-approx} Suppose Assumptions \ref{ass:BDD_Q}, \ref%
%{ass:qbar-sum} and \ref{ass:bdd}(ii) hold. Then, for any $\epsilon >0$,
%there exists a $T(\epsilon )$ such that
%\begin{equation*}
%\bar{P}_{\ast}^{\nu }\left( \sup_{\theta \in \Theta }\left\vert
%T^{-1}\sum_{t=1}^{T}\left\{\log  p_{t}^{\nu }(X_{t}\mid X_{0}^{t-1},\theta
%)- \log p^{\nu }(X_{t}\mid X_{-\infty }^{t-1},\theta )\right\} \right\vert
%>\epsilon \right) <\epsilon
%\end{equation*}%
%for all $t\geq T(\epsilon )$.
%\end{lemma}
The first lemma shows that the log-likelihood function $\ell_{T}^{\nu}%
(X_{0}^{T},\cdot)$ can be approximated by the sample average of $(\log p^{\nu
}(X_{t}\mid X_{-\infty}^{t-1},\cdot))_{t\in\mathbb{N}}$; this function is used
to construct the function $H^{\ast}$ that defines the pseudo-true parameter
set. The result relies on \textquotedblleft mixing\textquotedblright%
\ properties established in Theorem \ref{thm:Q-ergo} (see Lemma
\ref{lem:approx-l} in the Supplemental Material \ref{SM:consistent}).

\begin{lemma}
\label{lem:lik-approx} Suppose Assumptions \ref{ass:BDD_Q}, \ref{ass:bdd}(ii)
and \ref{ass:qbar-sum} hold. Then,
\[
\sup_{\theta\in\Theta}\left\vert T^{-1}\sum_{t=1}^{T}\left(  \log p_{t}^{\nu
}(X_{t}\mid X_{0}^{t-1},\theta)-\log p^{\nu}(X_{t}\mid X_{-\infty}%
^{t-1},\theta)\right)  \right\vert =o_{\bar{P}_{\ast}^{\nu}}(1).
\]

\end{lemma}

The second lemma essentially establishes a uniform law of large numbers for
the sample average of $(\log p^{\nu}(X_{t}\mid X_{-\infty}^{t-1},\cdot
))_{t\in\mathbb{N}}$.

\begin{lemma}
\label{lem:lik-conv} Suppose Assumptions \ref{ass:BDD_Q}, \ref{ass:fx-ergo},
\ref{ass:exist} and \ref{ass:bdd}(i) hold. Then: (i) For any compact
$K\subseteq\Theta$ and any $\epsilon>0$, there exists $T(\epsilon
)\in\mathbb{N}$ such that
\begin{equation}
\bar{P}_{\ast}^{\nu}\left(  \sup_{\theta\in K}T^{-1}\sum_{t=1}^{T}\left(  \log
p^{\nu}(X_{t}\mid X_{-\infty}^{t-1},\theta)-E_{\bar{P}_{\ast}^{\nu}}\left[
\log p^{\nu}(X_{t}\mid X_{-\infty}^{t-1},\theta)\right]  \right)
>\epsilon\right)  \leq\epsilon, \label{eqn:lik-conv-1}%
\end{equation}
for all $T\geq T(\epsilon)$.\newline

%(ii) For any $\epsilon >0$ and $\theta_{\ast}\in \Theta _{\ast}$, there exists a $%
%T(\epsilon ,\theta_{\ast})$ such that
%\begin{equation*}
%\bar{P}_{\ast}^{\nu }\left( \left\vert T^{-1}\sum_{t=1}^{T}\left(
%-\log p^{\nu }(X_{t}\mid X_{-\infty }^{t-1},\theta_{\ast})+E_{\bar{P}_{\theta
%^{\ast }}^{\nu }}\left[ \log p^{\nu }(X_{t}\mid X_{-\infty }^{t-1},\theta
%_{0})\right] \right) \right\vert >\epsilon \right) \leq \epsilon
%\end{equation*}%
%for all $T\geq T(\epsilon ,\theta_{\ast})$.
(ii) For any $\theta_{\ast}\in\Theta_{\ast}$,
\[
\left\vert T^{-1}\sum_{t=1}^{T}\left(  \log p^{\nu}(X_{t}\mid X_{-\infty
}^{t-1},\theta_{0})-E_{\bar{P}_{\ast}^{\nu}}\left[  \log p^{\nu}(X_{t}\mid
X_{-\infty}^{t-1},\theta_{0})\right]  \right)  \right\vert =o_{\bar{P}_{\ast
}^{\nu}}(1).
\]

\end{lemma}

\begin{proof}[Proof of Theorem \protect\ref{thm:consistent}]
	For simplicity, we set $\eta _{T}=0$ throughout the proof. Formally, we wish
	to establish that, for all $\epsilon >0$, there exists $T(\epsilon )\in
	\mathbb{N}$ such that
	\begin{equation*}
	\bar{P}_{\ast}^{\nu }\left( d_{\Theta }(\hat{\theta}_{\nu
		,T},\Theta _{\ast})\geq \epsilon \right) <\epsilon,
	\end{equation*}%
	for all $t\geq T(\epsilon )$. For this, it suffices to show that there exists a $\theta_{0}\in \Theta _{\ast}$ such that, for any $\epsilon>0$, there exists a $T(\theta_{0},\epsilon)$ such that
	\begin{equation*}
	\bar{P}_{\ast}^{\nu }\left( \sup_{\theta \in \Theta \setminus
		\Theta _{\ast}^{\epsilon }}\ell _{T}^{\nu }(X_{0}^{T},\theta )\geq \ell
	_{T}^{\nu }(X_{0}^{T},\theta_{0})\right) <\epsilon,
	\end{equation*}%
	for all $T \geq T(\theta_{0},\epsilon)$, where $\Theta _{\ast}^{\epsilon }=\{\theta \in
	\Theta \colon d_{\Theta }(\theta ,\Theta _{\ast})<\varepsilon \}$. Since, by Lemma \ref{lem:lik-approx}, $\ell^{\nu}_{T}(X^{T}_{0},\cdot)$ is well
	approximated by $\ell^{\nu}_{T}(X^{T}_{-\infty},\cdot) \equiv T^{-1}
	\sum_{t=1}^{T} \log p^{\nu}(X_{t} \mid X^{t-1}_{-\infty},\cdot)$, it
	suffices to work with the latter function.
	
	Let $A_{T}(\delta) = \left\{ X^{\infty}_{-\infty} \colon \sup_{\theta \in
		\Theta \setminus \Theta_{\ast}^{\epsilon}} T^{-1} \sum_{t=1}^{T} \left( \log
	p^{\nu}(X_{t} \mid X^{t-1}_{-\infty},\theta) - E_{\bar{P}^{\nu}_{\ast}} \left[ \log p^{\nu}(X_{t} \mid X^{t-1}_{-\infty},\theta) \right]
	\right) \leq \delta \right\}$ and $B_{T}(\delta) = \left\{
	X^{\infty}_{-\infty} \colon \left| T^{-1} \sum_{t=1}^{T} \left( - \log
	p^{\nu}(X_{t} \mid X^{t-1}_{-\infty},\theta_{0})+ E_{\bar{P}^{\nu}_{\ast}} \left[ \log p^{\nu}(X_{t} \mid
	X^{t-1}_{-\infty},\theta_{0}) \right] \right) \right| \leq \delta \right\}$, for any $\delta >0$ and any $\theta_{0}\in \Theta _{\ast}$. Observe that
	%\begin{align*}
	%&\bar{P}_{\ast}^{\nu} \left( \sup_{\theta \in \Theta \setminus
	%\Theta_{\ast}^{\epsilon} } \ell^{\nu}_{T}(X^{T}_{-\infty},\theta) \geq
	%\ell^{\nu}_{T}(X^{T}_{-\infty},\theta_{\ast}) \right) \\
	%\leq & \bar{P}_{\ast}^{\nu} \left( \sup_{\theta \in \Theta
	%\setminus \Theta_{\ast}^{\epsilon} } \ell^{\nu}_{T}(X^{T}_{-\infty},\theta)
	%\geq \ell^{\nu}_{T}(X^{T}_{-\infty},\theta_{\ast}) \cap A_{T}(\delta) \cap
	%B_{T}(\delta) \right) + \bar{P}^{\nu}_{\ast} \left(
	%A_{T}(\delta)^{C} \right) + \bar{P}^{\nu}_{\ast} \left(
	%B_{T}(\delta)^{C} \right) \\
	%\leq & \bar{P}_{\ast}^{\nu} \left( \sup_{\theta \in \Theta
	%\setminus \Theta_{\ast}^{\epsilon} } T^{-1} \sum_{t=1}^{T} \left( E_{\bar{P}_{\ast}^{\nu}} \left[ \log \frac{p^{\nu}_{\ast}(X_{t} \mid
	%X^{t-1}_{-\infty})}{p^{\nu}(X_{t} \mid
	%X^{t-1}_{-\infty},\theta)} \right] \right) \leq T^{-1} \sum_{t=1}^{T} \left(
	%E_{\bar{P}_{\ast}^{\nu}} \left[ \log \frac{p^{\nu}_{\ast}(X_{t} \mid
	%X^{t-1}_{-\infty})}{p^{\nu}(X_{t} \mid
	%X^{t-1}_{-\infty},\theta_{\ast})} \right] \right)- 2\delta \right) \\
	%& + \bar{P}_{\ast}^{\nu} \left( A_{T}(\delta)^{C} \right) + \bar{P}_{\ast}^{\nu} \left( B_{T}(\delta)^{C} \right).
	%\end{align*}
	\begin{align*}
	&\bar{P}_{\ast}^{\nu} \left( \sup_{\theta \in \Theta \setminus
		\Theta_{\ast}^{\epsilon} } \ell^{\nu}_{T}(X^{T}_{-\infty},\theta) \geq
	\ell^{\nu}_{T}(X^{T}_{-\infty},\theta_{0}) \right) \\
	\leq & \bar{P}_{\ast}^{\nu} \left( \sup_{\theta \in \Theta
		\setminus \Theta_{\ast}^{\epsilon} } \ell^{\nu}_{T}(X^{T}_{-\infty},\theta)
	\geq \ell^{\nu}_{T}(X^{T}_{-\infty},\theta_{0}) \cap A_{T}(\delta) \cap
	B_{T}(\delta) \right) + \bar{P}^{\nu}_{\ast} \left(
	A_{T}(\delta)^{c} \right) + \bar{P}^{\nu}_{\ast} \left(
	B_{T}(\delta)^{c} \right) \\
	\leq & \bar{P}_{\ast}^{\nu} \left( \sup_{\theta \in \Theta
		\setminus \Theta_{\ast}^{\epsilon} } T^{-1} \sum_{t=1}^{T} E_{\bar{P}^{\nu}_{\ast}} \left[ \log \frac{ p^{\nu}(X_{t} \mid X^{t-1}_{-\infty},\theta)}{p^{\nu}_{\ast}(X_{t} \mid X^{t-1}_{-\infty})} \right]
	\geq T^{-1} \sum_{t=1}^{T} E_{\bar{P}^{\nu}_{\ast}} \left[ \log \frac{ p^{\nu}(X_{t} \mid X^{t-1}_{-\infty},\theta_{0})}{p^{\nu}_{\ast}(X_{t} \mid X^{t-1}_{-\infty})} \right] - 2\delta \right) \\
	& + \bar{P}^{\nu}_{\ast} \left(
	A_{T}(\delta)^{c} \right) + \bar{P}^{\nu}_{\ast} \left(
	B_{T}(\delta)^{c} \right)\\
	\leq & \bar{P}_{\ast}^{\nu} \left( \inf_{\theta \in \Theta
		\setminus \Theta_{\ast}^{\epsilon} } H^{\ast}(\theta)  \leq H^{\ast}(\theta_{0}) + 2\delta \right) + \bar{P}_{\ast}^{\nu} \left( A_{T}(\delta)^{c} \right) + \bar{P}_{\ast}^{\nu} \left( B_{T}(\delta)^{c} \right),
	\end{align*}
	where the last line follows from the stationarity of $X_{-\infty}^{\infty}$ and the definition of $H^{\ast}$. By Assumption \ref{ass:exist} and the fact that, for any $\theta \in \Theta \setminus \Theta _{\ast}^{\epsilon }$, $H^{\ast }(\theta )>H^{\ast }(\theta_{0})$ (otherwise, $\theta $ would
	belong to $\Theta _{\ast}$), it follows that $\inf_{\theta \in \Theta
		\setminus \Theta_{\ast}^{\epsilon} } H^{\ast}(\theta) - H^{\ast}(\theta_{0}) \equiv \Delta > 0$. Hence, choosing $\delta < 0.5\Delta$, the first term of the right-hand side (RHS) vanishes. By Assumption \ref{ass:exist}(i), $\Theta \setminus \Theta_{\ast}^{\epsilon}$
	is compact; hence, by Lemma \ref{lem:lik-conv}, there exists a $T^{\prime}$
	(which may depend on $\epsilon$ and $\theta_{0}$) such that $\bar{P}_{\ast}^{\nu} \left( A_{T}(\delta)^{c} \right) + \bar{P}_{\ast}^{\nu} \left( B_{T}(\delta)^{c} \right) \leq \epsilon$
	for any $\delta \leq 0.5\epsilon$ and all $T\geq T^{\prime}$, and thus the desired result follows.
\end{proof}

\subsection{Mixing Results}

\label{app:mixing}

Throughout, fix $m$ and $j$ as in the statement of Theorem \ref{thm:Q-ergo}.
For any $n,n^{\prime}$ such that $-m\leq n,n^{\prime}\leq j+1$, we denote the
Dobrushin coefficient of $\bar{P}_{\theta}^{\nu}(S_{n^{\prime}}=\cdot\mid
S_{n}=\cdot,X_{-m}^{j})$ as {\small {
\begin{equation}
\label{eqn:Dobrushin}\alpha_{\theta,n^{\prime},n}(X_{-m}^{j})\equiv\tfrac
{1}{2}\max_{(a,b)\in\mathbb{S}^{2}}\left\Vert \bar{P}_{\theta}^{\nu
}(S_{n^{\prime}}=\cdot\mid S_{n}=a,X_{-m}^{j})-\bar{P}_{\theta}^{\nu
}(S_{n^{\prime}}=\cdot\mid S_{n}=b,X_{-m}^{j})\right\Vert _{1}.
\end{equation}
}}

Since $\alpha_{\theta,j+1,-m}(X_{-m}^{j})\leq\prod_{n=-m}^{j}\alpha
_{\theta,n+1,n}(X_{-m}^{j})$ (e.g., \cite{dobrushin1956central},
\cite{sethuraman2005martingale}), to prove Theorem \ref{thm:Q-ergo}, it
suffices to show the following.

\begin{lemma}
\label{lem:Dobrushin-l.bound} For any $l\in\{-m,\ldots,j\}$ and any $\theta
\in\Theta$, $\alpha_{\theta,l+1,l}(X_{-m}^{j})\leq1-\underline{q}(X_{l})$
a.s.-$\bar{P}_{\ast}^{\nu}$.
\end{lemma}

Lemma~\ref{lem:Dobrushin-l.bound} follows immediately from Lemmas
\ref{lem:couple} and \ref{lem:ergo-eta} below. To state these lemmas, we
construct the following processes that will be used for coupling. For any
$i\in\{1,2\}$ and any $\theta\in\Theta$, let $(X_{i,t},\eta_{i,t}%
,\upsilon_{i,t})_{t=-m}^{\infty}$, with $(X_{i,t},\eta_{i,t},\upsilon
_{i,t})\in\mathbb{X}\times\mathbb{S}\times\{0,1\}$, be defined as follows:
$(X_{i,-m},\eta_{i,-m})\sim\nu$; given $(X_{i,t})_{-m}^{\infty}$,
$(\upsilon_{i,t})_{t=-m}^{\infty}$ is i.i.d. with $\Pr(\upsilon_{i,t}=1\mid
X_{i,-m}^{\infty})=\Pr(\upsilon_{i,t}=1\mid X_{i,-m}^{t})\equiv\underline{q}%
(X_{i,t})$; for each $t\geq-m$, $\eta_{i,t+1}\sim\varrho(X_{i,t},\cdot)$ if
$\upsilon_{i,t}=1$, and $\eta_{i,t+1}\sim\frac{Q_{\theta}(X_{i,t},\eta
_{i,t},\cdot)-\underline{q}(X_{i,t})\varrho(X_{i,t},\cdot)}{1-\underline{q}%
(X_{i,t})}$ if $\upsilon_{i,t}=0$ (the last quotient expression is a valid
transition kernel under condition (\ref{eqn:Q-ergo}) in Theorem
\ref{thm:Q-ergo}); finally, $X_{i,t+1}\sim p_{\theta}(X_{i,t},\eta
_{i,t+1},\cdot)$.

This construction implies that the transition kernel of $(\eta_{i,t})_{t}$ is
given by
\begin{align*}
\mathrm{Pr}_{\theta}\left(  \eta_{i,t+1}=\cdot\mid\eta_{i,t},X_{i,t}\right)
=  &  \underline{q}(X_{i,t})\varrho(X_{i,t},\cdot)+(1-\underline{q}%
(X_{i,t}))\frac{Q_{\theta}(X_{i,t},\eta_{i,t},\cdot)-\underline{q}%
(X_{i,t})\varrho(X_{i,t},\cdot)}{1-\underline{q}(X_{i,t})}\\
=  &  Q_{\theta}(X_{i,t},\eta_{i,t},\cdot),
\end{align*}
and since the transition for $X_{i,t+1}$ given $(X_{i,t},\eta_{i,t+1})$ is
governed by $p_{\theta}$, the following result holds (its proof is relegated
to the Supplemental Material \ref{SM:mixing}).

\begin{lemma}
\label{lem:couple} For any $l\in\{-m,\ldots,j\}$ and any $\theta\in\Theta$,
\[
\bar{P}_{\theta}^{\nu}(S_{l+1}=\cdot\mid S_{l}=\cdot,X_{-m}^{j})=\mathrm{Pr}%
_{\theta}(\eta_{i,l+1}=\cdot\mid\eta_{i,l}=\cdot,X_{-m}^{j}),~\forall
i\in\{1,2\},
\]
$a.s.$\textrm{-}$\bar{P}_{\ast}^{\nu}$.
\end{lemma}

Furthermore, since $\eta_{i,t+1}$ becomes independent of its own past whenever
$\upsilon_{i,t}=1$, the following result can be established (its proof is
relegated to the Supplemental Material \ref{SM:mixing}).

\begin{lemma}
\label{lem:ergo-eta} For any $l\in\{-m,\ldots,j\}$ and any $\theta\in\Theta$,
\[
\tfrac{1}{2}\max_{(a,b)\in\mathbb{S}^{2}}\left\Vert \mathrm{Pr}_{\theta}%
(\eta_{1,l+1}=\cdot\mid\eta_{1,l}=a,X_{-m}^{j})-\mathrm{Pr}_{\theta}%
(\eta_{2,l+1}=\cdot\mid\eta_{2,l}=b,X_{-m}^{j})\right\Vert _{1}\leq
1-\underline{q}(X_{l}),
\]
$a.s.$\textrm{-}$\bar{P}_{\ast}^{\nu}$.
\end{lemma}

It is easy to see that Lemma \ref{lem:Dobrushin-l.bound} (and thus Theorem
\ref{thm:Q-ergo}) follows from the last two lemmas.

%		
%		
%		
%		
%		
%		For any $k,l,n$ in $\mathbb{N}$ and $-m \leq \min \{ l,n ,k \}$, let
%		\begin{align*}
%		x^{k}_{m} \mapsto \omega^{(i)}_{l,n}(x^{k}_{-m})(b|a) \equiv \Pr \left(  \eta_{i,l} = b \mid \eta_{i,n} = a , x^{k}_{-m} \right)
%		\end{align*}
%		for any $a,b$ in $\mathbb{S}$.

\subsection{Asymptotic Distribution Theory}

\label{app:LAR}

The next two lemmas are used to prove Theorems \ref{thm:LAN} and
\ref{thm:anormal} (their proofs are relegated to the Supplemental Material
\ref{SM:LAR1}). In what follows, for vector/matrix-valued functions $X\mapsto
f(X)$, $\left\Vert f\right\Vert _{L^{r}(P)}$ is short-hand notation for the
$L^{r}(P)$-norm of $x\mapsto\left\Vert f(x)\right\Vert $, where $\left\Vert
\cdot\right\Vert $ denotes the Euclidean/dual norm of $f$.

\begin{lemma}
\label{lem:score_approxV2} Suppose Assumptions \ref{ass:BDD_Q},
\ref{ass:fx-ergo}, \ref{ass:Theta-int}, \ref{ass:deriva-bdd}(i) and
\ref{ass:q-sum-p} hold. Then, there exists a stationary and ergodic (under
$\bar{P}_{\ast}^{\nu}$) process $(\Delta_{t}(\theta_{\ast}))_{t=-\infty
}^{\infty}$ in $L^{2}(\bar{P}_{\ast}^{\nu})$ such that
\[
\lim_{T\rightarrow\infty}\left\Vert T^{-1/2}\sum_{t=0}^{T}\{\nabla_{\theta
}\log p_{t}^{\nu}(\cdot|\cdot,\theta_{\ast})-\Delta_{t}(\theta_{\ast
})\}\right\Vert _{L^{2}(\bar{P}_{\ast}^{\nu})}=0.
\]

\end{lemma}

\begin{lemma}
\label{lem:Hess-approx} Suppose Assumptions \ref{ass:BDD_Q}, \ref{ass:fx-ergo}%
, \ref{ass:Theta-int}, \ref{ass:deriva-bdd} and \ref{ass:q-sum-p} hold. Then,
there exists a sequence of $\mathbb{R}^{q\times q}$-valued continuous
functions $(\theta\mapsto\xi_{t}(\theta))_{t}$ such that $\xi_{t}(\theta)$ is
negative definite for all $t$ and
\[
\lim_{T\rightarrow\infty}\left\Vert \sup_{\theta\in B(\delta,\theta_{\ast}%
)}||T^{-1}\sum_{t=0}^{T}\{\nabla_{\theta}^{2}\log p_{t}^{\nu}(\cdot\mid
\cdot,\theta)-\xi_{t}(\theta)\}||\right\Vert _{L^{1}(\bar{P}_{\ast}^{\nu}%
)}=0,
\]
where $\delta>0$ is the same as in Assumption \ref{ass:deriva-bdd}.
\end{lemma}

\begin{proof}[Proof of Theorem \protect\ref{thm:LAN}]
	Choose $K$ compact such that, for any $v \in K$, $||v|| \leq \delta$ for $\delta>0$ as in Lemma \ref{lem:Hess-approx}. For any $v \in K$, by Assumption \ref{ass:Theta-int},
	\begin{equation*}
	\ell^{\nu} _{T}(X_{0}^{T},\theta_{\ast}+v)-\ell^{\nu} _{T}(X_{0}^{T},\theta
	_{\ast})=v^{\intercal} \nabla _{\theta }\ell^{\nu} _{T}(X_{0}^{T},\theta
	_{\ast})+0.5v^{\intercal}\left( \int_{0}^{1}\nabla _{\theta }^{2}\ell^{\nu}
	_{T}(X_{0}^{T},\theta_{\ast}+wv)dw\right) v.
	\end{equation*}%
	By Lemmas \ref{lem:score_approxV2} and \ref{lem:Hess-approx}, and the fact that ($\theta_{\ast} + wv) \in B(v,\theta_{\ast})$,
	\begin{align*}
	\ell^{\nu} _{T}(X_{0}^{T},\theta_{\ast}+v)-\ell^{\nu} _{T}(X_{0}^{T},\theta
	_{\ast})=& v^{\intercal} \left( T^{-1}\sum_{t=0}^{T}\Delta _{t}(\theta_{\ast})+o_{\bar{P}_{\ast}^{\nu }}(T^{-1/2})\right) \\
	& +0.5v^{\intercal}\left( \int_{0}^{1} T^{-1}\sum_{t=0}^{T}\xi _{t}(\theta_{\ast} + wv ) dw + o_{\bar{P}_{\ast}^{\nu }}(1)\right) v.
	\end{align*}
	
	\noindent Now, let $R_{T}(v) \equiv v^{\intercal} \left( T^{-1} \sum_{t=0}^{T} \int_{0}^{1}
	\{\xi_{t}(\theta_{\ast} + w v) - \xi_{t}(\theta_{\ast}) \} dw \right) v $. Observe
	that $||v||^{-2}|R_{T}(v)| \leq \int_{0}^{1} \left \Vert T^{-1}
	\sum_{t=0}^{T} \{\xi_{t}(\theta_{\ast} + w v) - \xi_{t}(\theta_{\ast}) \}
	\right
	\Vert dw $, so, for any $\delta>0$,
	\begin{align*}
	\bar{P}_{\ast}^{\nu} \left( \sup_{v \in B(\delta,0) } \frac{|R_{T}(v)|}{||v||^{2} }\geq \epsilon
	\right) \leq & 	\bar{P}_{\ast}^{\nu} \left( \sup_{v \in B(\delta,0)} \int_{0}^{1} \left \Vert T^{-1}
	\sum_{t=0}^{T} \{\xi_{t}(\theta_{\ast} + w v) - \xi_{t}(\theta_{\ast}) \}
	\right 	\Vert dw  \geq \epsilon 	\right)\\
%	\leq &  \epsilon^{-1} E_{\bar{P}_{\ast}^{\nu}} \left[  \sup_{v \in B(\delta,0)} \int_{0}^{1} \left \Vert T^{-1}
%	\sum_{t=0}^{T} \{\xi_{t}(\theta_{\ast} + s v) - \xi_{t}(\theta_{\ast}) \}
%	\right 	\Vert ds 	\right] \\
%	\leq &  \epsilon^{-1} T^{-1}	\sum_{t=0}^{T} E_{\bar{P}_{\ast}^{\nu}} \left[  \sup_{v \in K} \int_{0}^{1}   \left \Vert \{\xi_{t}(\theta_{\ast} + s v) - \xi_{t}(\theta_{\ast}) \}\right 	\Vert ds 	\right]\\
	\leq & \epsilon^{-1} E_{\bar{P}_{\ast}^{\nu}} \left[ \sup_{v \in B(\delta,0)} \int_{0}^{1} \left \Vert  \xi_{1}(\theta_{\ast} + w v) - \xi_{1}(\theta_{\ast})
	\right 	\Vert dw 	\right],
	\end{align*}
	where the second line follows by the Markov inequality and stationarity. The desired result then follows by the continuity of $\xi_{1}$ (see Lemma \ref{lem:Hess-approx}) and the same arguments as in \cite[p. 1634]{bickel98}.
\end{proof}

\begin{proof}[Proof of Theorem \protect\ref{thm:anormal}]
	Henceforth, let $\bar{\Delta}_{T} \equiv T^{-1}
	\sum_{t=0}^{T}\Delta_{t}(\theta_{\ast}) + o_{\bar{P}_{\ast}^{\nu}}(
	T^{-1/2}) $. By Theorem~\ref{thm:consistent}, $\hat{\theta}_{\nu ,T} - \theta
	_{\ast}$ converges to zero with probability approaching one (w.p.a.1). Thus, $R_{T}(v_{T}) = o_{\bar{P}_{\ast}^{\nu}}(||v_{T}||^{2})$ and, by Theorem \ref{thm:LAN},
	\begin{align*}
	\ell _{T}^{\nu }(X_{0}^{T},\hat{\theta}_{\nu ,T})-\ell _{T}^{\nu
	}(X_{0}^{T},\theta_{\ast}) =& (\hat{\theta}_{\nu ,T}-\theta_{\ast})^{\intercal}\bar{\Delta}_{T} \\
	& +0.5(\hat{\theta}_{\nu ,T}-\theta_{\ast})^{\intercal}\left(  T^{-1}
	\sum_{t=0}^{T}\xi_{t}(\theta_{\ast}) +o_{\bar{P}_{\ast}^{\nu
	}}(1)\right) (\hat{\theta}_{\nu ,T}-\theta_{\ast}).
	\end{align*}	
	Ergodicity of $X_{-\infty}^{\infty}$ (Lemma~\ref{lem:sta-ergo}) implies ergodicity of $(\xi _{t}(\theta_{\ast}))_{t=-\infty}^{\infty} $; therefore, by Lemma~\ref{lem:Hess-approx} and Birkhoff's ergodic theorem,
	\begin{align}\notag
\ell _{T}^{\nu }(X_{0}^{T},\hat{\theta}_{\nu ,T})-\ell _{T}^{\nu
}(X_{0}^{T},\theta_{\ast}) =& (\hat{\theta}_{\nu ,T}-\theta_{\ast})^{\intercal}\bar{\Delta}_{T} \\ \label{eqn:Quad-1}
& +0.5(\hat{\theta}_{\nu ,T}-\theta_{\ast})^{\intercal}\left(  E_{\bar{P}_{\ast}^{\nu
}}[\xi_{1}(\theta_{\ast})]   +o_{\bar{P}_{\ast}^{\nu
}}(1)\right) (\hat{\theta}_{\nu ,T}-\theta_{\ast}),
\end{align}		
and $E_{\bar{P}_{\ast}^{\nu}}[\xi_{1}(\theta_{\ast})] $ is non-singular. The rest of the proof proceeds in two steps.
	
\bigskip

	\textsc{Step 1.}
%Let $r_{T} \equiv \min \{ o _{\bar{P}_{\ast}^{\nu }}(1) , o(T^{-1}) +  E_{\bar{P}_{\ast}^{\nu }} \left[  (\bar{\Delta}_{T})^{\prime}  (E_{\bar{P}_{\ast}^{\nu}}[\xi_{1}(\theta_{\ast})] )^{-1} (\bar{\Delta}_{T})     \right]   \}$.
	Let $r_{T} \equiv \min \{ o _{\bar{P}_{\ast}^{\nu }}(1) , o(T^{-1}) +  E_{\bar{P}_{\ast}^{\nu }} \left[  (\bar{\Delta}_{T})^{\intercal} (\bar{\Delta}_{T})     \right]   \}$.  We first establish that $||\hat{\theta}_{\nu ,T}-\theta
	_{\ast}||=O_{\bar{P}_{\ast}^{\nu }}(\sqrt{r_{T}})$; by Theorem \ref{thm:consistent}, the $o _{\bar{P}_{\ast}^{\nu }}(1) $ part of $r_{T}$ has been established.
	
	By (\ref{eqn:Quad-1}) and the fact that $\hat{\theta}_{\nu,T}$ is an (approximate) maximizer of the likelihood function,
	\begin{equation*}
-2\eta _{T}\leq 2(\hat{\theta}_{\nu ,T}-\theta_{\ast})^{\intercal} \bar{\Delta}_{T}-(\hat{\theta}_{\nu ,T}-\theta_{\ast})^{\intercal} A(\theta_{\ast}) (\hat{\theta}_{\nu
	,T}-\theta_{\ast}),
\end{equation*}
with $A(\theta_{\ast}) \equiv \left( - E_{\bar{P}_{\ast}^{\nu
}}[\xi_{1}(\theta_{\ast})]   +o_{\bar{P}_{\ast}^{\nu
}}(1)\right)$. Simple algebra yields
	\begin{equation*}
-2\eta _{T}\leq - \left \Vert  (\hat{\theta}_{\nu ,T}-\theta_{\ast})^{\intercal} A(\theta_{\ast})^{1/2} -  \bar{\Delta}_{T}^{\intercal}  A(\theta_{\ast})^{-1/2}  \right \Vert^{2}  + \bar{\Delta}_{T}^{\intercal}  A(\theta_{\ast})^{-1}  \bar{\Delta}_{T}.
\end{equation*}		
Moreover, by simple algebra and the Markov inequality,
\begin{align*}
	\left \Vert  (\hat{\theta}_{\nu ,T}-\theta_{\ast})^{\intercal} A(\theta_{\ast})^{1/2} \right \Vert = O_{\bar{P}_{\ast}^{\nu
	}}  \left( \sqrt{ \eta_{T}}  + \sqrt{E_{\bar{P}_{\ast}^{\nu }} \left[  (\bar{\Delta}_{T})^{\intercal}  (E_{\bar{P}_{\ast}^{\nu}}[\xi_{1}(\theta_{\ast})] )^{-1} (\bar{\Delta}_{T})     \right]}    \right).
\end{align*}
This expression, the fact that $E_{\bar{P}_{\ast}^{\nu}}[\xi_{1}(\theta_{\ast})]$ is non-singular, and $\eta_{T}= o(T^{-1})$, imply the desired result.
	\bigskip
	
\textsc{Step 2.} We now show that, for any $\epsilon>0$,
\begin{align*}
\bar{P}_{\ast}^{\nu} \left( r_{T}^{-1/2} \left \Vert (\hat{\theta}%
_{\nu,T} - \theta_{\ast}) - (-   E_{\bar{P}_{\ast}^{\nu }} [\xi_{1}(\theta_{\ast})]  + o_{\bar{P}^{\nu}_{\ast}}(1) )^{-1} \bar{\Delta}_{T} \right \Vert \geq \epsilon \right) \rightarrow
0.
\end{align*}
Since, by Step 1, $\left\Vert \hat{\theta}_{\nu ,T}-\theta_{\ast}\right\Vert =O_{\bar{%
		P}_{\ast}^{\nu }}(\sqrt{r_{T}})$, it suffices to show that
\begin{equation}
\bar{P}_{\ast}^{\nu }\left( \left\{ r^{-1/2}_{T} \left\Vert (\hat{\theta%
}_{\nu ,T}-\theta_{\ast})-(- E_{\bar{P}_{\ast}^{\nu }} [\xi_{1}(\theta_{\ast})]   +o_{\bar{P}_{\ast}^{\nu }}(1))^{-1} \bar{\Delta}_{T} \right\Vert \geq \epsilon \right\} \cap
\{   \left\Vert \hat{\theta}_{\nu ,T}-\theta_{\ast}\right\Vert \leq
\sqrt{r_{T}} M\}\right) \rightarrow 0,  \label{eqn:anorm-1}
\end{equation}%
where $M>0$. To this end, note that, by Theorem \ref{thm:LAN} and the fact that $T^{-1} \sum_{t=1}^{T} \xi_{t}(\theta_{\ast}) =  E_{\bar{P}_{\ast}^{\nu }} [\xi_{1}(\theta_{\ast})] + o_{\bar{P}_{\ast}^{\nu }}(1)$, it follows that
\begin{equation*}
\ell^{\nu} _{T}(X_{0}^{T},\theta_{\ast}+v)-\ell^{\nu} _{T}(X_{0}^{T},\theta
_{\ast})=(\bar{\Delta}_{T})^{\intercal}v-0.5v^{\intercal}(-E_{\bar{P}_{\ast}^{\nu }} [\xi_{1}(\theta_{\ast})] +o_{\bar{P}%
	_{\ast}^{\nu }}(1))v+R_{T}(v),
\end{equation*}%
for any $v\in K$.
Letting $\Lambda _{T}(v)\equiv \ell^{\nu}
_{T}(X_{0}^{T},\theta_{\ast}+v)-\ell^{\nu} _{T}(X_{0}^{T},\theta_{\ast})$ and $Q_{T}(v)\equiv (\bar{\Delta}_{T} )^{\intercal}v-0.5v^{\intercal}(- E_{\bar{P}_{\ast}^{\nu }} [\xi_{1}(\theta_{\ast})] +o_{\bar{P}_{\ast}^{\nu }}(1))v$, we show that $\sup_{v\in \{v\colon ||v||\leq
	\sqrt{r_{T}} M\}}  r^{-1}_{T}  \left\vert \Lambda _{T}(v)-Q_{T}(v)\right\vert  = o_{\bar{P}_{\ast}^{\nu }}(1)$. To do so, it suffices to prove that $\sup_{v\in \{v\colon ||v||\leq
	\sqrt{r_{T}} M\}} \left\vert R_{T}(v)\right\vert = o_{\bar{P}_{\ast}^{\nu }}(r_{T})$. But this follows from Theorem \ref{thm:LAN} and the fact that $\sqrt{r_{T}} = o_{\bar{P}_{\ast}^{\nu }}(1)$.
%Thus
%\begin{align*}
%\bar{P}_{\ast}^{\nu} \left( \sup_{v \in \{v \colon ||v|| \leq
%	\sqrt{r_{T}}  M \} } \left \vert \Lambda_{T}(v) - Q_{T}(v) \right \vert \geq \sqrt{r_{T}} \epsilon \right) \rightarrow 0.
%\end{align*}
Since $(\hat{\theta}_{\nu ,T}-\theta_{\ast})\in \{v\colon ||v||\leq \sqrt{r_{T}} M\} $ and maximizes $\Lambda _{T}(\cdot )$ (within a $\eta _{T}$
margin), the previous result implies that
\begin{align*}
\hat{\theta}_{\nu ,T}-\theta_{\ast}=& \arg \max_{v\in \{v\colon ||v||\leq
	T^{-1/2}M\}}Q_{T}(v)+o_{\bar{P}_{\ast}^{\nu}}(\sqrt{r_{T}})+\eta _{T} \\
=& (- E_{\bar{P}_{\ast}^{\nu }} [\xi_{1}(\theta_{\ast})]  +o_{\bar{P}_{\ast}^{\nu }}(1))^{-1} \bar{\Delta}_{T} +o_{\bar{P}_{\ast}^{\nu }}(\sqrt{r_{T}}),
\end{align*}%
and thus (\ref{eqn:anorm-1}) follows.	
\end{proof}

The proof of Corollary \ref{cor:anormal} uses the following lemma (whose proof
is relegated to the Supplemental Material \ref{SM:exa.Canon.Inference}).

\begin{lemma}
\label{lem:suff.HAC.Rate} Suppose there exists $\bar{L}\in\mathbb{N}$ such
that $\nabla_{\theta}\log p_{t}^{\nu}(X_{t}|X_{0}^{t-1},\theta)=\nabla
_{\theta}\log p_{t}^{\nu}(X_{t}|X_{t-\bar{L}}^{t-1},\theta)$ for all $t\geq0$. Then:

\begin{itemize}
\item[(a)] For all $t\geq0$,
\[
\Delta_{t}(\theta_{\ast})\equiv\Delta_{t,-\infty}(\theta_{\ast})(X_{-\infty
}^{t})=\nabla_{\theta}\log p_{t}^{\nu}(X_{t}|X_{t-\bar{L}}^{t-1},\theta_{\ast
}).
\]

\item[(b)] If, in addition, there exists $\delta>0$ such that $E_{\bar
{P}_{\ast}^{\nu}}[||\nabla_{\theta}\log p_{1}^{\nu}(X_{1}|X_{1-\bar{L}}%
^{0},\theta_{\ast})||^{4+4\delta}]<\infty$, then, for any $L\in\mathbb{N}$,
{\small {
\begin{align*}
\max_{j\in\{0,\ldots,L\}}||T^{-1}\sum_{t=1}^{T}\Delta_{t+j,-\infty}%
(\theta_{\ast})\Delta_{t,-\infty}(\theta_{\ast})^{\intercal}-E_{\bar{P}_{\ast
}^{\nu}}[\Delta_{j,-\infty}(\theta_{\ast})\Delta_{0,-\infty}(\theta_{\ast
})^{\intercal}]|| =O_{\bar{P}_{\ast}^{\nu}}\left(  \frac{L}{\sqrt{T}}\right)
.
\end{align*}
}}
\end{itemize}
\end{lemma}

\begin{proof}[Proof of Corollary \ref{cor:anormal}]
	
	Throughout the proof, we use $C$ to denote a universal constant that can take different values in different places. Also,  for any $k,T\geq 0$ and any $\theta \in \Theta $, we write $\Delta
	_{k,k-T}(\theta )\equiv \Delta _{k,k-T}(\theta )(X_{k-T}^{k}) = \nabla _{\theta }\log p_{k}^{\nu
	}(X_{k}|X_{k-T}^{k-1};\theta )$ (for the last equality, see Lemma \ref{lem:score-rep} in the Supplemental Material \ref{SM:LAR1}) and $\Delta _{k}(\theta)\equiv \Delta _{k,-\infty }(\theta)$.
	
	We first show that, under the conditions of part (a), $(\Delta_{t}(\theta_{\ast}))_{t}$ is a martingale difference sequence (MDS) with respect to the natural filtration of $(X_{t})_{t}$. To establish this, observe that
	\begin{align*}
	E_{\bar{P}^{\nu}_{\ast}} \left[  \Delta_{k,k-T}(\theta_{\ast}) \mid X^{k-1}_{-\infty}    \right] =&	E_{\bar{P}^{\nu}_{\ast}} \left[  \Delta_{k,k-T}(\theta_{\ast}) \mid X^{k-1}_{k-T}    \right] \\
	= & \int   \frac{  \nabla_{\theta} p^{\nu}_{k}(x_{k} \mid X^{k-1}_{k-T}; \theta^{\ast})   }{  p^{\nu}_{k}(x_{k} \mid X^{k-1}_{k-T}; \theta^{\ast}) } p^{\nu}_{\ast}(x_{k} \mid X^{k-1}_{k-T}) dx_{k}.
	\end{align*}
	By assumption, $p^{\nu}_{k}(x_{k} \mid X^{k-1}_{k-T}; \theta^{\ast})  = p^{\nu}_{\ast}(x_{k} \mid X^{k-1}_{k-T})$. Moreover, by using the representation of $\Delta _{k,k-T}(\cdot )$ in the Supplemental Material \ref{SM:LAR1} and Assumption \ref{ass:deriva-bdd}, it can be shown that $\nabla_{\theta} \log p^{\nu}_{k}(\cdot \mid X^{k-1}_{k-T}; \theta)$ is uniformly bounded for all $\theta \in B(\delta,\theta_{\ast})$ (where $\delta>0$ is as in Assumption \ref{ass:deriva-bdd}). This result allow us to use standard \textquotedblleft Fisher information equality\textquotedblright\ calculations and thus deduce that $E_{\bar{P}^{\nu}_{\ast}} \left[  \Delta_{k,k-T}(\theta_{\ast}) \mid X^{k-1}_{-\infty}    \right] = 0$.
	
	Next, we show that, for any $k \geq 0$, $A_{k} \equiv E_{\bar{P}^{\nu}_{\ast}} \left[  \Delta_{k,-\infty}(\theta_{\ast})  \mid X^{k-1}_{-\infty}    \right] = 0$. To do so, observe that, for any $k,T \geq 0$,
	\begin{align*}
	\left \Vert A_{k}  \right \Vert^{2}_{L^{2}(\bar{P}^{\nu}_{\ast})} = & 	E_{\bar{P}^{\nu}_{\ast} }  \left[ \left(  E_{\bar{P}^{\nu}_{\ast}} \left[  \Delta_{k,-\infty}(\theta_{\ast}) \mid X^{k-1}_{-\infty}   \right]  - E_{\bar{P}^{\nu}_{\ast}} \left[  \Delta_{k,k-T}(\theta_{\ast}) \mid X^{k-1}_{-\infty}  \right]  \right)^{2} \right] \\
	\leq & \left \Vert  \Delta_{k,-\infty}(\theta_{\ast})   -  \Delta_{k,k-T}(\theta_{\ast})  \right \Vert^{2}_{L^{2}(\bar{P}^{\nu}_{\ast})}\\
	= & \left \Vert  \Delta_{0,-\infty}(\theta_{\ast})   -  \Delta_{0,-T}(\theta_{\ast})  \right \Vert^{2}_{L^{2}(\bar{P}^{\nu}_{\ast})}\\
	\leq & C \left(  \max\{\sum_{j=[-T/2]}^{-1}%
	\varrho(j,-T),\sum_{j=-T}^{[-T/2]-1}\varrho(-1,j)\}\right),
	\end{align*}
	where the second line follows from the Jensen inequality, the third line follows from stationarity, and the fourth line follows from Lemma \ref{lem:score_approx}(i). Now, recall that, for any $j \geq k$, $\varrho (j,k)\equiv \left( E_{\bar{P}^{\nu}_{\ast}} \left[
	\prod_{i=k}^{j}(1-\underline{q}(X_{i}))^{\frac{2a}{1-a}}\right] \right)^{\frac{1-a}{2a}} = \left( E_{\bar{P}^{\nu}_{\ast}} \left[
	\prod_{i=0}^{j-k}(1-\underline{q}(X_{i}))^{\frac{2a}{1-a}}\right] \right)^{\frac{1-a}{2a}} = \varrho (j-k,0) $. Moreover, by Assumption \ref%
	{ass:q-sum-p}, $(\varrho (j,0))_{j}$ is $p$-summable with $p<2/3$, and thus $%
	\lim_{j\rightarrow \infty }\varrho (j,0)^{p}j=0$ (if not, then $\varrho
	(j,0)>c/j^{1/p}$ for some $c>0$ and all $j$ above certain point, and this
	violates the assumption). Hence, for sufficienly large $T$,
	\begin{align*}
	\left \Vert A_{k}  \right \Vert^{2}_{L^{2}(\bar{P}^{\nu}_{\ast})} \leq & C \left(  \max\{\sum_{j=[-T/2]}^{-1}%
	\varrho(j+T,0),\sum_{j=-T}^{[-T/2]-1}\varrho(-1-j,0)\}\right) \\
	= & C \left(  \max\{\sum_{l=[T/2]}^{T-1} \varrho(l,0),\sum_{l=[T/2]}^{T-1}\varrho(l,0)\}\right)\\
	\leq&  C \sum_{l=[T/2]}^{T-1} l^{-(1/p)}.	
\end{align*}	
	As $1/p > 1$, it follows that, as $T$ diverges, the RHS converges to zero, and thus $\left \Vert A_{k}  \right \Vert^{2}_{L^{2}(\bar{P}^{\nu}_{\ast})} =0$. Since $\Delta_{k}(\theta_{\ast}) = \Delta_{k,-\infty}(\theta_{\ast})$ for any $k$, the desired result follows.
	
	By Theorem \ref{thm:anormal} and the central limit theorem for MDS, we have, therefore,
	\begin{equation*}
	\sqrt{T}(\hat{\theta}_{\nu,T}-\theta_{\ast})\Rightarrow_{\bar
		{P}_{\ast}^{\nu}}\mathcal{N}%
	(0,(E_{\bar{P}_{\ast}^{\nu}}[\xi_{1}(\theta_{\ast})])^{-1}\Sigma(\theta_{\ast
	})(E_{\bar{P}_{\ast}^{\nu}}[\xi_{1}(\theta_{\ast})])^{-1}),
	\end{equation*}	
	where $\Sigma (\theta _{\ast })\equiv \lim_{T\rightarrow \infty }\Sigma_{T}(\theta _{\ast })$, with $\Sigma (\theta _{\ast }) = E_{\bar{P}_{\ast}^{\nu}}\left[  \left( \Delta_{1}(\theta_{\ast} )  \right)  \left( \Delta_{1}(\theta_{\ast}) \right)  ^{\intercal}\right]$ since $(\Delta_{t}(\theta_{\ast}))_{t}$ is a stationary MDS. Moreover, it can be shown that $E_{\bar{P}_{\ast}^{\nu}}[\xi_{1}(\theta_{\ast})] = -E_{\bar{P}_{\ast}^{\nu}}\left[  \left( \Delta_{1}(\theta_{\ast}) \right)  \left( \Delta_{1}(\theta_{\ast}) \right)  ^{\intercal}\right]$. This follows by standard \textquotedblleft Fisher information equality\textquotedblright\  calculations and derivations analogous to those above (so they are omitted). Hence, the proof of part (a) is complete once we show that $-H^{-1}_{T}(\hat{\theta}_{\nu,T})$ converges in $\bar{P}_{\ast }^{\nu }$-probability to $(E_{\bar{P}_{\ast}^{\nu}}[\xi_{1}(\theta_{\ast})])^{-1}\Sigma(\theta_{\ast
	})(E_{\bar{P}_{\ast}^{\nu}}[\xi_{1}(\theta_{\ast})])^{-1}$. But this follows from Theorem \ref{thm:StdErrors}(a).

	\medskip

	Under the conditions of part (b), and in view of Lemma \ref{lem:suff.HAC.Rate}(a), $\Delta _{k}(\theta _{\ast })=\nabla
	_{\theta }\log p_{k}^{\nu }(X_{k}|X_{k-\bar{L}}^{k-1};\theta _{\ast })$ for all $k\geq 0$.  This result and the fact that $(X_{k})_{k=-\infty}^{\infty}$ is $\beta$-mixing with mixing coefficients $\beta_{n} = O(\gamma^{n})$ (see Lemma \ref{lem:sta-ergo}) imply that $(\Delta_{k}(\theta_{\ast}))_{k=-\infty}^{\infty}$ is also $\beta$-mixing with mixing coefficients of the same order. Hence, to establish that $T^{-1/2} \sum_{t=0}^{T} \Delta_{t}(\theta_{\ast}) \Rightarrow_{\bar
		{P}_{\ast}^{\nu}}\mathcal{N}(0,\Sigma(\theta_{\ast}))$, it is enough to verify that, for some $\delta>0$, $E_{\bar{P}_{\ast}^{\nu}}[||\Delta_{1}(\theta_{\ast})||^{2+\delta}]<\infty$ and $\sum_{n=1}^{\infty }\alpha _{n}^{\delta /(2+\delta )}<\infty $, where $(\alpha _{n})_{n}$ are the $\alpha$-mixing coefficients of $(\Delta_{k}(\theta_{\ast}))_{k=-\infty}^{\infty}$ (see, e.g., \cite{doukhan1994invariance}). But, the summability condition on the $\alpha$-mixing coefficients is satisfied, because the $\beta$-mixing coefficients of the process decay at rate  $O(\gamma^{n})$, and the moment condition is directly assumed.

	The proof of part (b) is completed by showing that $\hat{\Omega}_{T}(\hat{\theta}_{\nu ,T})$ converges in $\bar{P}_{\ast }^{\nu }$-probability to $(E_{\bar{P}_{\ast}^{\nu}}[\xi_{1}(\theta_{\ast})])^{-1}\Sigma(\theta_{\ast
	})(E_{\bar{P}_{\ast}^{\nu}}[\xi_{1}(\theta_{\ast})])^{-1}$, which is non-singular by assumption. We do so by invoking Theorem \ref{thm:StdErrors}(b) and verifying its conditions. As stated above, $\alpha_{n}= O( \gamma^{n} )$. Moreover, by Corollary 6.17 in \cite{white2001}, there exists $C<\infty$ such that, for any $l>0$,
	\begin{align*}
		|| E_{\bar{P}_{\ast}^{\nu} }[ \Delta_{l}(\theta_{\ast}) \Delta_{0}(\theta_{\ast})^{\intercal} ] || \leq C ( \alpha_{l} )^{\frac{2}{2+2\delta}} \sqrt{ E_{\bar{P}_{\ast}^{\nu} }[ ||\Delta_{l}(\theta_{\ast})||^{2}   ] }
		(E_{\bar{P}_{\ast}^{\nu} }[ ||\Delta_{l}(\theta_{\ast})||^{2+2\delta}   ])^{\frac{1}{2+2\delta}}.
	\end{align*}
Under our assumptions, the RHS equals $C ( \alpha_{l} )^{\frac{2}{2+2\delta}} $ for some $C<\infty $. Hence, one can set $\bar{\upsilon}(l) \equiv C  ( \gamma ^{\frac{2}{2+2\delta}} )^{l} $. Since $\gamma < 1$, the function $l \mapsto \bar{\upsilon}(l) $ is integrable. Also, by Lemma \ref{lem:suff.HAC.Rate}(b),  it follows that $r_{T} = L_{T} T^{-1/2}$.
	
	Next, we show that $\ddot{\varpi}  (\delta') = C \delta'$ for some $C<\infty $ and any $\delta' \leq \delta$, where $\delta>0$ is as in Assumption  \ref{ass:deriva-bdd}. To do so, we note that for any $\theta$ and any $t \in \mathbb{N}$,
	\begin{align*}
		||\Delta_{t}(\theta) \Delta_{0}(\theta)^{\intercal} - \Delta_{t}(\theta_{\ast} ) \Delta_{0}(\theta_{\ast})^{\intercal}|| \leq & ||\Delta_{0}(\theta) || \times || \Delta_{t}(\theta) - \Delta_{t}(\theta_{\ast})|| \\
		& + ||\Delta_{t}(\theta_{\ast}) || \times || \Delta_{0}(\theta) - \Delta_{0}(\theta_{\ast})||.
	\end{align*}

	\noindent By the calculations in the proof of Lemma \ref{lem:score_approx} in the Supplemental Material \ref{SM:LAR}, $||\Delta_{t}(\theta) ||$ is a linear combination of terms involving $\Gamma(\cdot,\theta)$ and $\Lambda(\cdot,\theta)$. Since $\Delta_{t}(\theta)$ only depends on $X^{t}_{t-\bar{L}}$ for some finite $\bar{L}$, there are only finitely many terms in this linear combination. Moreover, by Assumption \ref{ass:deriva-bdd}, there exists $C<\infty $ such that $E_{\bar{P}_{\ast}^{\nu} }[\sup_{\theta \in B(\delta, \theta_{\ast})} ||\Delta_{t}(\theta) ||^{2}] \leq C$, where $\delta>0$ is as in Assumption  \ref{ass:deriva-bdd}. This bound is also uniform over $t$. Hence, on account of this result and stationarity, there exists $C<\infty $ such that
\begin{align*}
E_{\bar{P}_{\ast }^{\nu }}\left[ \sup_{||\theta -\theta _{\ast }||\leq
	\delta ^{\prime }}\left\Vert \Delta _{t}(\theta )\Delta _{0}(\theta
)^{\intercal }-\Delta _{t}(\theta _{\ast })\Delta _{0}(\theta _{\ast
})^{\intercal }\right\Vert \right] \leq C\left( E_{\bar{P}_{\ast }^{\nu }}%
\left[ \sup_{||\theta -\theta _{\ast }||\leq \delta ^{\prime }}\left\Vert
\Delta _{0}(\theta )-\Delta _{0}(\theta _{\ast })\right\Vert ^{2}\right]
\right) ^{1/2},
\end{align*}
for any $\delta' \leq \delta$ and any $t \in \mathbb{N}$. It remains to show that the RHS is of the form $C \delta'$ for some finite constant $C$. By the mean value theorem, $\left\Vert	\Delta _{0}(\theta )-\Delta _{0}(\theta _{\ast })\right\Vert \leq \sup_{ \theta \in B( \delta,  \theta_{\ast}) }  || \nabla_{\theta} \Delta_{0}(\theta)|| \times ||\theta - \theta_{\ast}||$, so it suffices to show that $ E_{\bar{P}_{\ast}^{\nu} } [ \sup_{ \theta \in B( \delta,  \theta_{\ast}) }  || \nabla_{\theta} \Delta_{0}(\theta)||^{2} ] \leq C $ for some $ C<\infty $. Since $\Delta_{0}(\theta)= \nabla_{\theta} \log p_{0}^{\nu}(X_{0} \mid X^{-1}_{-\bar{L}})$, this condition is equivalent to $ E_{\bar{P}_{\ast}^{\nu} } [ \sup_{ \theta \in B( \delta,  \theta_{\ast}) }  || \nabla^{2}_{\theta} \log p^{\nu}_{0} (X_{0} \mid X^{-1}_{-\bar{L}} , \theta )  ||^{2} ]  \leq C$. By the calculations in \cite[pp. 1627--1628]{bickel98}, $ \nabla^{2}_{\theta} \log p^{\nu}_{0} (X_{0} \mid X^{0}_{-\bar{L}} , \theta )$ is a linear combination of finitely many terms involving $\Gamma(\cdot,\theta)$, $\Lambda(\cdot,\theta)$, and their derivatives. By Assumption \ref{ass:deriva-bdd} and calculations analogous to those in the proof of Lemma \ref{lem:score_approx}, it then follows that $ E_{\bar{P}_{\ast}^{\nu} } [ \sup_{ \theta \in B( \delta,  \theta_{\ast}) }  || \nabla^{2}_{\theta} \log p^{\nu}_{0} (X_{0} \mid X^{-1}_{-\bar{L}} , \theta )  ||^{2} ]  \leq C$. Thus, there exists $C<\infty$ such that $\ddot{\varpi}(\delta') \leq C \delta'$ for any $\delta' \leq \delta$. Therefore, all the conditions of Theorem \ref{thm:StdErrors}(b) are satisfied.

\end{proof}

%\section{Summary}
%\label{sec:summary}
%In this paper we have considered ML estimation in a large class of models
%which includes, among others, hidden Markov models and autoregressive models
%with Markov regimes. Our results extend earlier work by allowing for:
%(i)~autoregressive dynamics in the observable process; (ii)~time
%heterogeneity in the hidden Markov transition mechanism; (iii)~possible
%model misspecification. None of the existing papers in the literature allow
%for more than two of these features simultaneously. We have established
%asymptotic results related to consistency and LAN behavior of the ML
%estimator.\ In a Monte Carlo study, we have investigated the finite-sample
%properties of the ML\ estimator in a Markov-switching autoregressive model
%in which the variables that determine the evolution of the covariate-dependent
%transition probabilities may not be exogenous. We have also discussed an
%application involving real-world data.

%\bibliographystyle{econometrica}
{\small {
\bibliographystyle{plainnat}
\bibliography{markov}
%\printbibliography
}}

\clearpage\newpage

%\appendix

\newpage

\setcounter{page}{1}

\renewcommand\thesection{SM.\arabic{section}}
\renewcommand\thesubsection{\thesection.\arabic{subsection}} \setcounter{section}{0}

\begin{center}
{\Huge {Online Supplemental Material} }
\end{center}

For any measure $P$, we use $L^{r}(P)$, $1\leq r<\infty$, to denote the class
of measurable functions integrable to order $r$ with respect to $P$;
$\left\Vert \cdot\right\Vert _{L^{r}(P)}$ denotes the usual $r$-norm in
$L^{r}(P)$. For a vector/matrix-valued functions $X\mapsto f(X)$,
$||f||_{L^{r}(P)}$ is short-hand notation for the $L^{r}(P)$-norm of
$x\mapsto||f(x)||$, where $||\cdot||$ denotes the Euclidean/dual norm of $f$.
For any two sequences of random variables $(X_{n})_{n}$ and $(Y_{n})_{n}$,
$X_{n}\precsim Y_{n}$, implies that $X_{n}\leq CY_{n}$ for some universal
positive finite constant $C$.

\section{Ergodicity and Stationarity}

\label{SM:Ergo}

Let $(\zeta_{t})_{t=-\infty}^{\infty}$ be a Markov chain with transition
kernel $\zeta\mapsto\mathbb{P}(\zeta, \cdot) \in\mathcal{P} (\mathbb{Z})$ and
$\zeta_{t}\in\mathbb{Z}\subseteq\mathbb{R}^{d}$ for some $d>0$. Also, for any
probability measure $P$ over $\mathbb{Z}$ and any $f:\mathbb{Z}\rightarrow
\mathbb{R}$, let $P[f]\equiv\int f(z)P(dz)$ and $P[f](z) \equiv\int f(u) P(z,
du )$ (if it exists).

\begin{assumption}
\label{ass:Harris-1} \label{ass:Harris-2} There exist constants $\gamma
\in(0,1)$, $\lambda\in(0,1)$, $b>0$ and $R>2b/(1-\gamma)$, a function
$\mathbf{V}:\mathbb{Z}\rightarrow\lbrack1,\infty)$, and a probability measure
$\varrho$ such that: (i)$~\mathbb{P}[\mathbf{V}](\zeta)\leq\gamma
\mathbf{V}(\zeta)+b1\{\zeta\in\mathcal{C}\}$ for all $\zeta\in\mathbb{Z} $
with $\mathcal{C}\equiv\{\zeta\in\mathbb{Z}:\mathbf{V}(\zeta)\leq R\}$ ;
(ii)~$\inf_{\zeta\in\mathcal{C}}\mathbb{P}(\zeta, \cdot)\geq\lambda
\varrho(\cdot)$, with $\varrho(\mathcal{C})>0$.
\end{assumption}

%\begin{assumption}
%	\label{ass:Harris-2}  There exists a constant $\lambda \in [0,1)$ and a
%	probability measure $\varrho$ such that
%	\begin{align}  \label{eq:2}
%	\inf_{\zeta \in \mathcal{C}(R)} \mathbb{P}(\cdot|\zeta) \geq \lambda
%	\varrho(\cdot)
%	\end{align}
%	with some $R > 2K/(1-\gamma)$, and $\varrho(\mathcal{C}(R)) > 0$.
%\end{assumption}

The next result is used for the proof of Lemma \ref{lem:sta-ergo}; it contains
well-known results that are stated and proved here for convenience. In
particular, the first part of Lemma \ref{lem:Harris} is a re-statement of
Theorem 1.2 in \citesupp{HairerMatt11}. The second part of
Lemma~\ref{lem:Harris} and Assumption \ref{ass:Harris-1}(ii) imply that
$\mathbb{P}$ is Harris recurrent (see \citesupp[Ch. 14]{athreya2006measure})
and aperiodic (see \citesupp[p. 65]{Thierney1996}). The proof follows from
standard arguments.

Let $v\mapsto||v||_{\mathbf{V}}\equiv\sup_{\zeta}\frac{|v(\zeta)|}{1+
\mathbf{V}(\zeta)}$. Also, for any $A\subseteq\mathbb{Z}$, let $T_{A}%
=\inf\{t\geq0\colon\zeta_{t}\in A\}$.

\begin{lemma}
\label{lem:Harris} If Assumption \ref{ass:Harris-1} holds, then:

(i)$~\mathbb{P}$ admits a unique invariant measure $\nu^{\ast}$, and there
exist constants $\gamma\in(0,1)$ and $C>0$ such that
\[
||\mathbb{P}^{n}[v]-\nu^{\ast}[v]||_{\mathbf{V}}\leq C\gamma^{n}||v-\nu^{\ast
}[v]||_{\mathbf{V}}%
\]
for every measurable function $v$ such that $||v||_{\mathbf{V}}<\infty$, where
$\nu^{\ast}[v]\equiv\int v(\zeta)\nu^{\ast}(d\zeta)$.

(ii)$~\mathbb{P}(\zeta, \{ T_{\mathcal{C}}<\infty\} )=1$ for all $\zeta
\in\mathbb{Z}$, and $\mathbb{P}(\zeta_{0}, \mathcal{C} )>0$ for all $\zeta_{0}
\in\mathcal{C}$.
\end{lemma}

\begin{proof}[Proof of Lemma \ref{lem:Harris}]
	Part~(i) is Theorem 1.2 in \citesupp{HairerMatt11}. Assumption \ref{ass:Harris-1}%
	(i) implies their Assumption 1 with $K=b$ and Assumption \ref{ass:Harris-1}%
	(ii) implies their Assumption 2.
	
	For part (ii), we first establish that $\mathbb{P}(\zeta _{0}, \mathcal{C} ) = \mathbb{P}(\zeta _{0},\{ \zeta _{1}\in \mathcal{C} \})>0$ for all $\zeta _{0}\in \mathcal{C}$. For this, note that $ \mathbb{P}( \zeta _{0} , \mathcal{C} ) \geq \inf_{\zeta \in \mathcal{C}} \mathbb{P}(\zeta , 	\mathcal{C}  ) \geq \lambda \varrho (C)>0$  by Assumption \ref{ass:Harris-1}(ii).
	.
	
	We now show that $\mathbb{P}(\zeta, \{ T_{\mathcal{C}}<\infty \} )=1$ for all $%
	\zeta \in \mathbb{Z}$. It suffices to show that $\mathbb{P}[T_{\mathcal{C}%
	}](\zeta )<\infty $ for all $\zeta \notin \mathcal{C}$. Under Assumption \ref%
	{ass:Harris-1}(i), $\mathbf{V}\geq 1$, so
	\begin{equation*}
	\mathbb{P}[T_{\mathcal{C}}](\zeta )\leq \mathbb{P}[\sum_{j=0}^{T_{\mathcal{C}%
		}-1}\mathbf{V}(\zeta _{j})](\zeta )=\sum_{T=0}^{\infty }\sum_{j=0}^{T}E_{%
		\mathbb{P}}[\mathbf{V}(\zeta _{j})\mid T_{\mathcal{C}}=T+1,\zeta ]\Pr (\zeta, \{T_{%
		\mathcal{C}}=T+1\} )
	\end{equation*}%
	for any $\zeta \in \mathbb{Z}\setminus \mathcal{C}$. To establish the
	desired result, it is sufficient to show that $\sup_{T}\sum_{j=0}^{T}E_{%
		\mathbb{P}}[\mathbf{V}(\zeta _{j})\mid T_{\mathcal{C}}=T+1,\zeta ]<\infty $.
	
	Take any $T\geq 0$ and any $j\leq T$, and note that
	\begin{align*}
	E_{\mathbb{P}}\left[ \mathbb{P}[\mathbf{V}](\zeta _{j})\mid \zeta _{l}\notin
	\mathcal{C},~\forall l\leq T+1\right] =& E_{\mathbb{P}}\left[ \int_{\zeta
		\notin \mathcal{C}}\mathbb{P}[\mathbf{V}](\zeta )\mathbb{P}(
	\zeta _{j-1}, d\zeta)\mid \zeta _{l}\notin \mathcal{C},~\forall l\leq j-1\right] \\
	\leq & \gamma E_{\mathbb{P}}\left[ \int_{\zeta \notin \mathcal{C}}\mathbf{V}%
	(\zeta )\mathbb{P}(\zeta _{j-1},d\zeta )\mid \zeta _{l}\notin \mathcal{C}%
	,~\forall l\leq j-1\right] \\
	\leq & \gamma E_{\mathbb{P}}\left[ \mathbb{P}[\mathbf{V}](\zeta _{j-1})\mid
	\zeta _{l}\notin \mathcal{C},~\forall l\leq j-1\right] \\
	\leq & \gamma ^{j}\mathbf{V}(\zeta _{0}),
	\end{align*}%
	where the second line follows from Assumption \ref{ass:Harris-1}(i) and the
	fact that $\zeta \notin \mathcal{C}$, the third line follows from the fact
	that $\mathbf{V}>0$, and the last line follows from repeated iteration of
	the first lines. Note that $T_{\mathcal{C}}=T+1$ is equivalent to $\zeta
	_{j}\notin \mathcal{C},~\forall j\leq T$ and $\zeta _{T+1}\in \mathcal{C}$.
	Thus, the previous display implies that
	\begin{equation*}
	E_{\mathbb{P}}\left[ \mathbf{V}(\zeta _{j})\mid T_{\mathcal{C}}=T+1\right]
	\leq \gamma ^{j}\mathbf{V}(\zeta _{0})
	\end{equation*}%
	for any $T\geq 0$ and any $j\leq T$. Consequently, $\sum_{j=0}^{T}E_{\mathbb{%
			P}}[\mathbf{V}(\zeta _{j})\mid T_{\mathcal{C}}=T+1,\zeta ]\leq \mathbf{V}%
	(\zeta )\sum_{j=0}^{T}\gamma ^{j}\leq \frac{\mathbf{V}(\zeta )}{1-\gamma }$,
	and thus the result follows.
\end{proof}

\begin{proof}[Proof of Lemma \protect\ref{lem:sta-ergo}]
	Let $(\zeta _{t})_{t=-\infty}^{\infty}$ be the stochastic process given by $\zeta _{t}\equiv
	(X_{t},S_{t})$. This process is a Markov chain with transition kernel $%
	\mathbb{X}\times \mathbb{S}\ni \zeta \mapsto \mathbb{P}(\zeta ,\cdot)\in
	\mathcal{P}(\mathbb{X}\times \mathbb{S})$ given by
	\begin{equation*}
	\mathbb{P}((x,s),\{ \zeta _{t+1}\in A_{1}\times A_{2}\} )=\sum_{s^{\prime }\in A_{2}}Q_{\ast}(x,s,s^{\prime})P_{\ast}( x,s^{\prime },A_{1}),
	\end{equation*}%
	for any Borel sets $A_{1}\subseteq \mathbb{X}$ and $A_{2}\subseteq \mathbb{S}
	$.
	
	By Lemma~\ref{lem:Harris}, there exists a
	unique invariant measure $\nu $, provided that the conditions of Assumption %
	\ref{ass:Harris-2} are met.	In order to verify the first part of Assumption \ref{ass:Harris-2}, consider
	$\mathbf{V}(\zeta )=\mathcal{U}(x)$, and $\mathcal{C}\equiv \mathcal{C}%
	_{1}\times \mathbb{S}$ with $\mathcal{C}_{1}\equiv \{x\in \mathbb{X}\colon
	\mathcal{U}(x)\leq R\}$. By Assumption \ref{ass:fx-ergo}(i),
	\begin{equation*}
	\mathbb{P}[\mathbf{V}](\zeta )=\int_{\mathbb{X}}\mathcal{U}(x^{\prime
	})\left\{ \sum_{s^{\prime }\in \mathbb{S}}Q_{\ast}(x,s,s^{\prime})P_{\ast}(x,s^{\prime },dx^{\prime })\right\} \leq \gamma
\mathcal{U}(x)+2b^{\prime }1\{x\in \mathcal{C}_{1}\}.
\end{equation*}%
Thus, $b\equiv 2b^{\prime }$.
Regarding Assumption \ref{ass:Harris-2}(ii), observe that, by Assumption \ref{ass:BDD_Q}(i), for $C$ and any $s\in \mathbb{S}$,
\begin{equation*}
\mathbb{P}((x,s) , C\times \{s^{\prime}\} )\geq \underline{q}(x)P_{\ast}(x,s^{\prime },C),
\end{equation*}%
and, by Assumption \ref{ass:fx-ergo}(iii), $P_{\ast
}(x,s^{\prime },C)\geq \lambda ^{\prime }\varpi (C)$ and $\lambda ^{\prime
}\in (0,1)$. Also note that, by Assumption~\ref{ass:BDD_Q}, $\underline{q}$
is continuous and $\underline{q}(x)>0$ for all $x\in \mathbb{X}$.
Furthermore, by Assumption~\ref{ass:fx-ergo}(ii), $\mathcal{U}$ is lower
semi-compact, because $\{x\in \mathbb{X}:\mathcal{U}(x)\leq R\}$ is closed ($%
x\mapsto \mathcal{U}(x)$ is lower semi-continuous), and is also bounded.
Therefore, $\inf_{x:\mathcal{U}(x)\leq R}\underline{q}(x)=\min_{x:\mathcal{U}%
	(x)\leq R}\underline{q}(x)\geq c>0$ (because it is a minimization of a
continuous function on compact set). Therefore,
\begin{equation*}
\mathbb{P}(\zeta , C\times \{s^{\prime}\}  )\geq c\lambda ^{\prime }\varpi (C)\frac{1}{|\mathbb{S}|},
\end{equation*}%
and, by putting $\varrho = \varpi (\cdot)\frac{1}{|\mathbb{S}|}$ and $\lambda \equiv
c\lambda ^{\prime }$, Assumption \ref{ass:Harris-2}(ii) follows since $%
\varpi (\mathcal{C}_{1})>0$.
Since $\nu $ is unique, it is trivially ergodic. Therefore, the process with
initial probability measure $\nu $ is stationary. Ergodicity of $(\zeta
_{t})_{t}$ follows from Theorem~14.2.11 in \citesupp{athreya2006measure} (recall
that $\mathbb{P}$ is Harris recurrent and aperiodic).
%XXX $\sup_{x\in \mathcal{C}_{1}}\int \mathcal{U}%
%(x^{\prime })P(dx^{\prime }\mid s,x)\leq \gamma R+b^{\prime }<\infty $ XXX;
%see XXX CITE XXX.
Since $X_{t}$ is a deterministic function of $\zeta _{t}$, $X_{0}^{\infty }$
is also stationary and ergodic.
Finally, observe that
\begin{equation*}
\int \sup_{0\leq f\leq 1}\left\vert \mathbb{P}^{n}[f](\zeta )-\nu \lbrack
f]\right\vert \nu (d\zeta ) \precsim   \gamma ^{n}\int \left\vert 1+\mathcal{U}%
(x)\right\vert \nu (d\zeta ).
\end{equation*}%
Since $\mathcal{U}$ satisfies Assumption \ref{ass:Harris-1}(i), it
follows that $\int \mathbb{P}[\mathcal{U}](\zeta )\nu (d\zeta )\leq \gamma
\nu \lbrack \mathcal{U}]+K$. Since $\nu $ is the invariant measure of $%
\mathbb{P}$ and $\gamma \in (0,1)$, this implies that $\nu (d\zeta )\leq
K/(1-\gamma )$. Therefore,
\begin{equation*}
\int \sup_{0\leq f\leq 1}\left\vert \mathbb{P}^{n}[f](\zeta )-\nu \lbrack
f]\right\vert \nu (d\zeta ) \precsim   \gamma ^{n} ,
\end{equation*}%
thereby implying that $(\zeta _{t})_{t}$ is $\beta $-mixing with rate $\beta
_{n}=O(\gamma ^{n})$ (see \citesupp{davydov1974mixing}). Since $X_{t}$ is a
deterministic function of $\zeta _{t}$, the same holds for $X_{0}^{\infty }$.
\end{proof}

\section{Proofs of Supplementary Lemmas in Appendix \ref{app:consistent}}

\label{SM:consistent}

To prove Lemmas \ref{lem:lik-approx} and \ref{lem:lik-conv} we use the
following result.

\begin{lemma}
\label{lem:approx-l} Suppose Assumptions \ref{ass:BDD_Q} and \ref{ass:bdd}(ii)
hold. Then, for all $t\in\mathbb{N}$ and $-n\leq-m\leq t-1$,
\[
\sup_{\theta\in\Theta}\left\vert \log p_{t}^{\nu}(X_{t}\mid X_{-m}%
^{t-1},\theta)-\log p_{t}^{\nu}(X_{t}\mid X_{-n}^{t-1},\theta)\right\vert \leq
C(X_{t-1},X_{t}) \prod_{i=-m}^{t-1}(1-\underline{q}(X_{i}))
\]
a.s.-$\bar{P}_{\ast}^{\nu}$.
\end{lemma}

\begin{proof}[Proof of Lemma \protect\ref{lem:approx-l}]
	Observe that, for any $n\in \mathbb{N}$,
	\begin{equation*}
	\log p_{t}^{\nu }(X_{t}\mid X_{-n}^{t-1},\theta )=\log \sum_{s\in \mathbb{S}%
	}p_{\theta}(X_{t-1},s,X_{t}) \bar{P}_{\theta}^{\nu} (s\mid X_{-n}^{t-1} ),
	\end{equation*}%
	and since $\log x-\log y\leq x/y-1$, it suffices to study $ \frac{\sum_{s\in \mathbb{S}}p_{\theta}(X_{t-1},s,X_{t})\left( \bar{P}_{\theta}^{\nu}
		(S_{t} = s\mid X_{-m}^{t-1} )-\bar{P}_{\theta}^{\nu} (S_{t} = s\mid X_{-n}^{t-1})\right) }{%
		\sum_{s\in \mathbb{S}}p_{\theta}(X_{t-1},s,X_{t}) \bar{P}_{\theta}^{\nu} (s\mid
		X_{-n}^{t-1})}$.
%	\begin{align*}
%	& \frac{\sum_{s\in \mathbb{S}}p_{\theta}(X_{t}\mid X_{t-1},s)\Pr (s\mid
%		X_{-m}^{t-1},\theta )-\sum_{s\in \mathbb{S}}p_{\theta}(X_{t}\mid
%		X_{t-1},s)\Pr (s\mid X_{-n}^{t-1},\theta )}{\sum_{s\in \mathbb{S}}f_{\theta
%		}(X_{t}\mid X_{t-1},s)\Pr (s\mid X_{-m-1}^{t-1},\theta )} \\
%	=& \frac{\sum_{s\in \mathbb{S}}p_{\theta}(X_{t}\mid X_{t-1},s)\left( \Pr
%		(s\mid X_{-m}^{t-1},\theta )-\Pr (s\mid X_{-n}^{t-1},\theta )\right) }{%
%		\sum_{s\in \mathbb{S}}p_{\theta}(X_{t}\mid X_{t-1},s)\Pr (s\mid
%		X_{-n}^{t-1},\theta )}.
%	\end{align*}
	This expression can be bounded above by
	\begin{equation*}
	\frac{\max_{s\in \mathbb{S}}p_{\theta}(X_{t-1},s,X_{t})}{\min_{s\in \mathbb{%
				S}}p_{\theta}(X_{t-1},s,X_{t})} \left\Vert
	\bar{P}_{\theta}^{\nu}  (S_{t}= \cdot \mid X_{-m}^{t-1} )- \bar{P}_{\theta}^{\nu}  (S_{t}= \cdot \mid X_{-n}^{t-1} )\right\Vert_{1}.
	\end{equation*}
	
	By Assumption \ref{ass:bdd}(ii),  $\sup_{\theta \in \Theta }\frac{\max_{s\in \mathbb{S}}p_{\theta
		}(X_{t-1},s,X_{t})}{\min_{s\in \mathbb{S}}p_{\theta}(X_{t-1},s,X_{t})}\leq
	C(X_{t-1},X_{t})$ a.s.-$\bar{P}_{\ast}^{\nu }$. So it suffices to bound $\left\Vert \bar{P}_{\theta}^{\nu}(S_{t} = \cdot \mid X_{-m}^{t-1} )- \bar{P}_{\theta}^{\nu}(S_{t} = \cdot \mid X_{-n}^{t-1} )\right\Vert_{1}$.
		
	By Lemma B.2.2 in \citesupp{stachurski2009}, 	
\begin{align*}
	& \left\Vert	\bar{P}_{\theta}^{\nu} (S_{t}= \cdot \mid X_{-m}^{t-1})-\bar{P}_{\theta}^{\nu} (S_{t}= \cdot \mid X_{-n}^{t-1}	)\right\Vert_{1} \\
	\leq &  \frac{1}{2} \sup_{b,c \in \mathbb{S}^{2} }	\left\Vert	\bar{P}_{\theta}^{\nu} (S_{t}= \cdot \mid S_{-m} = b , X_{-m}^{t-1})- \bar{P}_{\theta}^{\nu} (S_{t}= \cdot \mid S_{-m} = c , X_{-n}^{t-1}
)\right\Vert_{1} \\
	= & \frac{1}{2} \sup_{b,c \in \mathbb{S}^{2} }	\left\Vert	\bar{P}_{\theta}^{\nu} (S_{t}= \cdot \mid S_{-m} = b , X_{-m}^{t-1} )- \bar{P}_{\theta}^{\nu} (S_{t}= \cdot \mid S_{-m} = c , X_{-m}^{t-1}
	)\right\Vert_{1},
\end{align*}		
%		Observe that for any $- n \leq - m \leq k$ and any $a \in \mathbb{S}$ (we
%		omit the dependence on $\theta$ to simplify the notation)
%		%	 	  	\begin{align*}
%		%	 	  	& \left| \Pr( S_{k+1} = a \mid X^{k}_{-m}) - \Pr( S_{k+1} = a \mid X^{k}_{-n} ) \right| \\
%		%	 	  	= & \left| \sum_{b \in \mathbb{S}} \{ \Pr( S_{k+1} = a \mid S_{k}=b, X^{k}_{-m} )\Pr( S_{k}=b \mid X^{k}_{-m} ) - \Pr( S_{k+1} = a \mid S_{k}=b,X^{k}_{-n} )\Pr( S_{k}=b \mid X^{k}_{-n} ) \} \right|\\
%		%	 	  	= & \left| \sum_{b \in \mathbb{S}} \Pr( S_{k+1} = a \mid S_{k}=b, X_{k} ) \{ \Pr( S_{k}=b \mid X^{k}_{-m} ) - \Pr( S_{k}=b \mid X^{k}_{-n} ) \} \right|\\
%		%	 	  	\leq & |\mathbb{S}| \max_{b \in \mathbb{S}} \left| \Pr( S_{k}=b \mid X^{k}_{-m} ) - \Pr( S_{k}=b \mid X^{k}_{-n} )  \right|.
%		%	 	  	\end{align*}
%		
%		{\small {\
%				\begin{align*}
%				& \left\vert \Pr (S_{k+1}=a\mid X_{-m}^{k})-\Pr (S_{k+1}=a\mid
%				X_{-n}^{k})\right\vert \\
%				=& \left\vert \sum_{b\in \mathbb{S}}\{\Pr (S_{k+1}=a\mid
%				S_{-m}=b,X_{-m}^{k})\Pr (S_{-m}=b\mid X_{-m}^{k})-\Pr (S_{k+1}=a\mid
%				S_{-m}=b,X_{-n}^{k})\Pr (S_{-m}=b\mid X_{-n}^{k})\}\right\vert \\
%				\leq & \left\vert \max_{b,b^{\prime }}\{\Pr (S_{k+1}=a\mid
%				S_{-m}=b,X_{-m}^{k})-\Pr (S_{k+1}=a\mid S_{-m}=b^{\prime
%				},X_{-n}^{k})\}\right\vert \\
%				=& \left\vert \max_{b,b^{\prime }}\{\Pr (S_{k+1}=a\mid
%				S_{-m}=b,X_{-m}^{k})-\Pr (S_{k+1}=a\mid S_{-m}=b^{\prime
%				},X_{-m}^{k})\}\right\vert ,
%				\end{align*}%
%			}}
where the last line follows from the fact that, given $S_{-m}$, it is the
			same to condition on $X_{-m}^{t-1}$ and on $X_{-n}^{t-1}$.	Hence,
			\begin{align*}
			 \sup_{\theta \in \Theta }\left\vert \log p_{t}^{\nu }(X_{t}\mid
			X_{-m}^{t-1},\theta )-\log p_{t}^{\nu }(X_{t}\mid X_{-n}^{t-1},\theta
			)\right\vert 	\leq  C^{\prime } \alpha_{\theta,t,-m}(X^{t-1}_{-m}),
			\end{align*}
			where $\alpha_{\theta,t,-m}(X^{t-1}_{-m})$ is defined in expression (\ref{eqn:Dobrushin}). By Applying Lemmas \ref{lem:ergo-eta} and \ref{lem:couple} and the fact that $\alpha_{\theta,t,-m}(X^{t-1}_{-m}) \leq \prod_{j=-m}^{t-1} \alpha_{\theta,j,j+1}(X^{t-1}_{-m})$, it follows that
			\begin{equation*}
			\sup_{\theta \in \Theta }\left\vert \log p_{t}^{\nu }(X_{t}\mid
			X_{-m}^{t-1},\theta )-\log p_{t}^{\nu }(X_{t}\mid X_{-n}^{t-1},\theta
			)\right\vert \leq C(X_{t-1},X_{t}) \prod_{i=-m}^{t-1}(1-\underline{q}	(X_{i})),
			\end{equation*}%
			a.s.-$\bar{P}_{\ast}^{\nu }$.
		\end{proof}

We now prove Lemmas \ref{lem:lik-approx} and \ref{lem:lik-conv}.

\begin{proof}[Proof of Lemma \protect\ref{lem:lik-approx}]
			Fix any $\epsilon >0$. By Lemma \ref{lem:approx-l}, with $m=0$ and $n=M$,
			\begin{align*}
			& \sup_{\theta \in \Theta} \left\vert \log p_{t}^{\nu }(X_{t}\mid X_{0}^{t-1},\theta )-\log p^{\nu
			}(X_{t}\mid X_{-\infty }^{t-1},\theta )\right\vert \\
		\leq & \limsup_{M \rightarrow \infty} \sup_{\theta \in \Theta} \left\vert \log p_{t}^{\nu }(X_{t}\mid X_{0}^{t-1},\theta )-\log
			p_{t}^{\nu }(X_{t}\mid X_{-M}^{t-1},\theta )\right\vert \\
			\leq & C(X_{t-1},X_{t}) \prod_{i=0}^{t-1} (1-\underline{q}(X_{i})).
			\end{align*}
			
			Thus, it suffices to show that there exists a $T$ such that for all $t \geq T$, $\bar{P}^{\ast}_{\nu} \left( T^{-1} \sum_{t=1}^{T} C(X_{t-1},X_{t}) \prod_{i=0}^{t-1} (1-\underline{q}(X_{i})) \geq \epsilon \right) < \epsilon$. This follows from Assumption \ref{ass:qbar-sum}.		
	\end{proof}

\begin{proof}[Proof of Lemma \protect\ref{lem:lik-conv}]
		Recall that, by Lemma \ref{lem:sta-ergo}, the process $X_{-\infty }^{\infty }$ is ergodic and stationary under $\bar{P}^{\nu}_{\ast}$.
		
		\textbf{Part (i).} Consider a $\delta >0$ and an open cover $\{B(\theta
		,\delta )\colon \theta \in \Theta \}$ where $B(\theta ,\delta )$ is an open
		ball centered around $\theta $ with radius $\delta >0$. Since $\Theta $ is
		compact (Assumption \ref{ass:exist}), there exists a finite sub-cover $%
		B_{j}\equiv B(\theta _{j},\delta )$ with $j=1,\ldots ,J$. Also note that,
		pointwise in $\theta \in \Theta $, $\ell _{T}^{\nu }(X_{-\infty }^{T},\theta
		)-E_{\bar{P}_{\ast}^{\nu }}[\ell _{T}^{\nu }(X_{-\infty
		}^{T},\theta )]\rightarrow 0$ a.s.-$\bar{P}_{\ast}^{\nu }$ by the
		ergodic theorem and the fact that $X_{-\infty }^{\infty }\mapsto \ell
		_{T}^{\nu }(X_{-\infty }^{T},\theta )\in L^{1}(\bar{P}_{\theta ^{\ast
			}}^{\nu })$. Thus, it suffices to show that there exists a $T(j,\epsilon )$
			such that, for all $t\geq T(j,\epsilon )$,
			\begin{equation*}
			\bar{P}_{\ast}^{\nu }\left( \sup_{\theta \in
				B_{j}}T^{-1}\sum_{t=1}^{T}\left( l_{t}(X_{-\infty }^{t},\theta )-E_{\bar{P}%
				_{\ast}^{\nu }}[l_{t}(X_{-\infty }^{t},\theta )]\right) >\epsilon
			\right) \leq \epsilon ,
			\end{equation*}%
			where $l_{t}(X_{-\infty }^{t},\theta )\equiv \log \frac{p^{\nu }(X_{t}\mid
				X_{-\infty }^{t-1},\theta)}{p^{\nu }(X_{t}\mid X_{-\infty
				}^{t-1},\theta_{j} )}$. Observe that, for any $j$,
			\begin{align*}
			\sup_{\theta \in B_{j}}\sum_{t=1}^{T}\left( l_{t}(X_{-\infty }^{t},\theta
			)-E_{\bar{P}_{\ast}^{\nu }}[l_{t}(X_{-\infty }^{t},\theta
			)]\right) \leq & \sum_{t=1}^{T}\sup_{\theta \in B_{j}}\left(
			l_{t}(X_{-\infty }^{t},\theta )-E_{\bar{P}_{\ast}^{\nu
				}}[l_{t}(X_{-\infty }^{t},\theta )]\right) \\
				\equiv & \sum_{t=1}^{T}\bar{l}_{t}(X_{-\infty }^{t}).
				\end{align*}%
				Moreover, observe that
				\begin{equation*}
				\sup_{\theta \in B_{j}}\log \frac{p^{\nu }(X_{t}\mid X_{-\infty
					}^{t-1},\theta)}{p^{\nu }(X_{t}\mid X_{-\infty }^{t-1},\theta_{j} )}\leq
				\sup_{\theta \in B_{j}}\frac{p^{\nu }(X_{t}\mid X_{-\infty }^{t-1},\theta)}{p^{\nu }(X_{t}\mid X_{-\infty }^{t-1},\theta_{j} )}-1.
				\end{equation*}
				
				By Assumption \ref{ass:bdd}(i), for any $\epsilon >0$ there exists a $\delta
				>0$ such that $E_{\bar{P}_{\ast}^{\nu }}\left[ \sup_{\theta \in B_{j}}\frac{p^{\nu }(X_{0}\mid
					X_{-\infty }^{-1},\theta)}{p^{\nu }(X_{0}\mid X_{-\infty }^{-1},\theta_{j})%
				}\right] \leq 1+\epsilon $ for any $j\in \{1,\ldots ,J\}$ and any $t$.
				Therefore, we can choose a $\delta >0$ such that
				\begin{equation*}
				E_{\bar{P}_{\ast}^{\nu }}\left[ \sup_{\theta \in B_{j}}\log \frac{%
					p^{\nu }(X_{0}\mid X_{-\infty }^{-1},\theta)}{p^{\nu }(X_{0}\mid
					X_{-\infty }^{-1},\theta_{j} )}\right] \leq \epsilon /4.
				\end{equation*}%
				This in turn implies that $E_{\bar{P}_{\ast}^{\nu }}[\bar{l}%
				_{t}(X_{-\infty }^{t})]\leq \epsilon /2$. This result and the ergodic
				theorem establish that $\lim_{T\rightarrow \infty }T^{-1}\sum_{t=1}^{T}\bar{l%
				}_{t}(X_{-\infty }^{t})\leq \epsilon /2$ a.s.-$\bar{P}_{\ast}^{\nu }$. This implies the result in (\ref{eqn:lik-conv-1}).
				
				\textbf{Part (ii).} Follows directly from the ergodic theorem and the fact
				that $X_{-\infty }^{\infty }\mapsto \log p^{\nu }(X_{t}\mid X_{-\infty
				}^{t-1},\theta_{\ast})$ is in $L^{1}(\bar{P}_{\ast}^{\nu })$.
			\end{proof}

\section{Properties of $p_{\theta}(X_{1}|X^{0}_{-\infty})$}

\label{SM:PDFs}

For any $t\in\mathbb{N}_{0}\equiv\mathbb{N}\cup\{0\}$, any $X^{t}_{-\infty}$
and any $\theta\in\Theta$, $p^{\nu}(X_{t}|X^{t-1}_{-\infty},\theta)$ is
defined as $\liminf_{M \rightarrow\infty} p_{\theta}^{\nu}(X_{t}|X^{t-1}%
_{-M})$; $p^{\nu}_{\ast}(X_{t}|X^{t-1}_{-\infty})$ is defined analogously.

\begin{lemma}
\label{lem:prop-ergoPDF} Suppose Assumptions \ref{ass:bdd}(ii) and
\ref{ass:qbar-sum} hold. Then:

(1) For any $t \in\mathbb{N}_{0}$, $x \mapsto p^{\nu}(\cdot|X^{t-1}_{-\infty
},\theta)$ and $x \mapsto p^{\nu}_{\ast}(\cdot|X^{t-1}_{-\infty})$ are
densities, a.s.-$\bar{P}^{\nu}_{\ast}$.

(2) Suppose $\Theta$ is compact and that for each $n \in\mathbb{N}_{0}$,
$\theta\mapsto p^{\nu}_{\theta}(X_{1}|X^{0}_{-n})$ is uniformly continuous
a.s.-$\bar{P}^{\nu}_{\ast}$. Then, for any $\theta_{0} \in\Theta$ and
$\epsilon>0$, there exists a $\delta>0$ such that
\begin{align*}
&  \bar{P}^{\nu}_{\ast} \left(  \sup_{\theta_{0} \in\Theta} \sup_{\theta\in
B(\theta_{0},\delta)} | p^{\nu}(X_{1}|X^{0}_{-\infty},\theta) - p^{\nu}%
(X_{1}|X^{0}_{-\infty},\theta_{0}) | > \epsilon\right)  < \epsilon.
\end{align*}

(3) Suppose $\Theta$ is compact and that for each $n \in\mathbb{N}_{0}$,
$\theta\mapsto p^{\nu}_{\theta}(X_{1}|X^{0}_{-n})$ is uniformly continuous
a.s.-$\bar{P}^{\nu}_{\ast}$. Suppose also that there exists functions
$(x_{1},x_{0}) \mapsto(\overline{p}(x_{0},x_{1}),\underline{p}(x_{0},x_{1}))$
such that for any $p \in\{ p_{\theta} : \theta\in\Theta\} \cup p_{\ast} $,
$\underline{p}(x_{0},x_{1}) \leq p(x_{0},s,x_{1}) \leq\overline{p}(x_{0}%
,x_{1})$ for all $s \in\mathbb{S}$, and $E_{\bar{P}^{\nu}_{\ast}}[\overline
{p}(X_{0},X_{1})/\underline{p}(X_{0},X_{1})] < \infty$. Then, Assumptions
\ref{ass:exist} and \ref{ass:bdd} hold.
\end{lemma}

\begin{proof}
(1) We need to show that the functions integrate to 1. By analogous steps to those in the proof of Lemma \ref{lem:approx-l}, it follows that
\begin{align*}
& |\int \{p^{\nu}(x|X^{t-1}_{-\infty},\theta)  -  p_{\theta}^{\nu}(x|X^{t-1}_{-n})\}  dx  | \\
\leq & \int \sum_{s} p_{\theta}(X_{t-1},s,x) dx \limsup_{M \rightarrow \infty } || \bar{P}^{\nu}_{\theta}(S_{t} = \cdot \mid X^{t-1}_{-M} ) - \bar{P}^{\nu}_{\theta}(S_{t} = \cdot \mid X^{t-1}_{-n} )   ||_{1}\\
\leq & \int \sum_{s} p_{\theta}(X_{t-1},s,x) dx \limsup_{M \rightarrow \infty }  \prod_{i=-n}^{t-1} (1-\underline{q}(X_{i}))      dx \\
= &    \prod_{i=-n}^{t-1} (1-\underline{q}(X_{i}))  |\mathbb{S}|.
\end{align*}
Since this holds for any $n$ such that $ - n \leq t-1$, we can take averages and obtain, for any $\epsilon>0$, that
\begin{align*}
\bar{P}^{\nu}_{\ast} \left( \frac{1}{M+1} \sum_{n=0}^{M} |\int p^{\nu}(x|X^{t-1}_{-\infty},\theta) dx - 1| >  \epsilon  \right) \leq & 	\bar{P}^{\nu}_{\ast} \left( \frac{1}{M+1} \sum_{n=0}^{M} \prod_{i=-n}^{t-1} (1-\underline{q}(X_{i}))  |\mathbb{S}|  >   \epsilon  \right)\\
\leq & 	\bar{P}^{\nu}_{\ast} \left( \frac{1}{M+1} \sum_{n=0}^{M} \prod_{i=-n}^{0} (1-\underline{q}(X_{i}))  |\mathbb{S}|  >   \epsilon  \right).
\end{align*}
Since this holds for any $M$, by taking $M \rightarrow \infty$, stationarity and Assumption \ref{ass:qbar-sum} imply that the RHS vanishes. Thus, it follows that $	\bar{P}^{\nu}_{\ast} \left(  | \int p^{\nu}(x|X^{t-1}_{-\infty},\theta) dx  - 1 | \geq  \epsilon  \right) =0$. As the $\epsilon>0$ is arbitrary, this implies that $\int p^{\nu}(x|X^{t-1}_{-\infty},\theta) dx  = 1$, a.s.-$\bar{P}^{\nu}_{\ast}$.
%
%Lemma \ref{lem:approx-l} below and the fact that $C(X_{t-1},X_{t}) < \infty$ and $\lim_{m \rightarrow \infty} \prod_{i=-m}^{t-1} (1- \underline{q}(X_{i})) = 0$ a.s.-$\bar{P}^{\nu}_{\ast}$ (by XXX) imply that the sequence of functions $(\theta \mapsto p^{\nu}_{\theta}(X_{t}|X^{t-1}_{-n}))_{n = 1}^{\infty}$ is Cauchy under $L^{\infty}(\Theta)$. By completeness, this implies that, a.s.-$\bar{P}^{\nu}_{\ast}$, the limit exists an is equal to $\theta \mapsto p^{\nu}_{\theta}(X_{t}|X^{t-1}_{-\infty})$. Moreover, for each $\theta \in \Theta$, $(x \mapsto p^{\nu}_{\theta}(X_{t}|X^{t-1}_{-n}))_{n}$ is uniformly integrable because, for any
%\begin{align*}
%	\sup_{n} p^{\nu}_{\theta}(x|X^{t-1}_{-n}) \leq & \sup_{n} \sum_{s,s_{-}} p_{\theta}(X_{t-1},s,x) Q_{\theta}(X_{t-1},s_{-},s) Pr_{\theta}(S_{t-1}=s_{-} | X^{t-1}_{-n})\\
%	\leq & \sum_{s} p_{\theta}(X_{t-1},s,x)
%\end{align*}
%and the RHS is integrable (over $x$) for all $n$. Therefore, by the DCT, $p^{\nu}(\cdot|X^{t-1}_{-\infty},\theta)$ is in fact a density a.s.-$\bar{P}^{\nu}_{\ast}$.
Following the same logic, an analogous result can be shown for $p^{\nu}_{\ast}(\cdot|X^{t-1}_{-\infty})$.
\bigskip
(2) Similarly, for any $n$ such that $-n \leq t-1$,
\begin{align*}
	\sup_{\theta \in \Theta} | p^{\nu}(X_{t}|X^{t-1}_{-\infty},\theta)  -  p_{\theta}^{\nu}(X_{t}|X^{t-1}_{-n}) | \leq  \sum_{s \in \mathbb{S}} p_{\theta} (X_{t-1},s,X_{t}) \prod_{i=-n}^{t-1} (1-\underline{q}(X_{i})) .
\end{align*}
Hence, for any $\theta_{0}$ in $\Theta$,
\begin{align*}
	 & \sup_{\theta \in B(\theta_{0},\delta)} | p^{\nu}(X_{t}|X^{t-1}_{-\infty},\theta)  - p^{\nu}(X_{t}|X^{t-1}_{-\infty},\theta_{0}) |\\
	  \leq & \sup_{\theta \in B(\theta_{0},\delta)} \sum_{s \in \mathbb{S}} p_{\theta} (X_{t-1},s,X_{t}) \frac{1}{1+M} \sum_{n=0}^{M} \prod_{i=-n}^{t-1} (1-\underline{q}(X_{i}))\\
	 & + \frac{1}{1+M} \sum_{n=0}^{M} \sup_{\theta \in B(\theta_{0},\delta)} | p^{\nu}_{\theta}(X_{t}|X^{t-1}_{-n},\theta)  - p^{\nu}_{\theta_{0}}(X_{t}|X^{t-1}_{-n}) |.
\end{align*}
For any $\gamma>0$, choose $M$ such that $\bar{P}^{\nu}_{\ast} \left( \frac{1}{1+M} \sum_{n=0}^{M} \prod_{i=-n}^{t-1} (1-\underline{q}(X_{i})) \geq  \gamma \right) < \gamma $; such $M$ exists by Assumption \ref{ass:qbar-sum}. Given such $M$, for any $\epsilon>0$ and any $\delta>0$, it follows that
\begin{align*}
&\bar{P}^{\nu}_{\ast} \left( \sup_{\theta_{0} \in \Theta} \sup_{\theta \in B(\theta_{0},\delta)} | p^{\nu}(X_{t}|X^{t-1}_{-\infty},\theta)  - p^{\nu}(X_{t}|X^{t-1}_{-\infty},\theta_{0}) | > \epsilon \right) \\
\leq & \bar{P}^{\nu}_{\ast} \left( \sup_{\theta \in \Theta} \sum_{s \in \mathbb{S}} p_{\theta} (X_{t-1},s,X_{t}) \geq 0.5 \epsilon/\gamma  \right) + \gamma/3 \\
& +  \bar{P}^{\nu}_{\ast} \left( \frac{1}{1+M} \sum_{n=0}^{M} \sup_{\theta_{0} \in \Theta} \sup_{\theta \in B(\theta_{0},\delta)} | p^{\nu}_{\theta}(X_{t}|X^{t-1}_{-n})  - p^{\nu}_{\theta_{0}}(X_{t}|X^{t-1}_{-n}) | \geq 0.5 \epsilon \right).
\end{align*}
Let $\delta>0$ be such that for each $n \in \{ 0,...,M  \}$, $\sup_{\theta_{0} \in \Theta} \sup_{\theta \in B(\theta_{0},\delta)} | p^{\nu}_{\theta}(X_{t}|X^{t-1}_{-n})  - p^{\nu}_{\theta_{0}}(X_{t}|X^{t-1}_{-n}) | < 0.5 \epsilon$ a.s.-$\bar{P}^{\nu}_{\ast}$; such $\delta>0$ exists by our conditions. So the third term in the RHS is 0. Also, under our conditions, it follows that $\sup_{\theta \in \Theta } p_{\theta} (X_{t-1},s,X_{t}) = O_{\bar{P}^{\nu}_{\ast}}(1)$ and thus the first term in the RHS can be made smaller than $\epsilon/3$ by a suitably chosen $\gamma \in (0,\epsilon)$. Thus, the desired result holds.
\bigskip
(3) By the definition of $p^{\nu}(.|.,\theta)$, it follows that
\begin{align*}
	p^{\nu}(X_{1}|X^{0}_{-\infty},\theta) \leq \max_{s \in \mathbb{S}} p_{\theta}(X_{0},s,X_{1}),
\end{align*}
and, for some $n \in \mathbb{N}_{0}$,
\begin{align*}
	p^{\nu}(X_{1}|X^{0}_{-\infty},\theta) \geq & 0.5 \sum_{s \in \mathbb{S}} p_{\theta}(X_{0},s,X_{1}) Pr_{\theta}(S_{1} = s \mid X^{0}_{-n} ) \\
	\geq & 0.5 \min_{s \in \mathbb{S}} p_{\theta}(X_{0},s,X_{1}).
\end{align*}
Thus, by our conditions, for any $\delta>0$ and $\theta_{0} \in \Theta$,
\begin{align*}
	\sup_{\theta \in B(\delta,\theta_{0})} \frac{p^{\nu}(X_{1}|X^{0}_{-\infty},\theta)}{p^{\nu}(X_{1}|X^{0}_{-\infty},\theta_{0})} \leq \frac{\overline{p}(X_{0},X_{1})}{\underline{p}(X_{0},X_{1})},
	\end{align*}
	and the RHS is in $L^{1}(\bar{P}^{\nu}_{\ast})$. Thus, by part (2) and the dominated convergence theorem, for any $\epsilon>0$, there exists a $\delta>0$ such that
	\begin{align*}
		E_{\bar{P}^{\nu}_{\ast}} \left[  	\sup_{\theta \in B(\theta_{0},\delta)} \frac{p^{\nu}(X_{1}|X^{0}_{-\infty},\theta)}{p^{\nu}(X_{1}|X^{0}_{-\infty},\theta_{0})}   \right] \leq 1 + \epsilon,
	\end{align*}
	for any $\theta_{0} \in \Theta$. This readily implies Assumption \ref{ass:bdd}(i). Part (ii) also follows with $(x_{0},x_{1}) \mapsto C(x_{0},x_{1}) = \frac{\overline{p}(x_{0},x_{1})}{\underline{p}(x_{0},x_{1})}$.
Similarly, by noting that
\begin{align*}
1 - \frac{p^{\nu}(X_{1}|X^{0}_{-\infty},\theta)}{p^{\nu}_{\ast}(X_{1}|X^{0}_{-\infty})}	\leq  \log \frac{p^{\nu}(X_{1}|X^{0}_{-\infty},\theta)}{p^{\nu}_{\ast}(X_{1}|X^{0}_{-\infty})} \leq \frac{p^{\nu}(X_{1}|X^{0}_{-\infty},\theta)}{p^{\nu}_{\ast}(X_{1}|X^{0}_{-\infty})} - 1,
\end{align*}
it follows that, for any $\theta \in \Theta$,
\begin{align*}
	\left|  \log \frac{p^{\nu}_{\ast}(X_{1}|X^{0}_{-\infty})}{p^{\nu}(X_{1}|X^{0}_{-\infty},\theta)}  \right| \leq 1 + \frac{\overline{p}(X_{0},X_{1})}{\underline{p}(X_{0},X_{1})}.
\end{align*}
Since the RHS is in $L^{1}(\bar{P}^{\nu}_{\ast})$, the results in part (2) and the dominated convergence theorem imply that $\theta \mapsto H(\theta)$ is continuous.
\end{proof}

\section{Sufficient Conditions for Assumptions \ref{ass:qbar-sum} and
\ref{ass:q-sum-p}}

\label{SM:SuffAssumptions}

By exploiting the fact that $(X_{t})_{t=-\infty}^{\infty}$ is $\beta$-mixing
and stationary, the following lemma provides sufficient conditions for
Assumption \ref{ass:qbar-sum}.

\begin{lemma}
\label{lem:suff-A5andA8} Suppose Assumptions \ref{ass:BDD_Q} and
\ref{ass:fx-ergo} hold.

(1) Suppose further there exists $l \geq1$ such that
\begin{align*}
E_{\bar{P}^{\ast}_{\nu}} \left[  (1-\underline{q}(X_{1}))^{l^{\prime}}
\right]  < 1 ~and~ E_{\bar{P}^{\ast}_{\nu}} \left[  C(X_{1},X_{0})^{l}
\right]  <\infty,
\end{align*}
where $1/l^{\prime}+ 1/l = 1$. Then,%
\begin{align*}
\lim_{T\rightarrow\infty} E_{\bar{P}^{\ast}_{\nu}} \left[  T^{-1} \sum
_{t=1}^{T} C(X_{t-1},X_{t}) \prod_{i=0}^{t-1} (1-\underline{q}(X_{i}))
\right]  = 0
\end{align*}

(2) Suppose further $E_{\bar{P}^{\ast}_{\nu}} \left[  (1-\underline{q}%
(X_{1}))^{\frac{2a}{1-a}} \right]  < 1$. Then, Assumption \ref{ass:q-sum-p} holds.
\end{lemma}

Clearly, by the Markov inequality, Part (1) implies that
\begin{align*}
\lim_{T \rightarrow\infty} \bar{P}^{\ast}_{\nu} \left(  T^{-1} \sum_{t=1}^{T}
C(X_{t-1},X_{t}) \prod_{i=0}^{t-1} (1-\underline{q}(X_{i})) \geq
\epsilon\right)  = 0 .
\end{align*}

\begin{proof}
	(1) We use well-known coupling results for $\beta$-mixing sequences (see \citesupp{yu1994rates}). For any $q \in \mathbb{N}$, let $(\tilde{X}_{t})_{t=-\infty}^{\infty}$ be such that: (a) for any $i \in \mathbb{N}_{0}$, $\bar{U}_{i} \equiv (\tilde{X}_{iq+1},...,\tilde{X}_{iq+q})$ has the same distribution as $U_{i} \equiv (X_{iq+1},...,X_{iq+q})$; (b) the sequence $(\bar{U}_{i})_{i~even}$ is i.i.d. and so is $(\bar{U}_{i})_{i~odd}$; (c) for any $ i \in \mathbb{N}_{0}$, $\mathbb{P}(\bar{U}_{i} \ne U_{i}) \leq \beta_{q}$, where $\mathbb{P}$ is the product measure of the random processes $(X_{t})_{t=-\infty}^{\infty}$ and $(\tilde{X}_{t})_{t=-\infty}^{\infty}$.
	
By stationarity, for any $t \in \mathbb{N}$,
\begin{align*}
	E_{\bar{P}^{\ast}_{\nu}} \left[ C(X_{t-1},X_{t}) \prod_{i=0}^{t-1} (1-\underline{q}(X_{i}))\right] = & E_{\bar{P}^{\ast}_{\nu}} \left[ C(X_{t-1},X_{t}) (1-\underline{q}(X_{t})) ... (1-\underline{q}(X_{0}))\right]\\
	= & E_{\bar{P}^{\ast}_{\nu}} \left[ C(X_{-1},X_{0}) (1-\underline{q}(X_{-1})) ... (1-\underline{q}(X_{-t}))\right]\\
	= & E_{\bar{P}^{\ast}_{\nu}} \left[ C(X_{-1},X_{0}) \prod_{i=1}^{t}(1-\underline{q}(X_{-i})) \right].
\end{align*}
By the conditions of the lemma, $E[C(X_{-1},X_{0})^{l}] < \infty$, so by Jensen's inequality (and the fact that $1/l' \leq 1$, where $1/l +1/l' = 1$), it suffices to show that
\begin{align*}
\lim_{T\rightarrow \infty} T^{-1} \sum_{t=1}^{T} E_{\bar{P}^{\ast}_{\nu}} \left[  \prod_{i=0}^{t-1}(1-\underline{q}(X_{i}))^{l'} \right] = 0.
\end{align*}
Note that, for any $1 \leq m \leq T$,
\begin{align*}
	T^{-1} \sum_{t=1}^{T} E_{\bar{P}^{\ast}_{\nu}} \left[  \prod_{i=0}^{t-1}(1-\underline{q}(X_{i}))^{l'} \right] \leq \frac{m}{T} + T^{-1} \sum_{t=m+1}^{T} E_{\bar{P}^{\ast}_{\nu}} \left[  \prod_{i=0}^{t-1}(1-\underline{q}(X_{i}))^{l'} \right].
\end{align*}
We now examine the term $E_{\bar{P}^{\ast}_{\nu}} \left[  \prod_{i=0}^{t-1}(1-\underline{q}(X_{i}))^{l'} \right]$ for any $t \equiv q k$ and some positive integers $q$ and $k$; if this does not hold, $t$ in the product is simply replaced with $qk(t)$, where $k(t)$ is the largest integer such that $t \geq qk$.
Letting $M_{j} \equiv \{jq +1,...,(j+1)q  \}$, it follows that $\{1,...,t\} = \cup_{j=0}^{k-1} M_{j} $. Thus,
	\begin{align*}
	E_{\bar{P}^{\ast}_{\nu}} \left[  \prod_{i=0}^{t-1}(1-\underline{q}(X_{i}))^{l'} \right]  = & E_{\bar{P}^{\ast}_{\nu}} \left[  \prod_{j=0}^{k-1} \prod_{i \in M_{j} }(1-\underline{q}(X_{i}))^{l'} \right] \\
	\leq & E_{\bar{P}^{\ast}_{\nu}} \left[  \prod_{j \in E_{k}} \prod_{i \in M_{j} }(1-\underline{q}(X_{i}))^{l'} \right],
	\end{align*}
	where $E_{k} \equiv \{  0 \leq j \leq k-1 \colon j~is~even  \}$. Observe that
		\begin{align*}
	E_{\bar{P}^{\ast}_{\nu}} \left[  \prod_{j \in E_{k}} \prod_{i \in M_{j} }(1-\underline{q}(X_{i}))^{l'} \right] = & E_{\mathbb{P}} \left[  |\prod_{j \in E_{k}} \prod_{i \in M_{j} }(1-\underline{q}(X_{i}))^{l'} - \prod_{j \in E_{k}} \prod_{i \in M_{j} }(1-\underline{q}(\bar{X}_{i}))^{l'}| \right] \\
	& + E \left[ \prod_{j \in E_{k}} \prod_{i \in M_{j} }(1-\underline{q}(\bar{X}_{i}))^{l'} \right]\\
	\leq & \mathbb{P} \left(  \exists j \in E_{k} \colon \bar{U}_{j}  \ne U_{j}  \right) + E \left[ \prod_{j \in E_{k}} \prod_{i \in M_{j} }(1-\underline{q}(\bar{X}_{i}))^{l'} \right] \\
	\leq & \beta_{q} k + E \left[ \prod_{j \in E_{k}} \prod_{i \in M_{j} }(1-\underline{q}(\bar{X}_{i}))^{l'} \right].
	\end{align*}
	Because of the properties of the process $(\bar{U}_{j})_{j~even}$,
	\begin{align*}
		E \left[ \prod_{j \in E_{k}} \prod_{i \in M_{j} }(1-\underline{q}(\bar{X}_{i}))^{l'} \right] \leq & \prod_{j \in E_{k}} E \left[  (1-\underline{q}(\bar{X}_{jq+1}))^{l'} \right] \\
		= &  \prod_{j \in E_{k}} E_{\bar{P}^{\ast}_{\nu}} \left[  (1-\underline{q}(X_{jq+1}))^{l'} \right] \\
		= &  \left( E_{\bar{P}^{\ast}_{\nu}}  \left[  (1-\underline{q}(X_{1}))^{l'} \right] \right)^{(k-1)/2}.
	\end{align*}
	
	Noting that $t=qk$, these results imply that
	\begin{align*}
	  E_{\bar{P}^{\ast}_{\nu}} \left[  \prod_{i=0}^{t-1}(1-\underline{q}(X_{i}))^{l'} \right] \leq  \beta_{q} t/q +  \left( E_{\bar{P}^{\ast}_{\nu}}  \left[  (1-\underline{q}(X_{1}))^{l'} \right] \right)^{(t/q-1)/2}.
	\end{align*}
	
	Therefore,	for any $1 \leq m \leq T$,
	\begin{align*}
	T^{-1} \sum_{t=1}^{T} E_{\bar{P}^{\ast}_{\nu}} \left[  \prod_{i=0}^{t-1}(1-\underline{q}(X_{i}))^{l'} \right] \leq \frac{m}{T} + T^{-1} \sum_{t=m+1}^{T} \{ \beta_{q} t/q +  \left( E_{\bar{P}^{\ast}_{\nu}}  \left[  (1-\underline{q}(X_{1}))^{l'} \right] \right)^{(t/q-1)/2}\}.
	\end{align*}
	
	By Lemma \ref{lem:sta-ergo}, $\beta_{q} = \exp \{ q \log \gamma    \}$. This fact and the condition $E_{\bar{P}^{\ast}_{\nu}}  \left[  (1-\underline{q}(X_{1}))^{l'} \right] < 1$ imply that we can take, for instance, $q \equiv t^{1/2}$ and $m = \sqrt{T}$, so that the RHS vanishes as $T$ diverges.
	
	\bigskip
	
	(2) We wish to show that 		
	\begin{align*}
		\sum_{t=1}^{\infty} \left( E_{\bar{P}^{\nu}_{\ast}} \left[ \prod_{i=0}^{t-1} (1 - \underline{q}(X_{i}))^{\frac{2a}{1-a}} \right]	\right)^{p\frac{1-a}{2a}} < \infty.
	\end{align*}
	
	By our previous calculations,
	\begin{align*}
		E_{\bar{P}^{\nu}_{\ast}} \left[ \prod_{i=0}^{t-1} (1 - \underline{q}(X_{i}))^{\frac{2a}{1-a}} \right] \leq \beta_{q} t/q +  \left( E_{\bar{P}^{\ast}_{\nu}}  \left[  (1-\underline{q}(X_{1}))^{\frac{2a}{1-a}} \right] \right)^{(t/q-1)/2};
	\end{align*}
	thus, it suffices to show that there exists a choice $t \mapsto q(t)$ such that $\sum_{t=1}^{\infty} \{\beta_{q(t)} t/q(t) +  \left( E_{\bar{P}^{\ast}_{\nu}}  \left[  (1-\underline{q}(X_{1}))^{\frac{2a}{1-a}} \right] \right)^{(t/q(t)-1)/2}\} < \infty$. Let $q(t) = \sqrt{t}$; then as $E_{\bar{P}^{\ast}_{\nu}}  \left[  (1-\underline{q}(X_{1}))^{\frac{2a}{1-a}} \right] < 1$ by assumption, it follows that the second term in the sum is finite. Since $t \mapsto \beta_{q(t)} t/q(t) = \exp\{ \sqrt{t} \log \gamma  \} \sqrt{t} $ is summable (note that $\log \gamma < 0$), the same holds true for the first term.
	
\end{proof}

\section{Proofs and Results for Example \ref{exa:Canon}}

\label{SM:examples}

In what follows $e_{max}(M)$ and $e_{min}(M)$ denote the maximal and minimal
eigenvalues of a matrix $M$.

\begin{lemma}
\label{lem:fx-ergo-holds} Assumption \ref{ass:fx-ergo} holds.
\end{lemma}

\begin{proof}
	For each $s \in \mathbb{S}$, we apply Theorem 3.3 in \citesupp{douc2004practical}. To do so, we first verify their Assumptions 3.3 and 3.4. In our case, $\epsilon \sim N(0,\Sigma(s))$, so their Assumption 3.3 is satisfied for any $z_{0}$ and $\gamma_{0} = 1$. In their notation, $g(x) \equiv \mu(s) + \Phi^{\intercal}x$. Observe that $||g(x)||\leq ||\mu(s)|| + e_{max}(\Phi \Phi ^{\intercal})||x||$. By assumption, $e_{max}(\Phi \Phi ^{\intercal}) \equiv \gamma < 1$ and
	\begin{align*}
	||\mu(s)||  \leq (1 -e_{max}(\Phi \Phi ^{\intercal}) ) ||x||(1 - ||x||^{-0.5})
	\end{align*}
	for all $x$ such that $||x|| \geq R_{0}$. Such an $R_{0}$ exists because $||\mu(s)||$ is bounded and $e_{max}(\Phi \Phi ^{\intercal}) < 1$. This choice ensures that $||g(x)|| \leq ||x||(1 - ||x||^{0.5})$, which, in turn, ensures the validity of their Assumption 3.4 with $r=1$ and $\rho = 0.5$. By their Theorem 3.3, Assumption \ref{ass:fx-ergo}(i) holds.
	
	Assumption \ref{ass:fx-ergo}(ii),(iii) is satisfied because $\inf_{x \in A} P_{\ast}(s,x,C) \geq \int_{C} \inf_{x \in A} \underline{p}(x,a) da $ and, since $A$ is bounded, it follows that $\inf_{x \in A} \underline{p}(x,a) \geq \exp\{ \tilde{D} + (x^{\prime})^{T} \tilde{F} a + \tilde{G} a   \}$, so the RHS plays the role of the measure $\varpi$, which clearly is such that $\varpi(A) > 0$.
\end{proof}

Let $\overline{\kappa}\equiv\max_{s\in\mathbb{S}}\sup_{\Sigma(s)\in\Theta
}e_{max}(\Sigma(s))$ and $\underline{\kappa}\equiv\min_{s\in\mathbb{S}}%
\inf_{\Sigma(s)\in\Theta}e_{min}(\Sigma(s))$; $\overline{m}=\max
_{s\in\mathbb{S}}\sup_{\mu(s)\in\Theta}||\mu(s)||$; $\overline{M}=\sup
_{\Phi\in\Theta}e_{max}(\Phi\Phi^{\intercal})$ and $\underline{M}=\min
_{\Phi\in\Theta}e_{min}(\Phi\Phi^{\intercal})$; $\overline{\kappa}_{\ast
}\equiv\max_{s\in\mathbb{S}}e_{max}(\Sigma_{\ast}(s))$ and $\underline{\kappa
}_{\ast}\equiv\min_{s\in\mathbb{S}}e_{min}(\Sigma_{\ast}(s))$; $\overline
{m}_{\ast}=\max_{s\in\mathbb{S}}||\mu_{\ast}(s)||$; $\overline{M}_{\ast
}=e_{max}(\Phi_{\ast}\Phi_{\ast}^{\intercal})$ and $\underline{M}_{\ast
}=e_{min}(\Phi_{\ast}\Phi_{\ast}^{\intercal})$. By the assumptions in the
text, $\overline{\kappa},\overline{\kappa}_{\ast},\underline{\kappa
},\underline{\kappa}_{\ast}$, $\underline{M},\overline{M},\underline{M}_{\ast
},\overline{M}_{\ast}$ are all in $(0,\infty)$.

\begin{lemma}
\label{lem:exa-CBound} There exists a $C\in\lbrack1,\infty)$ such that, for
any $(x,y)$ and any $s$,
\begin{align*}
&  f_{\mathcal{N}}((y-\Phi^{\intercal}x-\mu(s))\Sigma^{-1/2}(s))\\
\leq &  C\exp\{-0.5\underline{\kappa}(||y||^{2}+\underline{M}||x||^{2}%
-2\sqrt{\overline{M}}||x||||y||)+\underline{\kappa}(||y||+\sqrt{\overline{M}%
}||x||)\overline{m}\}
\end{align*}
and
\begin{align*}
&  f_{\mathcal{N}}((y-\Phi^{\intercal}x-\mu(s))\Sigma^{-1/2}(s))\\
\geq &  C^{-1}\exp\{-0.5\overline{\kappa}(||y||^{2}+\overline{M}%
||x||^{2}+2\sqrt{\overline{M}}||x||||y||)-\overline{\kappa}(||y||+\sqrt
{\overline{M}}||x||)\overline{m}\}.
\end{align*}

An analogous bound holds for $f_{\mathcal{N}}((y-\Phi_{\ast}^{\intercal}%
x-\mu_{\ast}(s))\Sigma_{\ast}^{-1/2}(s))$.
\end{lemma}

\begin{proof}
	We have $(a,s,b) \mapsto p_{\theta}(a,s,b) \asymp \phi( (b - \mu(s) - \Phi^{\intercal} a ) \Sigma^{-1/2}(s) ) $. Thus, there exists a $C \in [1,\infty)$ such that, for any $u \in \mathbb{R}^{p}$,
	\begin{align*}
	& f_{\mathcal{N}}( (u - \mu(s)  ) \Sigma^{-1/2}(s) ) \leq C \exp\{ -0.5 \underline{\kappa} ||u - \mu(s)||^{2}    \}, \\
	&  f_{\mathcal{N}}( (u - \mu(s)  ) \Sigma^{-1/2}(s) ) \geq C^{-1}  \exp\{ -0.5 \overline{\kappa} ||u - \mu(s)||^{2}    \},
	\end{align*}
	where $ \overline{\kappa} \in (0, \infty)$ and  $ \underline{\kappa} \in (0, \infty)$ are the largest and smallest eigvenvalues within the collection $(\Sigma(s) )_{s \in \mathbb{S}}$ in $\Theta$.
	
	Observe that
	\begin{align*}
	||x-y||^{2} \leq ||x||^{2} + ||y||^{2} +  2 ||x||  ||y||,
	\end{align*}
	and, if $||x|| > ||y||$,
	\begin{align*}
	||x-y||^{2} \geq (||x|| - ||y|| )^{2} = ||x||^{2} + ||y||^{2} - 2 ||x|| ||y||,
	\end{align*}
	and the same for $||x|| < ||y||$.
	
	Hence,
	\begin{align*}
	f_{\mathcal{N}}( (u - \mu(s)  ) \Sigma^{-1/2}(s) ) \leq & C \exp\{ -0.5 \underline{\kappa} (||u||^{2} + ||\mu(s)||^{2} - 2||u|| || \mu(s)|| )    \}\\
	\leq & C \exp\{ -0.5 \underline{\kappa}  ||u||^{2}  + \underline{\kappa}  ||u|| \overline{m}    \}
	\end{align*}	
	and
	\begin{align*}
	f_{\mathcal{N}}( (u - \mu(s)  ) \Sigma^{-1/2}(s) ) \geq & C^{-1}  \exp\{ -0.5 \overline{\kappa}  (||u||^{2} + ||\mu(s)||^{2} + 2||u|| || \mu(s)|| )     \} \\
	\geq & C^{-1}  \exp\{ -0.5 \overline{\kappa}  ||u||^{2} - 0.5 \overline{\kappa} m^{2} - \overline{\kappa}||u|| m      \} \\
	\equiv & C^{-1} \exp\{ -0.5 \overline{\kappa}  ||u||^{2}
	- \overline{\kappa}||u|| m      \},
	\end{align*}
	where $\overline{m} = \sup_{\mu(.) \in \Theta} \max_{s \in \mathbb{S}} ||\mu(s)||$.
	
	Observe that
	\begin{align*}
	||y-\Phi^{\intercal}x||^{2} \leq & ||y||^{2} + ||\Phi^{\intercal} x ||^{2} +  2 ||\Phi^{\intercal} x ||  ||y|| \\
	\leq & ||y||^{2} + \overline{M} ||x ||^{2} +  2 \sqrt{\overline{M}} ||x||  ||y||
	\end{align*}
	and
	\begin{align*}
	||y-\Phi^{\intercal}x||^{2} \geq & ||y||^{2} + ||\Phi^{\intercal} x ||^{2} - 2 ||\Phi^{\intercal} x ||  ||y|| \\
	\geq & ||y||^{2} + \underline{M} ||x ||^{2} -  2 \sqrt{\overline{M}} ||x||  ||y||,
	\end{align*}
	where $\overline{M} = \sup_{\Phi \in \Theta} e_{max}(\Phi \Phi^{\intercal} )$ and $\underline{M} = \min_{\Phi \in \Theta} e_{min}(\Phi \Phi^{\intercal} )$, which, by our conditions, are finite and positive. Hence,
	\begin{align*}
	& f_{\mathcal{N}}( (y - \Phi^{\intercal} x - \mu(s)  ) \Sigma^{-1/2}(s) )\\
	\leq  & C \exp\{ -0.5 \underline{\kappa} (||y||^{2} + \underline{M} ||x ||^{2} -  2 \sqrt{\overline{M}} ||x||  ||y||  )  + \underline{\kappa} (||y|| + \sqrt{\overline{M}}||x||  )  \overline{m}    \}
	\end{align*}
	and
	\begin{align*}
	& f_{\mathcal{N}}( (y - \Phi^{\intercal} x - \mu(s)  ) \Sigma^{-1/2}(s) ) \\
	\geq  &  C^{-1} \exp\{ -0.5 \overline{\kappa} (||y||^{2} + \overline{M} ||x ||^{2} +  2 \sqrt{\overline{M}} ||x||  ||y||  )  - \overline{\kappa} (||y|| + \sqrt{\overline{M}}||x||  )  \overline{m}    \}.
	\end{align*}
\end{proof}

\begin{lemma}
\label{lem:exa-CBound2} Suppose there exists a $l \geq1$ such that
$l\overline{\kappa} < l \underline{\kappa} + \underline{\kappa}_{\ast} $ and
$l \overline{\kappa} \overline{M} < l \underline{\kappa} \underline{M} +
\underline{\kappa}_{\ast}$. Then,
\[
E_{\bar{P}^{\ast}_{\nu}}[\exp\{ -0.5l ((b_{1} - a_{1}) ||Y||^{2} + (b_{2} -
a_{2}) ||X ||^{2} - 2 (b_{3} + a_{3}) ||X|| ||Y|| ) + l(b_{4} +a_{4}) ||Y|| +
l(b_{5} + a_{5}) ||X|| \}] < \infty,
\]
where
\begin{align*}
&  a_{1} = \overline{\kappa},~b_{1} = \underline{\kappa},\\
&  a_{2} = \overline{\kappa} \overline{M},~b_{2} = \underline{\kappa}
\underline{M},\\
&  a_{3} = \overline{\kappa} \sqrt{ \overline{M} },~b_{3} = \underline{\kappa}
\sqrt{ \overline{M} },\\
&  a_{4} = \underline{\kappa},~b_{4} = \overline{\kappa},\\
&  a_{5} = \underline{\kappa} \sqrt{ \overline{M} },~ b_{5} = \overline
{\kappa} \sqrt{ \overline{M} }.
\end{align*}

\end{lemma}

\begin{remark}
Before going to the proof, we discuss the conditions in the lemma. They
basically require that the \textquotedblleft spread\textquotedblright\ of the
eigenvalues of the matrices $\Sigma(\cdot)$ and $\Phi\Phi^{\ast}$ is not too
large relative to the eigenvalues in $\Sigma_{\ast}(\cdot)$. This condition
comes naturally since we are essentially requiring that the ratio of two
exponential functions is integrable with respect to a Gaussian measure. For
instance, if $\Sigma(\cdot)$, $\Sigma_{\ast}(\cdot)$ and $\Phi\Phi^{\intercal
}$, $\Phi_{\ast}\Phi_{\ast}^{\intercal}$ are matrices with eigenvalues bounded
between $0<a$ and $a+\Delta$, then sufficient conditions are given by
$l\Delta<a$ and $l(2a\Delta+(\Delta)^{2})<a$, which is equivalent to
$\frac{l\Delta^{2}}{1-2l\Delta}<a$. $\triangle$
\end{remark}

\begin{proof}
	It is enough to show that
	\begin{align*}
	E_{\bar{P}^{\ast}_{\nu}}[T_{1}(Y)] \equiv &E_{\bar{P}^{\ast}_{\nu}}[\exp\{  -0.5l (b_{1} - a_{1}) ||Y||^{2}    + l(b_{4} +a_{4}) ||Y||   \}] < \infty, \\
	E_{\bar{P}^{\ast}_{\nu}}[T_{2}(X)] \equiv &E_{\bar{P}^{\ast}_{\nu}}[\exp\{  -0.5l (b_{2} - a_{2}) ||X ||^{2}  + l(b_{5} + a_{5}) ||X||    \}] < \infty,\\
	E_{\bar{P}^{\ast}_{\nu}}[T_{3}(X,Y)] \equiv & E_{\bar{P}^{\ast}_{\nu}}[\exp\{   l (b_{3} + a_{3}) ||X||  ||Y||  )  \}] < \infty.
	\end{align*}

	For any $d \in \{1,2\}$, suppose there exists $\varphi$ such that $\int T_{d}(b) p_{\ast}(x,s,b) db \leq \varphi(x) T_{d}(x) $ for any $x$, and, for any $\gamma>0$, $\{ x : \varphi(x) \geq \gamma   \}$ is either empty or compact.  Then, for any  $\gamma > 0$,
	\begin{align*}
	\int T_{d}(x) \nu(dx) = & \int 1\{ \varphi(x) \leq \gamma\}  T_{d}(x) \nu(dx)  + \int 1\{ \varphi(x) > \gamma\}  T_{d}(x) \nu(dx) \\
	= & \int \int 1\{ \varphi(x) \leq \gamma\}  T_{d}(b) p_{\ast}(x,s,b) db \nu(dx,ds)  + \int 1\{ \varphi(x) > \gamma\}  T_{d}(x) \nu(dx)\\
	\leq & \gamma \int T_{d}(x) \nu(dx)  + \sup_{x : \varphi(x) \geq \gamma} T_{d}(x),
	\end{align*}
	where second line follows because $\nu$ is the invariant probability distribution. Since $\{x : \varphi(x) \geq \gamma \}$ is bounded and compact (if it is non-empty), $\sup_{x : \varphi(x) \geq \gamma} T_{d}(x) \leq M < \infty$. Choosing $\gamma < 1$, it follows that $\int T_{d}(x) \nu(dx) \leq \frac{M}{1-\gamma} < \infty$,  as desired.
	
	We now show that $\int T_{1}(b) p_{\ast}(x,s,b) db \leq \varphi(x) T_{1}(x) $ for any $x$. To do this, note that
	\begin{align*}
	\int T_{1}(b) p_{\ast}(x,s,b) db \leq &  \int  \exp \{  0.5 l (\overline{\kappa} - \underline{\kappa})||y||^{2}   +  l (\overline{\kappa} + \underline{\kappa} ) ||y||     \}  p_{\ast}(x,s,y) dy,
	\end{align*}
	and, by Lemma \ref{lem:exa-CBound},
	\begin{align*}
	f_{\mathcal{N}}( (y - \Phi^{\intercal}_{\ast} x - \mu_{\ast}(s)  ) \Sigma^{-1/2}_{\ast}(s) ) \leq &  \exp\{ -0.5 \underline{\kappa}_{\ast} (||y||^{2} -  2 \sqrt{\overline{M}_{\ast}} ||x||  ||y||  )  + \underline{\kappa}_{\ast} ||y||  \overline{m}_{\ast}    \} \\
	& \times B_{\ast}(x),
	\end{align*}
	with $B_{\ast}(x) \equiv  C_{\ast} \exp\{ -0.5 \underline{\kappa}_{\ast}  \underline{M}_{\ast} ||x||^{2}   + \underline{\kappa}_{\ast}  \sqrt{\overline{M}_{\ast} }||x||  \overline{m}_{\ast}     \}  $. Therefore,
	
	\begin{align*}
	& \frac{\int T_{1}(b) p_{\ast}(x,s,y) dy}{B_{\ast}(x)}\\
	\leq &  \int \exp \{  0.5 (\overline{\kappa} - \underline{\kappa})||y||^{2} +  l (\overline{\kappa} + \underline{\kappa} ) ||y||   \} \exp\{ -0.5 \underline{\kappa}_{\ast}  (||y||^{2} -  2 \sqrt{\overline{M}_{\ast} } ||x||  ||y||  )  + \underline{\kappa}_{\ast}   \overline{m}_{\ast}  ||y||    \} dy \\
	= &  \int \exp \{  0.5 (l(\overline{\kappa} - \underline{\kappa}) - \underline{\kappa}_{\ast})||y||^{2}  +  \underline{\kappa}_{\ast}  \sqrt{\overline{M}_{\ast} } ||x||  ||y||  + (\underline{\kappa}_{\ast}   \overline{m}_{\ast}  + l (\overline{\kappa} + \underline{\kappa} ) ) ||y||    \} dy.
	\end{align*}
	By our conditions, $l(\overline{\kappa} - \underline{\kappa}) - \underline{\kappa}_{\ast} <0$. Hence, the expression above is an integral of an exponential function of a quadratic form with negative leading coefficient, and is thus finite. Moreover, after some algebra, there exists a finite constant $C$, such that  $ \int T_{1}(b) p_{\ast}(x,s,y) dy \leq C B_{\ast}(x) \exp  \{   \underline{\kappa}_{\ast}  \sqrt{\overline{M}_{\ast} } ||x|| /(-(l(\overline{\kappa} - \underline{\kappa}) - \underline{\kappa}_{\ast}))  \})$; we re-define the RHS as $C B_{\ast}(x) \exp  \{  D \underline{\kappa}_{\ast}  \sqrt{\overline{M}_{\ast} } ||x||   \})$ with $D>0$. Therefore, the result holds with
	\begin{align*}
	x \mapsto \varphi(x) \equiv    \exp  \{ - 0.5 (\underline{\kappa}_{\ast}  \underline{M}_{\ast}  +  l ( \overline{\kappa} - \underline{\kappa})  ) ||x||^{2}   +   ( \underline{\kappa}_{\ast}  \sqrt{\overline{M}_{\ast} } (D+\overline{m}_{\ast} )  + l (b_{4} + a_{4}) ) ||x||     \}.
	\end{align*}
	As the coefficient $||x||^{2}$ is $- 0.5 (\underline{\kappa}_{\ast}  \underline{M}_{\ast} +  l ( \overline{\kappa} - \underline{\kappa})  )$, which is negative, the function satisfies the required conditions.
	
	The case for $d=2$ is analogous and is thus omitted; for this case, we use the restriction that $l \overline
	{\kappa} \overline{M} < l \underline{\kappa} \underline{M} + \underline{\kappa
	}_{\ast}$, instead of $l(\overline{\kappa} - \underline{\kappa}) - \underline{\kappa}_{\ast} <0$.
	
	Finally, observe that, if $\int T_{3}(x,y) p_{\ast}(x,s,y) \nu (dx, ds) dy \leq C \int \exp\{ c_{1} ||x||   \} \nu(dx)$ for some $c_{1} < \infty$, then we can follow the same approach as before to show that $\int T_{3}(x,y) p_{\ast}(x,s,y) \nu (dx, ds) dy < \infty$.  Observe that, by the same calculations as before,
	\begin{align*}
	& \int T_{3}(x,y) p_{\ast}(x,s,y) dy \\
	\leq & B_{\ast}(x) \int    \exp\{ l (a_{3} + b_{3} ) ||x|| ||y||  -0.5 \underline{\kappa}  (||y||^{2} -  2 \sqrt{\overline{M} } ||x||  ||y||  )  + \underline{\kappa}  ||y||  \overline{m}     \}  dy \\
	\leq & C B_{\ast}(x) \exp\{ (l (a_{3} + b_{3} ) + \underline{\kappa}  \sqrt{\overline{M} } )  ||x|| \}.
	\end{align*}
	The RHS is of the form $C \exp\{ c_{1} ||x|| \}$, as desired.
\end{proof}

\section{Pseudo-True Parameter Set}

\label{SM:examplemixture}

In this Section, we present an example in which it is possible to characterize
the pseudo-true parameter set $\Theta_{\ast}$ of Theorem~\ref{thm:consistent}.
The example shows that for a subclass of the models in Example \ref{exa:MSR},
namely hidden Markov models with covariate-dependent transition probabilities,
$\Theta_{\ast}$ can be characterized when the (misspecified) model is a simple
mixture. Specifically, we show that, even though the mixture model
misspecifies the dependence structure of the hidden state -- and thus of the
overall system -- parameters related to the outcome equation are consistently estimated.

\begin{example}
\label{exa:HMM-miss} Suppose the true model is a version of the hidden Markov
model of Example \ref{exa:HMM} with%
\begin{align*}
Y_{t}  &  =\mu^{\ast}(S_{t})+\sigma^{\ast}(S_{t})U_{1,t},\\
Z_{t}  &  =\mu_{2}^{\ast}+\psi^{\ast}Z_{t-1}+\sigma_{2}^{\ast}U_{2,t},
\end{align*}
where $((U_{1,t},U_{2,t}))_{t}$ are i.i.d. with zero mean and covariance
matrix $\Sigma=\left[
\begin{array}
[c]{cc}%
1 & \rho^{\ast}\\
\rho^{\ast} & 1
\end{array}
\right]  $, and $S_{t}\sim Q_{\ast}(Z_{t-1},S_{t-1},\cdot)$. The only
additional requirement is that the process $(S_{t},Z_{t})_{t}$ is stationary
with invariant distribution $\nu_{ZS}$. Since the focus of the example is
misspecification of the Markov transition mechanism, we also assume that the
process for $(Z_{t})_{t}$ is correctly specified. The researcher, however,
postulates a mixture model (with i.i.d. regimes) of the form%
\begin{align*}
Y_{t}  &  =\mu(S_{t})+\sigma(S_{t})\varepsilon_{t},\\
S_{t}  &  \sim Q_{\bar{\vartheta}},
\end{align*}
where $Q_{\bar{\vartheta}}$ is parameterized as $Q_{\bar{\vartheta}}%
(s)=\bar{\vartheta}_{s}$ for $s\in\{0,1\}$, with $\bar{\vartheta}_{s}\in(0,1)$
(see, e.g., \citesupp{McLachlan2000}).\footnote{For our results to go through,
it suffices that a weaker condition holds, namely there exists a
$\bar{\vartheta}~$such that$~s\mapsto Q_{\bar{\vartheta}}(s)=E_{\nu_{ZS}%
}[Q_{\ast}(Z,S,s)]$. The idea is that the parametrization is \textquotedblleft
rich enough\textquotedblright\ to mimic the average behavior of the true
transitions.}

%We assume that $\bar{\vartheta}%
%\in\mathbb{R}^{m}$, for some $m\geq1$, and that
%\[
%\exists\bar{\vartheta}~:~s\mapsto Q_{\bar{\vartheta}}(s)=E_{\nu_{ZS}}[Q_{\ast
%}(Z,S,s)].
%\]
%The idea is that the parametrization is \textquotedblleft rich
%enough\textquotedblright\ to mimic the average behavior of the true
%transitions. For instance, the assumption is trivially satisfied if
%$Q_{\bar{\vartheta}}(s)=\bar{\vartheta}_{s}$ for any $s$.

%The only restriction on the function $m$ is that the process $(Z_{t},S_{t})_{t}$ is stationary with invariant probability given by $\nu_{ZS}$.
%XXX Finally, we assume that $f$ and $\beta \mapsto Q_{\beta}$ are such that the
%maximum of
%\begin{align*}
%\theta \mapsto E_{\bar{P}^{\nu}_{\ast}} \left[ \log \sum_{s \in \mathbb{S}}
%Q_{\beta}(s) \sigma(s)^{-1} f\left(\frac{Y_{1} - \mu(s) - \phi Y_{0}}{%
%\sigma(s) } \right) \right]
%\end{align*}
%is in the interior of the parameter space and it is unique. XXX

The following proposition shows that, although the mixture model is
misspecified, ML estimates correctly the parameters of the outcome equation.

\begin{proposition}
\label{pro:Exa-ID} The choice $\mu=\mu_{\ast}$, $\sigma=\sigma_{\ast}$, and
$\bar{\vartheta}$ such that $Q_{\bar{\vartheta}}=E_{\nu_{ZS}}[Q_{\ast
}(Z,S,\cdot)]$ is a pseudo-true parameter, i.e., minimizes $\theta\mapsto
E_{\bar{P}_{\ast}^{\nu}}\left[  \log\sum_{s\in\mathbb{S}}Q_{\bar{\vartheta}%
}(s)\sigma^{-1}(s)f\left(  \frac{Y_{1}-\mu(s)-\phi Y_{0}}{\sigma(s)}\right)
\right]  $, where $f$ is the probability density function of $U_{1,1}$.
\end{proposition}

\begin{proof}[Proof]
	Throughout the proof, we omit the sub-index "1" from the parameters in the outcome equation to ease the notational burden.
	
	Observe that $H^{\ast}$ is proportional to
	\begin{align*}
	\theta \mapsto - \int \log \frac{\sum_{s\in \mathbb{S}} Q_{\bar{\vartheta}}(s) \sigma(s)^{-1} f\left(  (y - \mu(s) )/\sigma(s)  \right)}{f_{\ast}(y)} f_{\ast}\left( y  \right) dy,
	\end{align*}
	where \begin{align*}
	y \mapsto f_{\ast}(y) = \sum_{s \in \mathbb{S}} Pr_{\ast}(S_{1} = s)  \sigma_{\ast}(s)^{-1} f\left(  (y - \mu_{\ast}(s) )/\sigma_{\ast}(s)  \right),
	\end{align*}
	and $Pr_{\ast}$ is the true probability over the hidden state, given by
	\begin{align*}
	s \mapsto Pr_{\ast}(S_{1} = s) = \int  \sum_{s_{0} \in \mathbb{S}}  Q_{\ast}(z,s_{0},s) \nu_{ZS}(dz,ds_{0}).
	\end{align*}
	It is well-known that the minimizers of this function are all $\theta$ such that $$\sum_{s\in \mathbb{S}} Q_{\bar{\vartheta}}(s) \sigma(s)^{-1} f\left(  (\cdot - \mu(s) )/\sigma(s)  \right) = f_{\ast}(\cdot).$$ It is straightforward to check that by choosing $\bar{\vartheta}$ such that $Q_{\bar{\vartheta}}(\cdot) = 	Pr_{\ast}(S_{1} = \cdot)$ and $\mu=\mu_{\ast}$ and $\sigma = \sigma_{\ast}$, the equality holds.
\end{proof}

\end{example}

\begin{remark}
This result may be of practical interest because it indicates that dependence
of the regimes in a hidden Markov model may be ignored if the parameters of
interest are only those of the outcome equation. Note, however, that the
result does not hold in cases where the right-hand side of the outcome
equation contains lagged values of the dependent variable.
%, i.e., for the regime-switching autoregressive model in Example \ref{exa:MSR}.
In this case, the misspecification of the transition function will affect
estimation of all parameters. \citesupp{Carter2012} make a similar point in
the context of a simple data-generating process with time-invariant transition
probabilities. $\triangle$
\end{remark}

\section{Proofs of Supplemental Lemmas in Appendix \ref{app:mixing}}

\label{SM:mixing}

\begin{proof}[Proof of Lemma \ref{lem:couple}]
	For any $a,b$ in $\mathbb{S}$,
	\begin{align*}
		\bar{P}^{\nu}_{\theta} \left(  S_{l+1} = b \mid S_{l} = a , X^{j}_{-m} \right) = &  \frac{ \bar{P}^{\nu}_{\theta} \left(  S_{l+1} = b , S_{l} = a , X^{j}_{-m} \right)    }{ \sum_{c \in \mathbb{S}} \bar{P}^{\nu}_{\theta} \left(  S_{l+1} = c , S_{l} = a , X^{j}_{-m} \right) } \\
		%		= & \frac{ \Pr \left( X^{j}_{l+1} \mid   S_{l+1} = b , S_{l} = a , X^{l}_{-m} \right)  \Pr \left(  S_{l+1} = b \mid  S_{l} = a , X^{l}_{-m} \right)   }{ \sum_{c \in \mathbb{S}} \Pr \left( X^{j}_{l+1} \mid   S_{l+1} = c , S_{l} = a , X^{l}_{-m} \right)  \Pr \left(  S_{l+1} = c \mid  S_{l} = a , X^{l}_{-m} \right)  }  \\
		= & \frac{ \bar{P}^{\nu}_{\theta} \left( X^{j}_{l+1} \mid   S_{l+1} = b , S_{l} = a , X^{l}_{-m} \right)  Q_{\theta} \left(  X_{l} , a,b \right)   }{ \sum_{c \in \mathbb{S}} \bar{P}^{\nu}_{\theta} \left( X^{j}_{l+1} \mid   S_{l+1} = c , S_{l} = a , X^{l}_{-m} \right)  Q_{\theta} \left(  X_{l} , a,c \right)  }.		
	\end{align*}
	
	The expression $\bar{P}^{\nu}_{\theta} \left( X^{j}_{l+1} \mid   S_{l+1} = b , S_{l} = a , X^{l}_{-m} \right)$ equals $\bar{P}^{\nu}_{\theta} \left( X^{j}_{l+1} \mid   S_{l+1} = b , X_{l} \right)  $ because of the Markov property. The latter probability depends on the transitions of $X_{t+1}$ given $(X_{t},S_{t+1})$ and $S_{t+1}$ given $(X_{t},S_{t})$ for each $t \geq l+1$. Since these are the same for the process with $i=1$ and $i=2$ and the "original" process" $(S_{t},X_{t})_{t=-m}^{\infty}$, it follows that the last line of the previous display equals $	Pr_{\theta} \left(  \eta_{1,l+1} = b \mid \eta_{1,l} = a , X^{j}_{-m} \right)= 	Pr_{\theta} \left(  \eta_{2,l+1} = b \mid \eta_{2,l} = a , X^{j}_{-m} \right)$, as desired.		
\end{proof}

\begin{proof}[Proof of Lemma \ref{lem:ergo-eta}]
	Throughout this proof, we omit the dependence on $\theta$ in the probability terms and on other quantities. For any $a,c$ in $\mathbb{S}$,
	\begin{align*}
		& \left \Vert  \Pr \left(  \eta_{1,l+1} = \cdot  \mid \eta_{1,l} = a , X^{j}_{-m} \right) -  \Pr \left(  \eta_{2,l+1} = \cdot \mid \eta_{1,l} = c , X^{j}_{-m} \right)  \right \Vert_{1} \\
		\leq & \left \Vert  \Pr \left(  \eta_{1,l+1} = \cdot , \upsilon_{1,l} = 0 \mid \eta_{1,l} = a , X^{j}_{-m} \right)  - \Pr \left(  \eta_{2,l+1} = \cdot ,  \upsilon_{2,l} = 0 \mid \eta_{2,l} = c , X^{j}_{-m} \right)   \right \Vert_{1} \\
		& + \left \Vert  \Pr \left(  \eta_{1,l+1} = \cdot , \upsilon_{1,l} = 1 \mid \eta_{1,l} = a , X^{j}_{-m} \right)  - \Pr \left(  \eta_{2,l+1} = \cdot ,  \upsilon_{2,l} = 1 \mid \eta_{2,l} = c , X^{j}_{-m} \right)   \right \Vert_{1} \\
		\equiv & Term_{1} + Term_{2}
	\end{align*}	
	
	To bound the second term, note that
	\begin{align*}
		\Pr \left(  \eta_{1,l+1} = \cdot , \upsilon_{1,l} = 1 \mid \eta_{1,l} = a , X^{j}_{-m} \right)   = &  \Pr \left(  \eta_{1,l+1} = \cdot  \mid \upsilon_{1,l} = 1 , \eta_{1,l} = a , X^{j}_{-m} \right)\\
		& \times \Pr \left( \upsilon_{1,l} = 1 \mid \eta_{1,l} = a , X^{j}_{-m} \right).  \end{align*}
	It follows that $\Pr \left( \upsilon_{1,l} = 1 \mid \eta_{1,l} = a , X^{j}_{-m} \right) = \underline{q}(X_{l})$ because, given $X^{j}_{-m}$, $\upsilon_{1,l}$ is drawn independently according to a probability function that only depends on $X_{l}$ (in particular, it does not depend on $\eta_{1,l}$), and is given by $\underline{q}(X_{l})$. By some algebra, the Markov property, and the fact that, given $\upsilon_{1,l}=1$ and $X^{l}_{m}$, the random variable $\eta_{1,l+1}$ is independent of its past, it follows that $a \mapsto \Pr \left(  \eta_{1,l+1} = \cdot  \mid \upsilon_{1,l} = 1 , \eta_{1,l} = a , X^{j}_{-m} \right)$ is constant (i.e., does not depend on $\eta_{1,l} =a$). Thus, $a \mapsto \Pr \left(  \eta_{1,l+1} = \cdot , \upsilon_{1,l} = 1 \mid \eta_{1,l} = a , X^{j}_{-m} \right)$ is constant (i.e., does not depend on the value of $a$); since one can obtain the exact result for $c \mapsto \Pr \left(  \eta_{2,l+1} = \cdot , \upsilon_{2,l} = 1 \mid \eta_{2,l} = c , X^{j}_{-m} \right)$ and, moreover, the laws for $i=1$ and $i=2$ coincide (see the proof of Lemma \ref{lem:couple}),  it follows that $Term_{2}= 0$.
	
	To bound $Term_{1}$, it follows that from the previous arguments that
	\begin{align*}
		Term_{1} =  & \sum_{s \in \mathbb{S}} \left| \Pr \left(  \eta_{1,l+1} = s \mid \upsilon_{1,l} = 0, \eta_{1,l} = a , X^{j}_{-m} \right) -  \Pr \left(  \eta_{2,l+1} = s \mid \upsilon_{2,l} = 0, \eta_{2,l} = c , X^{j}_{-m} \right) \right| \\
		& \times (1-\underline{q}(X_{l})) \\
		\leq & 2 (1-\underline{q}(X_{l})),
	\end{align*}
	and thus the desired result follows.
\end{proof}

\section{Proofs of Supplementary Lemmas in Appendix \ref{app:LAR}}

\label{SM:LAR}

\subsection{Proofs of Lemmas \ref{lem:score_approxV2} and
\ref{lem:Hess-approx}}

\label{SM:LAR1}

In this section, we provide the proofs of Lemmas \ref{lem:score_approxV2} and
\ref{lem:Hess-approx}. To do this, we use a series of lemmas which we state
below (their proofs are relegated to the end of this section).

Henceforth, for any $j \geq m$, let
\begin{align}
\varrho(j,m) \equiv\left(  E_{\bar{P}_{\ast}^{\nu}} \left[  \prod_{i=m}^{j}
(1-\underline{q} (X_{i}))^{\frac{2a}{1-a}} \right]  \right)  ^{\frac{1-a}{2a}%
},
\end{align}
where the constant $a$ is the same as in Assumption \ref{ass:deriva-bdd}. We
also introduce the following notation: for any $\theta\in\Theta$, $(x^{\prime
},x,s)\mapsto\Gamma(x^{\prime}|x,s;\theta)\equiv\nabla_{\theta}\log p_{\theta
}(x,s,x^{\prime})$ and $(s,x,s)\mapsto\Lambda(s^{\prime}|s,x;\theta
)\equiv\nabla_{\theta}\log Q_{\theta}(x,s,s^{\prime})$. Furthermore, for any
$k\geq n$ and any $l\geq m$, let
\begin{align*}
&  \Phi_{\theta}(k,n,l,m)\equiv E_{\bar{P}_{\theta}^{\nu}}\left[  \sum
_{j=n}^{k}\Gamma(X_{j}|X_{j-1},S_{j};\theta)\mid X_{m}^{l}\right] \\
and~  &  \Psi_{\theta}(k,n,l,m)\equiv E_{\bar{P}_{\theta}^{\nu}}\left[
\sum_{j=n}^{k}\Lambda(S_{j}|S_{j-1},X_{j-1};\theta)\mid X_{m}^{l}\right]  .
\end{align*}

To state the first lemma, for any $k,T$ and $X_{k-T}^{k}$ and any $\theta$,
let
%\begin{align}
%\Delta _{k,k-T}(\theta )(X^{k}_{k-T}) \equiv & E_{\bar{P}_{\theta }^{\nu }}\left[
%\sum_{j=k-T}^{k-1}\Gamma (X_{j}|X_{j-1},S_{j};\theta )\mid X_{k-T}^{k}\right]
%+E_{\bar{P}_{\theta }^{\nu }}\left[ \sum_{l=k-T-1}^{k-1}\Lambda
%(S_{l}|S_{l-1},X_{l-1};\theta )\mid X_{k-T}^{k}\right]  \notag \\
%& -E_{\bar{P}_{\theta }^{\nu }}\left[ \sum_{j=k-T}^{k-1}\Gamma
%(X_{j}|X_{j-1},S_{j};\theta )\mid X_{k-T}^{k-1}\right] -E_{\bar{P}_{\theta }^{\nu
%	}}\left[ \sum_{l=k-T-1}^{k-1}\Lambda (S_{l}|S_{l-1},X_{l-1};\theta )\mid
%	X_{k-T}^{k-1}\right]  \notag \\
%	& +E_{\bar{P}_{\theta }^{\nu }}\left[ \Gamma (X_{k}|X_{k-1},S_{k};\theta )\mid
%	X_{k-T}^{k}\right] +E_{\bar{P}_{\theta }^{\nu }}\left[ \Lambda
%	(S_{k}|S_{k-1},X_{k-1};\theta )\mid X_{k-T}^{k}\right]  \notag \\
%	=& \sum_{j=k-T}^{k-1}E_{\bar{P}_{\theta }^{\nu }}\left[ \Gamma
%	(X_{j}|X_{j-1},S_{j};\theta )\mid X_{k-T}^{k}\right] -E_{\bar{P}_{\theta }^{\nu }}%
%	\left[ \Gamma (X_{j}|X_{j-1},S_{j};\theta )\mid X_{k-T}^{k-1}\right]  \notag
%	\\
%	& +\sum_{j=k-T}^{k-1}E_{\bar{P}_{\theta }^{\nu }}\left[ \Lambda
%	(S_{l}|S_{l-1},X_{l-1};\theta )\mid X_{k-T}^{k}\right] -E_{\bar{P}_{\theta }^{\nu
%		}}\left[ \Lambda (S_{l}|S_{l-1},X_{l-1};\theta )\mid X_{k-T}^{k-1}\right]
%		\notag \\
%		& +E_{\bar{P}_{\theta }^{\nu }}\left[ \Gamma (X_{k}|X_{k-1},S_{k};\theta )\mid
%		X_{k-T}^{k}\right] +E_{\bar{P}_{\theta }^{\nu }}\left[ \Lambda
%		(S_{k}|S_{k-1},X_{k-1};\theta )\mid X_{k-T}^{k}\right] ,  \label{eqn:Delta}
%		\end{align}%
\begin{align}
\Delta_{k,k-T}(\theta)(X^{k}_{k-T}) \equiv &  \Phi_{\theta}(k-1,k-T,k,k-T) +
\Psi_{\theta}(k-1,k-T-1,k,k-T)\nonumber\\
&  - \Phi_{\theta}(k-1,k-T,k-1,k-T) - \Psi_{\theta}%
(k-1,k-T-1,k-1,k-T)\nonumber\\
&  + \Phi_{\theta}(k,k,k,k-T) + \Psi_{\theta}(k,k,k,k-T)\nonumber\\
=  &  \sum_{j=k-T}^{k-1}E_{\bar{P}_{\theta}^{\nu}}\left[  \Gamma(X_{j}%
|X_{j-1},S_{j};\theta)\mid X_{k-T}^{k}\right]  -E_{\bar{P}_{\theta}^{\nu}}
\left[  \Gamma(X_{j}|X_{j-1},S_{j};\theta)\mid X_{k-T}^{k-1}\right]
\nonumber\\
&  +\sum_{j=k-T}^{k-1}E_{\bar{P}_{\theta}^{\nu}}\left[  \Lambda(S_{l}%
|S_{l-1},X_{l-1};\theta)\mid X_{k-T}^{k}\right]  -E_{\bar{P}_{\theta}^{\nu}%
}\left[  \Lambda(S_{l}|S_{l-1},X_{l-1};\theta)\mid X_{k-T}^{k-1}\right]
\nonumber\\
&  +E_{\bar{P}_{\theta}^{\nu}}\left[  \Gamma(X_{k}|X_{k-1},S_{k};\theta)\mid
X_{k-T}^{k}\right]  +E_{\bar{P}_{\theta}^{\nu}}\left[  \Lambda(S_{k}%
|S_{k-1},X_{k-1};\theta)\mid X_{k-T}^{k}\right]  , \label{eqn:Delta}%
\end{align}

The next lemma is analogous to the results in \citesupp{douc04} and
\citesupp{bickel98}, and uses ideas of missing data models.

\begin{lemma}
\label{lem:score-rep} Suppose Assumption \ref{ass:Theta-int} holds. Then, for
any $k,T \geq0$ and for any $\theta\in\Theta$,
\[
\nabla_{\theta}\log p_{k}^{\nu}(X_{k}|X_{k-T}^{k-1};\theta)=\Delta
_{k,k-T}(\theta)(X^{k}_{k-T})
\]
a.s.-$\bar{P}_{\ast}^{\nu}$.\footnote{When there is no risk of confusion, we
will omit the dependence of $\Delta_{k,k-T}$ on the data $X^{k}_{k-T}$.}
%where
%\begin{align*}
%\Delta_{k,k-T}(\theta) = & \sum_{j=k-T}^{k-1} \left\{ E_{P_{0}} \left[ \Gamma(\mathbf{x}_{j}|\mathbf{X}_{j-1},S_{j};\theta^{ S_{j}}_{x}) \mid \mathbf{X}^{k}_{k-T} \right] - E_{P_{0}} \left[  \Gamma(\mathbf{x}_{j}|\mathbf{X}_{j-1},S_{j};\theta^{ S_{j}}_{x}) \mid \mathbf{X}^{k-1}_{k-T} \right]\right\}\\
%& + \sum_{l=k-T-1}^{k-1} \left\{ E_{P_{0}} \left[  \Lambda(S_{l},S_{l+1}|\mathbf{z}_{l};\theta_{s}) \mid \mathbf{X}^{k}_{k-T} \right] -  E_{P_{0}} \left[  \Lambda(S_{l},S_{l+1}|\mathbf{z}_{l};\theta_{s}) \mid \mathbf{X}^{k-1}_{k-T} \right] \right\}\\
%& + E_{P_{0}} \left[ \Gamma(\mathbf{x}_{k}|\mathbf{X}_{k-1},S_{k};\theta^{ S_{k}}_{x}) \mid \mathbf{X}^{k}_{k-T} \right].
%\end{align*}

\end{lemma}

This Lemma characterizes the asymptotic behavior of the score functions; in
particular, it shows that they are well-approximated by $(\Delta_{t,-\infty
}(\theta_{\ast}))_{t,}$, which is to be defined below, but at this stage is
worth pointing out that it is stationary and ergodic; this last fact is
established in Lemma \ref{lem:Delta-mds} below.

\begin{lemma}
\label{lem:score_approx} Suppose Assumptions \ref{ass:BDD_Q},
\ref{ass:fx-ergo}, \ref{ass:Theta-int}, \ref{ass:deriva-bdd}(i) and
\ref{ass:q-sum-p} hold. Then:

(i) There exists a finite constant $C>0$ such that for any $k$ and $T \geq0$,
\[
||\Delta_{k,k-T}(\theta_{\ast})-\Delta_{k,-\infty}(\theta_{\ast})||_{L^{2}(
\bar{P}_{\ast}^{\nu})} \leq C \left(  \max\{\sum_{j=[k-T/2]}^{k-1}%
\varrho(j,k-T),\sum_{j=k-T}^{[k-T/2]-1}\varrho(k-1,j)\}\right)  ;
\]

(ii)
\[
\lim_{T\rightarrow\infty}\left\Vert T^{-1/2}\sum_{t=0}^{T}\{\Delta_{t,-\infty
}(\theta_{\ast})-\nabla_{\theta}\log p_{t}^{\nu}(\cdot|\cdot, \theta_{\ast
})\}\right\Vert _{L^{2}(\bar{P}_{\ast}^{\nu})}=0.
\]

\end{lemma}

\begin{lemma}
\label{lem:Delta-mds} Suppose Assumptions \ref{ass:BDD_Q}, \ref{ass:fx-ergo}
and \ref{ass:qbar-sum} hold. Then, $(\Delta_{t,-\infty}(\theta_{\ast
}))_{t=-\infty}^{\infty}$ is a stationary and ergodic $L^{2}(\bar{P}_{\ast
}^{\nu})$ process (under $\bar{P}_{\ast}^{\nu}$).
\end{lemma}

\begin{lemma}
\label{lem:Q-diff} Suppose Assumption \ref{ass:BDD_Q} holds. Then, there
exists a finite constant $L>0$ such that: \newline(i) For $-m\leq j<k$ and any
$\theta\in\Theta$,
\[
\left\Vert \bar{P}^{\nu}_{\theta}(S_{j}= \cdot|X_{-m}^{k})-\bar{P}^{\nu
}_{\theta}(S_{j}= \cdot|X_{-m}^{k-1}) \right\Vert _{1} \leq L \prod
_{i=j}^{k-1}(1-\underline{q}(X_{i})),
\]
a.s.-$\bar{P}_{\ast}^{\nu}$.

(ii) For $-n \leq-m \leq j < k$ and any $\theta\in\Theta$,
\begin{align*}
\left\Vert \bar{P}_{\theta}^{\nu} ( S_{j} = \cdot| X^{k}_{-m} ) - \bar
{P}_{\theta}^{\nu} ( S_{j} = \cdot| X^{k-1}_{-n} ) \right\Vert _{1} \leq L
\prod_{i=-m}^{j} (1-\underline{q}(X_{i})),
\end{align*}
a.s.-$\bar{P}_{\ast}^{\nu}$.
\end{lemma}

\begin{proof}[Proof of Lemma \ref{lem:score_approxV2}]
	Follows directly from Lemmas \ref{lem:score_approx} and \ref{lem:Delta-mds}.
\end{proof}

\begin{proof}[Proof of Lemma \protect\ref{lem:Hess-approx}]
					Lemma \ref{lem:Hess-approx} is analogous to Lemma~10 in \citesupp{bickel98}. The
					proof follows by their Lemma~9, which in turn holds by analogous steps to
					theirs and by invoking Lemma~\ref{lem:Q-diff} (which is analogous to their
					Lemma~7).
				\end{proof}

%The next lemma offers an intermediate result towards the LAN representation.
%\begin{lemma}
%\label{lem:LAN} Suppose Assumptions \ref{ass:BDD_Q}, \ref{ass:Theta-int}, %
%\ref{ass:deriva-bdd} and \ref{ass:q-sum-p} hold. Then, there exists a
%stationary and ergodic sequence $(\Delta _{t}(\theta_{\ast}))_{t}$ and a
%sequence of $\mathbb{R}^{q\times q}$-valued continuous functions $(\theta
%\mapsto \xi _{t}(\theta ))_{t}$ such that
%\begin{align*}
%\ell^{\nu} _{T}(X_{0}^{T},\theta_{\ast}+v)-\ell^{\nu} _{T}(X_{0}^{T},\theta
%_{0})=& v^{\prime} \left( T^{-1}\sum_{t=0}^{T}\Delta _{t}(\theta_{\ast})+o_{%
%\bar{P}_{\ast}^{\nu }}(T^{-1/2})\right) \\
%& +0.5v^{\prime }\left( T^{-1}\sum_{t=0}^{T}\int_{0}^{1}\xi _{t}(\theta
%_{0}+sv)ds+o_{\bar{P}_{\ast}^{\nu }}(1)\right) v
%\end{align*}%
%for any $v$ such that $\theta_{\ast}+v\in B(\delta ,\theta_{\ast})$.
%\end{lemma}
%\begin{proof}
%	See Supplemental Material XXX.
%\end{proof}		

\subsubsection{Proofs of Lemmas}

Throughout this section, in cases where the expectations are taken with
respect to $\bar{P}_{\ast}^{\nu}$, we omit the probability from the notation
and simply use $E[\cdot]$.

\begin{proof}[Proof of Lemma \protect\ref{lem:score-rep}]
	By \citesupp[p. 227]{louis1982finding},
	\begin{align*}
	\nabla _{\theta }\log p_{k}^{\nu }(X_{k}|X_{k-T}^{k-1};\theta )=& \nabla
	_{\theta }\log p_{k}^{\nu }(X_{k-T}^{k};\theta )-\nabla _{\theta }\log
	p_{k-1}^{\nu }(X_{k-T}^{k-1};\theta ) \\
	=& E_{\bar{P}_{\theta }^{\nu }}[\nabla _{\theta }\log p_{k}^{\nu
	}(X_{k-T}^{k},S_{k-T}^{k};\theta )\mid X_{k-T}^{k}] \\
	& -E_{\bar{P}_{\theta }^{\nu }}[\nabla _{\theta }\log p_{k-1}^{\nu
	}(X_{k-T}^{k-1},S_{k-T}^{k-1};\theta )\mid X_{k-T}^{k-1}].
	\end{align*}%
	(Note that the expectation is with respect to $S_{k-T}^{k}$, which takes
	finitely many values; thus interchanging differentiation and integration is
	allowed).
	
	Since $p_{k}^{\nu }(X_{k-T}^{k},S_{k-T}^{k};\theta )=p_{\theta
	}(X_{k-1},S_{k},X_{k})Q_{\theta }(X_{k-1},S_{k-1},S_{k})\times
	p_{k}^{\nu }(X_{k-T}^{k-1},S_{k-T}^{k-1};\theta )$ (and an analogous result
	holds for $p_{k}^{\nu }(X_{k-T}^{k-1},S_{k-T}^{k-1};\theta )$), it follows
	that
	\begin{align*}
	& \nabla _{\theta }\log p_{k}^{\nu }(X_{k}|X_{k-T}^{k-1};\theta ) \\
	=& E_{\bar{P}_{\theta }^{\nu }}\left[ \sum_{j=k-T}^{k}\nabla _{\theta }\log
	p_{\theta}(X_{j-1},S_{j},X_{j})\mid X_{k-T}^{k}\right] +E_{\bar{P}_{\theta }^{\nu
		}}\left[ \sum_{j=k-T}^{k}\nabla _{\theta }\log Q_{\theta
	}(X_{j-1},S_{j-1},S_{j})\mid X_{k-T}^{k}\right] \\
	& -E_{\bar{P}_{\theta }^{\nu }}\left[ \sum_{j=k-T}^{k-1}\nabla _{\theta }\log
	p_{\theta}(X_{j-1},S_{j},X_{j})\mid X_{k-T}^{k-1}\right] -E_{\bar{P}_{\theta
		}^{\nu }}\left[ \sum_{j=k-T}^{k-1}\nabla _{\theta }\log Q_{\theta
	}(X_{j-1},S_{j-1},S_{j})\mid X_{k-T}^{k-1}\right] \\
	=& E_{\bar{P}_{\theta }^{\nu }}\left[ \sum_{j=k-T}^{k-1}\Gamma
	(X_{j}|X_{j-1},S_{j};\theta )\mid X_{k-T}^{k}\right] +E_{\bar{P}_{\theta }^{\nu }}%
	\left[ \sum_{l=k-T-1}^{k-1}\Lambda (S_{l}|S_{l-1},X_{l-1};\theta )\mid
	X_{k-T}^{k}\right] \\
	& -E_{\bar{P}_{\theta }^{\nu }}\left[ \sum_{j=k-T}^{k-1}\Gamma
	(X_{j}|X_{j-1},S_{j};\theta )\mid X_{k-T}^{k-1}\right] -E_{\bar{P}_{\theta }^{\nu
		}}\left[ \sum_{l=k-T-1}^{k-1}\Lambda (S_{l}|S_{l-1},X_{l-1};\theta )\mid
		X_{k-T}^{k-1}\right] \\
		& +E_{\bar{P}_{\theta }^{\nu }}\left[ \Gamma (X_{k}|X_{k-1},S_{k};\theta )\mid
		X_{k-T}^{k}\right] +E_{\bar{P}_{\theta }^{\nu }}\left[ \Lambda
		(S_{k}|S_{k-1},X_{k-1};\theta )\mid X_{k-T}^{k}\right] .
		\end{align*}
	\end{proof}

%-----------------------------------------------------------

The proof of Lemma \ref{lem:score_approx} requires the following lemma.

\begin{lemma}
\label{lem:G-L-BDD} Suppose that Assumptions \ref{ass:BDD_Q} and
\ref{ass:deriva-bdd}(i) hold. Then, there exists a finite constant $C >0$ such that:

(i) for any $-n\leq-m\leq-m^{\prime}\leq l\leq k$,
%		\begin{align*}
%		& \left\Vert E_{\bar{P}_{\theta_{\ast}}^{\nu }}\left[	\sum_{j=-m^{\prime
%			}}^{l}\Gamma (X_{j}|X_{j-1},S_{j};\theta_{\ast})\mid X_{-m}^{k}\right]
%			-E_{\bar{P}_{\theta_{\ast}}^{\nu }}\left[ \sum_{j=-m^{\prime }}^{l}\Gamma
%			(X_{j}|X_{j-1},S_{j};\theta_{\ast})\mid X_{-n}^{k}\right] \right\Vert _{L^{2}(%
%				\bar{P}_{\ast}^{\nu })} \\
%			=& O\left( \sum_{j=-m^{\prime }}^{l}\varrho (j,-m)\right) ;
%			\end{align*}%
\begin{align*}
\left\Vert \Phi_{\theta_{\ast}}(l,-m^{\prime},k,-m) - \Phi_{\theta_{\ast}%
}(l,-m^{\prime},k,-n) \right\Vert _{L^{2}(\bar{P}_{\ast}^{\nu})} \leq C
\left(  \sum_{j=-m^{\prime}}^{l}\varrho(j,-m)\right)  ;
\end{align*}

(ii) for any $-m\leq-m^{\prime}<l\leq k-1$,
\begin{align*}
\left\Vert \Phi_{\theta_{\ast}}(l,-m^{\prime},k,-m) - \Phi_{\theta_{\ast}%
}(l,-m^{\prime},k-1,-m) \right\Vert _{L^{2}(\bar{P}_{\ast}^{\nu})} \leq C
\left(  \sum_{j=-m^{\prime}}^{l}\varrho(k-1,j)\right)  ;
\end{align*}
%			\begin{align*}
%			& \left\Vert E_{\bar{P}_{\theta_{\ast}}^{\nu }}\left[ \sum_{j=-m^{\prime
%				}}^{l}\Gamma (X_{j}|X_{j-1},S_{j};\theta_{\ast})\mid X_{-m}^{k}\right]
%				-E_{\bar{P}_{\theta_{\ast}}^{\nu }}\left[ \sum_{j=-m^{\prime }}^{l}\Gamma
%				(X_{j}|X_{j-1},S_{j};\theta_{\ast})\mid X_{-m}^{k-1}\right] \right\Vert
%				_{L^{2}(\bar{P}_{\ast}^{\nu })} \\
%				=& O\left( \sum_{j=-m^{\prime }}^{l}\varrho (k-1,j)\right) ;
%				\end{align*}

(iii) for any $-n\leq-m\leq-m^{\prime}<l\leq k$,
\begin{align*}
\left\Vert \Psi_{\theta_{\ast}}(l,-m^{\prime},k,-m) - \Psi_{\theta_{\ast}%
}(l,-m^{\prime},k,-n) \right\Vert _{L^{2}(\bar{P}_{\ast}^{\nu})} \leq C
\left(  \sum_{j=-m^{\prime}}^{l}\varrho(j-1,-m)\right)  ;
\end{align*}
%				\begin{align*}
%				& \left\Vert E_{\bar{P}_{\theta _{0}}^{\nu }}\left[ \sum_{j=-m^{\prime
%					}}^{l}\Lambda (S_{j}|S_{j-1},X_{j-1};\theta _{0})\mid X_{-m}^{k}\right]
%					-E_{\bar{P}_{\theta _{0}}^{\nu }}\left[ \sum_{j=-m^{\prime }}^{l}\Lambda
%					(S_{j}|S_{j-1},X_{j-1};\theta _{0})\mid X_{-n}^{k}\right] \right\Vert
%					_{L^{2}(\bar{P}_{\ast}^{\nu })} \\
%					=& O\left( \sum_{j=-m^{\prime }}^{l}\varrho (j-1,-m)\right) ;
%					\end{align*}

(iv) for any $-m\leq-m^{\prime}<l\leq k-1$,
\begin{align*}
\left\Vert \Psi_{\theta_{\ast}}(l,-m^{\prime},k,-m) - \Psi_{\theta_{\ast}%
}(l,-m^{\prime},k-1,-m) \right\Vert _{L^{2}(\bar{P}_{\ast}^{\nu})} \leq C
\left(  \sum_{j=-m^{\prime}}^{l}\varrho(k-1,j)\right)  .
\end{align*}
%					\begin{align*}
%					& \left\Vert E_{\bar{P}_{\theta_{\ast}}^{\nu }}\left[ \sum_{j=-m^{\prime
%						}}^{l}\Lambda (S_{j}|S_{j-1},X_{j-1};\theta_{\ast})\mid X_{-m}^{k}\right]
%						-E_{\bar{P}_{\theta_{\ast}}^{\nu }}\left[ \sum_{j=-m^{\prime }}^{l}\Lambda
%						(S_{j}|S_{j-1},X_{j-1};\theta_{\ast})\mid X_{-m}^{k-1}\right] \right\Vert
%						_{L^{2}(\bar{P}_{\ast}^{\nu })} \\
%						=& O\left( \sum_{j=-m^{\prime }}^{l}\varrho (k-1,j)\right) .
%						\end{align*}

\end{lemma}

\begin{proof}[Proof of Lemma \protect\ref{lem:G-L-BDD}]
						Throughout the proof we omit the dependence of $E[\cdot]$ on $\bar{P}_{\theta_{\ast}}^{\nu} $. Also, let $L$ denote the constant in Lemma \ref{lem:Ker-rev}.
						
						\textbf{Part (i).} Observe that, for any $j\leq k$,
						\begin{align*}
						& \left\Vert E\left[ \Gamma (X_{j}|X_{j-1},S_{j};\theta_{\ast})\mid X_{-m}^{k}%
						\right] -E\left[ \Gamma (X_{j}|X_{j-1},S_{j};\theta_{\ast})\mid X_{-n}^{k}%
						\right] \right\Vert \\
						& =\left\Vert \sum_{a\in \mathbb{S}}\Gamma (X_{j}|X_{j-1},a;\theta
						_{\ast})\{ \bar{P}_{\theta_{\ast}}^{\nu}  (S_{j}=a\mid X_{-m}^{k})- \bar{P}_{\theta_{\ast}}^{\nu}  (S_{j}=a\mid X_{-n}^{k})\}\right\Vert
						%						\\
						%						\leq & \sqrt{\sum_{a\in \mathbb{S}}||\Gamma (X_{j}|X_{j-1},a;\theta
						%							_{0})||^{2}}\sqrt{\sum_{a\in \mathbb{S}}\{\Pr (S_{j}=a\mid X_{-m}^{k})-\Pr
						%							(S_{j}=a\mid X_{-n}^{k})\}^{2}},
						\\
						\leq & \max_{a\in \mathbb{S}} ||\Gamma (X_{j}|X_{j-1},a;\theta
						_{\ast})|| \left \Vert \bar{P}_{\theta_{\ast}}^{\nu}  (S_{j}= \cdot \mid X_{-m}^{k})- \bar{P}_{\theta_{\ast}}^{\nu} (S_{j}=\cdot \mid X_{-n}^{k}) \right \Vert_{1}.
						\end{align*}%
						By Lemma \ref{lem:Q-diff}(ii),
						% it follows that
						%						\begin{equation*}
						%						\left( \Pr (S_{j}=a\mid X_{-m}^{k})-\Pr (S_{j}=a\mid X_{-n}^{k})\right)
						%						^{2}\leq \prod_{i=-m}^{j}(1-\underline{q}(X_{i}))^{2},
						%						\end{equation*}%
						%						which in turn implies that
						\begin{align*}
						& \left \Vert  E\left[ \Gamma (X_{j}|X_{j-1},S_{j};\theta_{\ast})\mid X_{-m}^{k}\right] -E%
						\left[ \Gamma (X_{j}|X_{j-1},S_{j};\theta_{\ast})\mid X_{-n}^{k}\right] \right \Vert \\
						\leq & L \max_{a\in \mathbb{S}}||\Gamma
						(X_{j}|X_{j-1},a;\theta_{\ast})|| \prod_{i=-m}^{j}(1-\underline{q}(X_{i})),
						\end{align*}
						Therefore, by the H\"{o}lder inequality, it follows that, for $a^{-1} + b^{-1} = 1$ (with $a$ as in Assumption \ref{ass:deriva-bdd}),
						\begin{align*}
						& \left\Vert \sum_{j=-m^{\prime}}^{l} \left\{ E\left[ \Gamma (X_{j}|X_{j-1},S_{j};\theta_{\ast})\mid X_{-m}^{k}%
						\right] -E\left[ \Gamma (X_{j}|X_{j-1},S_{j};\theta_{\ast})\mid X_{-n}^{k}%
						\right] \right\} \right\Vert _{L^{2}(\bar{P}_{\ast}^{\nu })} \\
						& \leq \sum_{j=-m^{\prime}}^{l} \left\Vert E\left[ \Gamma (X_{j}|X_{j-1},S_{j};\theta_{\ast})\mid X_{-m}^{k}%
						\right] -E\left[ \Gamma (X_{j}|X_{j-1},S_{j};\theta_{\ast})\mid X_{-n}^{k}%
						\right] \right\Vert _{L^{2}(\bar{P}_{\ast}^{\nu })} \\
						\leq & L \left( \sum_{a\in \mathbb{S}}E_{\bar{P}_{\ast}^{\nu }} \left[ ||\Gamma (X_{1}|X_{0},a;\theta_{\ast})||^{2 a }\right] \right)^{1/(2a)} \sum_{j=-m^{\prime}}^{l}  \left( E_{\bar{P}_{\ast}^{\nu }}\left[ \prod_{i=-m}^{j}(1-\underline{q}(X_{i}))^{2b}\right] \right)^{1/(2b)},
						\end{align*}
						where the second line follows from the triangle inequality and the third follows from stationarity (Lemma \ref{lem:sta-ergo}). The
						fact that $\Gamma (X_{1}|X_{0},a;\theta_{\ast})=\nabla _{\theta }\log
						p_{\theta_{\ast}}(X_{0},a,X_{1})$, Assumption~\ref{ass:deriva-bdd}(i) and definition of $\varrho$ imply
						the desired result.
						
						\bigskip
						
						\textbf{Part (ii).} Follows from analogous calculations to those in part (i) and Lemma \ref{lem:Q-diff}(i).
						
						\bigskip
						
						\textbf{Parts (iii) and (iv).} We only work out part (iii) since (iv) is
						analogous.					
					
					   Observe that
						\begin{align*}
						& \left\Vert E\left[ \Lambda (S_{j}|S_{j-1},X_{j-1};\theta_{\ast})\mid
						X_{-m}^{k}\right] -E\left[ \Lambda (S_{j}|S_{j-1},X_{j-1};\theta_{\ast})\mid
						X_{-n}^{k}\right] \right\Vert \\
						 = &\left\Vert \sum_{(a,b)\in \mathbb{S}^{2}}\Lambda (b|a,X_{j-1};\theta
						_{\ast})\{\bar{P}_{\theta^{\ast}}^{\nu} (S_{j}=b,S_{j-1}=a\mid X_{-m}^{k})-\bar{P}_{\theta^{\ast}}^{\nu} (S_{j}=b,S_{j-1}=a\mid
						X_{-n}^{k})\}\right\Vert \\
						\leq & \sum_{(a,b)\in \mathbb{S}^{2}}||\Lambda (b|a,X_{j-1};\theta
							_{\ast})|| \left | \bar{P}_{\theta^{\ast}}^{\nu} (S_{j}=b,S_{j-1}=a\mid
							X_{-m}^{k})-\bar{P}_{\theta^{\ast}}^{\nu} (S_{j}=b,S_{j-1}=a\mid X_{-n}^{k}) \right|\\
%						= & \sum_{(a,b)\in \mathbb{S}^{2}}||\Lambda (b|a,X_{j-1};\theta
%						_{0})|| \\
%						& \times \left | \bar{P}_{\theta^{\ast}}^{\nu} (S_{j}=b \mid S_{j-1}=a,
%						X_{-m}^{k})\bar{P}_{\theta^{\ast}}^{\nu} (S_{j-1}=a\mid X_{-m}^{k}) -\bar{P}_{\theta^{\ast}}^{\nu} (S_{j}=b \mid S_{j-1}=a, X_{-n}^{k}) \bar{P}_{\theta^{\ast}}^{\nu} (S_{j-1}=a\mid X_{-n}^{k}) \right|\\	
						= & \sum_{(a,b)\in \mathbb{S}^{2}}||\Lambda (b|a,X_{j-1};\theta
						_{\ast})|| \bar{P}_{\theta^{\ast}}^{\nu} (S_{j}=b \mid S_{j-1}=a ,	X_{j-1}^{k}) \left | \bar{P}_{\theta^{\ast}}^{\nu} (S_{j-1}=a\mid X_{-m}^{k}) - \bar{P}_{\theta^{\ast}}^{\nu} (S_{j-1}=a\mid X_{-n}^{k}) \right|		
						\end{align*}
						where the last line follows from the Markov property of the model and the
						fact that $-m\leq j$. Thus
	\begin{align*}
	& \left\Vert E\left[ \Lambda (S_{j}|S_{j-1},X_{j-1};\theta_{\ast})\mid
	X_{-m}^{k}\right] -E\left[ \Lambda (S_{j}|S_{j-1},X_{j-1};\theta_{\ast})\mid
	X_{-n}^{k}\right] \right\Vert \\
	\leq & \max_{a,b} ||\Lambda (b|a,X_{j-1};\theta
	_{\ast})|| \left \Vert  \bar{P}_{\theta^{\ast}}^{\nu} (S_{j-1}= \cdot \mid X_{-m}^{k}) - \bar{P}_{\theta^{\ast}}^{\nu} (S_{j-1}= \cdot \mid X_{-n}^{k}) \right \Vert_{1} \\
	\leq & L \max_{a,b} ||\Lambda (b|a,X_{j-1};\theta_{\ast})|| \prod_{i=-m}^{j-1} (1-\underline{q}(X_{i}))		
	\end{align*}
	where the last line follows from Lemma \ref{lem:Q-diff}(ii). Hence	
	\begin{align*}
		 & \left\Vert \sum_{j=-m^{\prime}}^{l} \{ E\left[ \Lambda (S_{j}|S_{j-1},X_{j-1};\theta_{\ast})\mid
		 X_{-m}^{k}\right] -E\left[ \Lambda (S_{j}|S_{j-1},X_{j-1};\theta_{\ast})\mid
		 X_{-n}^{k}\right] \} \right\Vert_{L^{2}(\bar{P}^{\nu}_{\ast})} \\
		 \leq & L \left( \sum_{(a,b)\in \mathbb{S}^{2}} E_{\bar{P}_{\ast}^{\nu }} \left[ ||\Lambda (b|a,X_{0};\theta_{\ast})||^{2a }\right] \right)^{1/(2a)} \sum_{j=-m^{\prime}}^{l}  \left( E_{P_{\ast}^{\nu }}\left[ \prod_{i=-m}^{j-1}(1-\underline{q}(X_{i}))^{2b}\right] \right)^{1/(2b)},
	\end{align*}				
	 so the result follows from Assumption \ref{ass:deriva-bdd}(i).
\end{proof}

\begin{proof}[Proof of Lemma \protect\ref{lem:score_approx}]
						Throughout the proof we denote $||.||_{L^{2}(\bar{P}_{\ast}^{\nu
							})}$ as $||.||_{L^{2}}$. Also, we use $\Phi$ and $\Psi$ to denote $\Phi_{\theta_{\ast}}$ and $\Psi_{\theta_{\ast}}$, resp.
						
						\textbf{Part (i):} Observe that $\Phi
						(k-1,k-T,l,k-T)=\Phi (k-1,[k-T/2],l,k-T)+\Phi ([k-T/2]-1,k-T,l,k-T)$ and an
						analogous result holds for $\Psi $. Therefore, by the definition of $\Delta
						_{k,k-T}$ and analogous calculations to those in \citesupp[pp.
						1624--1626]{bickel98},
						\begin{align*}
						& ||\Delta _{k,k-T}(\theta_{\ast})-\Delta _{k,-\infty }(\theta_{\ast})||_{L^{2}}
						\\
						\leq & \left\Vert \Phi (k-1,[k-T/2],k,k-T)-\Phi (k-1,[k-T/2],k,-\infty
						)\right\Vert _{L^{2}} \\
						& +\left\Vert \Phi (k-1,[k-T/2],k-1,k-T)-\Phi (k-1,[k-T/2],k-1,-\infty
						)\right\Vert _{L^{2}} \\
						& +\left\Vert \Phi ([k-T/2]-1,k-T,k,k-T)-\Phi
						([k-T/2]-1,k-T,k-1,k-T)\right\Vert _{L^{2}} \\
						& +\left\Vert \Psi (k-1,[k-T/2],k,k-T)-\Psi (k-1,[k-T/2],k,-\infty
						)\right\Vert _{L^{2}} \\
						& +\left\Vert \Psi (k-1,[k-T/2],k-1,k-T)-\Psi (k-1,[k-T/2],k-1,-\infty
						)\right\Vert _{L^{2}} \\
						& +\left\Vert \Psi ([k-T/2]-1,k-T-1,k,k-T)-\Psi
						([k-T/2]-1,k-T-1,k-1,k-T)\right\Vert _{L^{2}} \\
						& +\left\Vert \Phi (k,k,k,k-T)-\Phi (k,k,k,-\infty )\right\Vert
						_{L^{2}}+\left\Vert \Psi (k,k,k,k-T)-\Psi (k,k,k,-\infty )\right\Vert
						_{L^{2}} \\
						\equiv & \sum_{i=1}^{8}Term_{i}.
						\end{align*}
						
						Terms $1$ and $2$ are analogous of the form
						\begin{align*}
						& \left\Vert \Phi (k-1,[k-T/2],k,k-T)-\Phi (k-1,[k-T/2],k,-\infty
						)\right\Vert _{L^{2}} \\
						=& \left\Vert \sum_{j=[k-T/2]}^{k-1}E_{\bar{P}_{\theta_{\ast}}^{\nu }}[\Gamma
						(X_{j}\mid X_{j-1},S_{j};\theta_{\ast})\mid
						X_{k-T}^{m}]-\sum_{j=[k-T/2]}^{k-1}E_{\bar{P}_{\theta_{\ast}}^{\nu }}[\Gamma
						(X_{j}\mid X_{j-1},S_{j};\theta_{\ast})\mid X_{-\infty }^{m}]\right\Vert
						_{L^{2}}
						\end{align*}%
						for $m\in \{k,k-1\}$. By Lemma \ref{lem:G-L-BDD}(i), for $i\in \{1,2\}$, $%
						Term_{i} \precsim \left( \sum_{j=[k-T/2]}^{k-1}\varrho (j,k-T)\right) $. $Term_{7}$ has the same bound by
						analogous calculations.
						
						The $Term_{3}$ is of the form
						\begin{align*}
						\left \Vert \sum_{j=k-T}^{[k-T/2]-1} E_{\bar{P}^{\nu}_{\theta_{\ast}}}[\Gamma(X_{j}
						\mid X_{j-1},S_{j};\theta_{\ast} ) \mid X^{k}_{k-T}] - \sum_{j=k-T}^{[k-T/2]-1}
						E_{\bar{P}^{\nu}_{\theta_{\ast}}}[\Gamma(X_{j} \mid X_{j-1},S_{j};\theta_{\ast} ) \mid
						X^{k-1}_{k-T}] \right \Vert_{L^{2}}.
						\end{align*}
						By Lemma \ref{lem:G-L-BDD}(ii), $Term_{3}  \precsim \left( \sum_{j=k-T}^{[k-T/2]-1}
						\varrho(k-1,j) \right)$. The term $Term_{8}$ has the same bound by analogous
						calculations.
						
						The terms $Term_{4}$ and $Term_{5}$ are of the form
						\begin{align*}
						\left \Vert \sum_{j=[k-T/2]}^{k-1} E_{\bar{P}^{\nu}_{\theta_{\ast}}}[\Lambda(S_{j}
						\mid S_{j-1},X_{j-1};\theta_{\ast} ) \mid X^{m}_{k-T}] - \sum_{j=[k-T/2]}^{k-1}
						E_{\bar{P}^{\nu}_{\theta_{\ast}}}[\Lambda(S_{j} \mid S_{j-1},X_{j-1};\theta_{\ast} )
						\mid X^{m}_{-\infty}] \right \Vert_{L^{2}}
						\end{align*}
						and by Lemma \ref{lem:G-L-BDD}(iii) are bounded by $C \left(
						\sum_{j=[k-T/2]}^{k-1} \varrho(j-1,k-T) \right)$ for some universal constant $C>0$.
						
						Finally, by analogous calculations to those for $Term_{3}$, $Term_{6}$ is
						bounded by $C \left( \sum_{j=k-T}^{[k-T/2]-1} \varrho(k-1,j) \right)$ by
						Lemma \ref{lem:G-L-BDD}(iv).\newline
						
						\textbf{Part (ii).} By part (i) and Lemma \ref{lem:score-rep},
						\begin{align*}
						 \left\Vert T^{-1/2}\sum_{t=1}^{T}\{\Delta _{t,-\infty }(\theta
						_{\ast})-\nabla _{\theta }\log p_{t}^{\nu }(\cdot |\cdot;\theta
						_{\ast})\}\right\Vert _{L^{2}} \leq & T^{-1/2}\sum_{t=1}^{T}\left\Vert \{\Delta _{t,-\infty }(\theta
						_{\ast})-\Delta _{t,0}(\theta_{\ast})\}\right\Vert _{L^{2}} \\
						\precsim & \left(  T^{-1/2}\sum_{t=1}^{T}\sum_{j=[t/2]}^{t-1}\varrho
						(j,0)+T^{-1/2}\sum_{t=1}^{T}\sum_{j=0}^{[t/2]-1}\varrho (t,j) \right).
						\end{align*}%
						By Kronecker's lemma, it suffices to show that
						\begin{equation}
						\sum_{t=1}^{T}t^{-1/2}\sum_{j=[t/2]}^{t-1}\varrho (j,0)~\quad \mathrm{%
							and\quad }~\sum_{t=1}^{T}t^{-1/2}\sum_{j=0}^{[t/2]-1}\varrho (t,j)
						\label{eqn:score_approx-1}
						\end{equation}%
						are bounded uniformly in $T$, where, recall, $\varrho (j,k)\equiv \left( E\left[
							\prod_{i=k}^{j}(1-\underline{q}(X_{i}))^{\frac{2a}{1-a}}\right] \right)^{\frac{1-a}{2a}}$.
						
						Moreover, $j\mapsto \varrho (j,k)$ is non-increasing and $k\mapsto \varrho
						(j,k)$ is non-decreasing since $1-\underline{q}(\cdot)\leq 1$. By Assumption \ref%
						{ass:q-sum-p}, $(\varrho (j,0))_{j}$ is $p$-summable with $p<2/3$, thus $%
						\lim_{j\rightarrow \infty }\varrho (j,0)^{p}j=0$ (if not, then $\varrho
						(j,0)>c/j^{1/p}$ for some $c>0$ and all $j$ above certain point and this
						violates the assumption). Hence,
						\begin{equation*}
						\sum_{j=[t/2]}^{t-1}\varrho (j,0)<\sum_{j=[t/2]}^{t-1}\frac{1}{j^{1/p}}\leq
						\int_{\lbrack t/2]+1}^{t}x^{-1/p}dx\leq \frac{p}{1-p}(t/2)^{1-1/p},
						\end{equation*}%
						for all $t\geq \tau $ and some $\tau >0$, and this implies that, for
						some constant $const>0$,
						\begin{equation*}
						\sum_{t=1}^{T}t^{-1/2}\sum_{j=[t/2]}^{t-1}\varrho (j,0)\leq C(\tau
						)+const\times \sum_{t=\tau +1}^{T}\frac{p}{1-p}t^{1-1/p-1/2}\leq C<\infty ,
						\end{equation*}%
						because $1-1/p-1/2<-1\Leftrightarrow p<2/3$ ($C$ is a finite constant, which
						may depend on $\tau $).
						
						By stationarity of $X_{-\infty}^{\infty}$ (Lemma \ref{lem:sta-ergo}),
						\begin{align*}
						\sum_{j=0}^{[t/2]-1}\varrho (t,j)=& \left( E\left[
						\prod_{i=0}^{t}(1-\underline{q}(X_{i}))^{\frac{2a}{1-a}}\right] \right)^{\frac{1-a}{2a}}+\left( E\left[
						\prod_{i=1}^{t}(1-\underline{q}(X_{i}))^{\frac{2a}{1-a}}\right] \right)^{\frac{1-a}{2a}}\\
						& +...+\left( E\left[
						\prod_{i=[t/2]-1}^{t}(1-\underline{q}(X_{i}))^{\frac{2a}{1-a}}\right] \right)^{\frac{1-a}{2a}} \\
						=& \left( E\left[
						\prod_{i=0}^{t}(1-\underline{q}(X_{i}))^{\frac{2a}{1-a}}\right] \right)^{\frac{1-a}{2a}}+\left( E\left[
						\prod_{i=0}^{t-1}(1-\underline{q}(X_{i}))^{\frac{2a}{1-a}}\right] \right)^{\frac{1-a}{2a}}\\
						& +...+\left( E\left[
						\prod_{i=0}^{[t/2]-1}(1-\underline{q}(X_{i}))^{\frac{2a}{1-a}}\right] \right)^{\frac{1-a}{2a}} \\
						=& \sum_{j=0}^{[t/2]+1}\varrho (t-j,0).
						\end{align*}
						
						Thus $\sum_{j=0}^{[t/2]+1}\varrho (t-j,0)\leq \sum_{j=0}^{[t/2]+1}\frac{1}{%
							(t-j)^{1/p}}\leq \int_{\lbrack t/2]-1}^{t}\frac{1}{u^{1/p}}du$ and by our
						previous calculation the result follows. Thus, the terms in (\ref%
						{eqn:score_approx-1}) are uniformly bounded.
					\end{proof}

%-------------------------------------------------------

\begin{proof}[Proof of Lemma \protect\ref{lem:Delta-mds}]
						It is easy to see that $\Delta _{t,-\infty }(\theta_{\ast})$ is adapted to the
						filtration associated with the $\sigma $-algebra generated by $X_{-\infty
						}^{t}$. Since $X_{-\infty }^{\infty }$ is stationary and ergodic (by Lemma \ref%
						{lem:sta-ergo}), so is $(\Delta _{t,-\infty }(\theta_{\ast}))_{t=-\infty
						}^{\infty }$.
						%Ergodicity of $(\Delta _{t,-\infty }(\theta_{\ast}))_{t=-\infty}^{\infty }$ follows from Lemma \ref{lem:sta-ergo}.
					\end{proof}

To prove Lemma \ref{lem:Q-diff}, we need the following results.

\begin{lemma}
\label{lem:Ker-rev} Suppose Assumption \ref{ass:BDD_Q} holds. Then, there
exists a finite constant $L>0$, such that, for any $-m\leq j < n \leq k$ and
any $\theta\in\Theta$,
\[
\max_{b,c}\left\Vert \bar{P}^{\nu}_{\theta}\left(  S_{j}=\cdot|S_{n}%
=b,X_{-m}^{k}\right)  -\bar{P}^{\nu}_{\theta}\left(  S_{j}=\cdot
|S_{n}=c,X_{-m}^{k}\right)  \right\Vert _{1} \leq L \prod_{i=j}^{n}%
(1-q(X_{i}))
\]
a.s.-$\bar{P}_{\ast}^{\nu}$.
\end{lemma}

\begin{lemma}
\label{lem:Rev-MP} For any $-m<i<l\leq r \leq n$, let $\mathbf{S}_{l}%
^{r}\equiv(S_{l},...,S_{r})$. Then, for any $\theta\in\Theta$,
\begin{align*}
\bar{P}_{\theta}^{\nu} \left(  S_{i}|\mathbf{S}_{l}^{r},X_{-m+1}^{n}\right)  =
\bar{P}_{\theta}^{\nu} \left(  S_{i}|S_{l},X_{-m+1}^{n}\right)  ,
\end{align*}
i.e., the Markov property holds backward in time.
\end{lemma}

\begin{proof}[Proof of Lemma \protect\ref{lem:Ker-rev}]				
						Observe that, for
						any $b,c\in \mathbb{S}^{2}$,
						\begin{align*}
						& \left \Vert \bar{P}_{\theta}^{\nu} \left( S_{j}= \cdot|S_{n}=b,X_{-m}^{k}\right) -\bar{P}_{\theta}^{\nu} \left(
						S_{j}=\cdot|S_{n}=c,X_{-m}^{k}\right) \right\Vert_{1} \\
						=& \left\Vert \sum_{s\in \mathbb{S}}\bar{P}_{\theta}^{\nu} \left(
						S_{j}=\cdot|S_{j+1}=s,X_{-m}^{j+1}\right) \left( \bar{P}_{\theta}^{\nu} \left(
						S_{j+1}=s|S_{n}=b,X_{-m}^{k}\right) - \bar{P}_{\theta}^{\nu} \left(
						S_{j+1}=s|S_{n}=c,X_{-m}^{k}\right) \right) \right\Vert_{1} \\
						\leq & \alpha_{\theta,j+1,j}(X^{k}_{-m}) \left \Vert \bar{P}_{\theta}^{\nu} \left( S_{j+1}= \cdot|S_{n}=b,X_{-m}^{k}\right) -\bar{P}_{\theta}^{\nu} \left(
						S_{j+1}=\cdot|S_{n}=c,X_{-m}^{k}\right) \right\Vert_{1}
						\end{align*}
						where the second line follows from Lemma \ref{lem:Rev-MP} with $i=j$, $r=l=j+1$ and $n=k$, and the third follows from Lemma B.2.1 in \citesupp{stachurski2009} and the definition of $\alpha_{\theta,j+1,j}(X^{k}_{-m})$ in expression (\ref{eqn:Dobrushin}). Iterating in this fashion it follows that
						\begin{align*}
							\left \Vert \bar{P}_{\theta}^{\nu} \left( S_{j}= \cdot|S_{n}=b,X_{-m}^{k}\right) -\bar{P}_{\theta}^{\nu} \left(
							S_{j}=\cdot|S_{n}=c,X_{-m}^{k}\right) \right\Vert_{1} \leq 2 \prod_{l=j}^{n} \alpha_{\theta,l+1,l}(X^{k}_{-m}).
						\end{align*}
						
						Thus, it suffices to show that $\alpha_{\theta,l+1,l}(X^{k}_{-m}) \leq 1-\underline{q}(X_{l})$. Since
						\begin{align*}
							\alpha_{\theta,l+1,l}(X^{k}_{-m}) = 1 - \min_{a,b}\sum_{s' \in \mathbb{S}} \min \{\bar{P}_{\theta}^{\nu} \left( S_{l}=s'|S_{l+1}=a,X_{-m}^{k}\right),\bar{P}_{\theta}^{\nu} \left( S_{l}=s'|S_{l+1}=b,X_{-m}^{k}\right)\}
						\end{align*}
						(see \citesupp[p. 344]{stachurski2009}), it suffices to show that, for any $(a,b) \in \mathbb{S}^{2}$,
						\begin{equation*}
						\bar{P}_{\theta}^{\nu} \left( S_{l}=a|S_{l+1}=b,X_{-m}^{k}\right) \geq
						\underline{q}(X_{l})\varpi (X_{-m+1}^{l},a),
						\end{equation*}
						where $\varpi (X_{-m+1}^{k-1},\cdot) \in \mathcal{P}(\mathbb{S})$.
						
						To do this, first note that $\bar{P}_{\theta}^{\nu} \left( S_{l}=a|S_{l+1}=b,X_{-m}^{k}\right) = \bar{P}_{\theta}^{\nu} \left( S_{l}=a|S_{l+1}=b,X_{-m}^{l+1}\right)$, by Lemma \ref{lem:Rev-MP}, and
						\begin{align*}
						\bar{P}_{\theta}^{\nu} \left( S_{l}=a|S_{l+1}=b,X_{-m+1}^{l+1}\right) %=& \frac{\bar{P}_{\theta}^{\nu} \left( S_{l}=a,X_{-m+1}^{l+1},S_{l+1}=b\right) }{\bar{P}_{\theta}^{\nu} \left(	X_{-m+1}^{l+1},S_{l+1}=b\right) } \\
						=& \frac{p_{\theta}( X_{l}, b , X_{l+1})Q_{\theta}(X_{l},a,b)\bar{P}_{\theta}^{\nu}
							(X_{-m+1}^{l},S_{l}=a)}{\sum_{s\in \mathbb{S}}p_{\theta}(X_{l},b,X_{l+1})Q_{\theta}(X_{l},s,b)\bar{P}_{\theta}^{\nu}
							(X_{-m+1}^{l},S_{l}=s)}\\
						\geq & \underline{q}(X_{l}) \frac{\bar{P}_{\theta}^{\nu}
							(X_{-m+1}^{l},S_{l}=a)}{\sum_{s\in \mathbb{S}}\bar{P}_{\theta}^{\nu}
							(X_{-m+1}^{l},S_{l}=s)},
						\end{align*}
						where the last line follows from Assumption \ref{ass:BDD_Q}. Letting $\varpi(\cdot,X^{l}_{-m+1}) \equiv \frac{\bar{P}_{\theta}^{\nu}
							(S_{l}=\cdot \mid X_{-m+1}^{l} )}{\sum_{s\in \mathbb{S}}\bar{P}_{\theta}^{\nu}
							(S_{l}=s \mid X_{-m+1}^{l})} $, the desired result is obtained.
					\end{proof}

\begin{proof}[Proof of Lemma \ref{lem:Rev-MP}]
					   	Throughout the proof, we omit $\theta $ from the notation. Let $\mathbf{S}_{i:l}^{r}\equiv
					   	(S_{i},S_{l},S_{l+1},...,S_{r-1},S_{r})$ and note that
					   	\begin{equation*}
					   	\bar{P}_{\theta}^{\nu} \left( S_{i}|\mathbf{S}_{l}^{r} , X_{-m+1}^{n}\right) =\frac{\bar{P}_{\theta}^{\nu} \left( X^{n}_{r}\mid \mathbf{S}%
					   		_{i:l}^{r},X_{-m+1}^{r-1}\right) \bar{P}_{\theta}^{\nu} \left( \mathbf{S}%
					   		_{i:l}^{r},X_{-m+1}^{r-1}\right) }{\bar{P}_{\theta}^{\nu}\left( X^{n}_{r}\mid \mathbf{S}%
					   		_{l}^{r},X_{-m+1}^{r-1}\right) \bar{P}_{\theta}^{\nu} \left( \mathbf{S}_{l}^{r},X_{-m+1}^{r-1}%
					   		\right) },
					   	\end{equation*}%
					   	by Bayes' rule. By the Markov property, it follows that $\bar{P}_{\theta}^{\nu} \left( X^{n}_{r}\mid \mathbf{S}					   	_{i:l}^{r},X_{-m+1}^{r-1}\right) = \bar{P}_{\theta}^{\nu} \left( X^{n}_{r}\mid X_{r-1},S_{r}\right) $, so
					   	\begin{equation*}
					   	\bar{P}_{\theta}^{\nu} \left( S_{i}|\mathbf{S}_{l}^{r},X_{-m+1}^{n}\right) =\frac{\bar{P}_{\theta}^{\nu} \left(
					   		\mathbf{S}_{i:l}^{r},X_{-m+1}^{r-1}\right) }{\bar{P}_{\theta}^{\nu} \left( \mathbf{S}%
					   		_{l}^{r},X_{-m+1}^{r-1}\right) }=\frac{\bar{P}_{\theta}^{\nu} \left( S_{r}\mid \mathbf{S}%
					   		_{i:l}^{r-1},X_{-m+1}^{r-1}\right) \bar{P}_{\theta}^{\nu} \left( \mathbf{S}%
					   		_{i:l}^{r-1},X_{-m+1}^{r-1}\right) }{\bar{P}_{\theta}^{\nu} \left( S_{r}\mid \mathbf{S}%
					   		_{l}^{r-1},X_{-m+1}^{r-1}\right) \bar{P}_{\theta}^{\nu} \left( \mathbf{S}%
					   		_{l}^{r-1},X_{-m+1}^{r-1}\right) }.
					   	\end{equation*}%
					   	Observe that $\bar{P}_{\theta}^{\nu} \left( S_{r}\mid \mathbf{S}_{l}^{r-1},X_{-m+1}^{r-1}%
					   	\right) = Q_{\theta}(X_{r-1},S_{r-1},S_{r})$, and thus $\bar{P}_{\theta}^{\nu} \left( S_{i}|X_{-m+1}^{r},\mathbf{S}%
					   	_{l}^{r}\right) =\frac{\bar{P}_{\theta}^{\nu} \left( \mathbf{S}_{i:l}^{r-1},X_{-m+1}^{r-1}%
					   		\right) }{\bar{P}_{\theta}^{\nu} \left( \mathbf{S}_{l}^{r-1},X_{-m+1}^{r-1}\right) }$ and, by
					   	iterating, it follows that
					   	\begin{align*}
					   	\bar{P}_{\theta}^{\nu} \left( S_{i}|\mathbf{S}_{l}^{r},X_{-m+1}^{n}\right) = \frac{\bar{P}_{\theta}^{\nu} \left(
					   		S_{i},S_{l},X_{-m+1}^{l}\right) }{\bar{P}_{\theta}^{\nu} \left( \mathbf{S}_{l}^{l},X_{-m+1}^{l}%
					   		\right) } = \bar{P}_{\theta}^{\nu} \left( S_{i}|S_{l},X_{-m+1}^{l}\right),
					   	\end{align*}
					   	as desired.
					\end{proof}

\begin{proof}[Proof of Lemma \protect\ref{lem:Q-diff}]
												
						% Observe also that by similar steps
						% \begin{align*}
						%   P_{0} \left( S_{i} | \mathbf{X}^{j}_{-m+1} , S_{l} , S_{-m}  \right) = P_{0} \left( S_{i} | \mathbf{X}^{l}_{-m+1} , S_{l} , S_{-m}  \right)
						% \end{align*}
						% for any $j > l$.
						
						\textbf{Part (i).}  By Lemma \ref{lem:Rev-MP} with $l=r=k-1$, $n=k$, it follows that
						\begin{align*}
						\bar{P}_{\theta}^{\nu} \left( S_{i}|X_{-m+1}^{k}\right) =& \sum_{s\in \mathbb{S}}\bar{P}_{\theta}^{\nu}
						\left( S_{i}|S_{k-1}  = s,X_{-m+1}^{k}\right) \bar{P}_{\theta}^{\nu} \left(
						S_{k-1} = s|X_{-m+1}^{k}\right) \\
						=& \sum_{s \in \mathbb{S}}\bar{P}_{\theta}^{\nu} \left(
						S_{i}|S_{k-1}=s,X_{-m+1}^{k-1}\right) \bar{P}_{\theta}^{\nu} \left( S_{k-1} = s |X_{-m+1}^{k}\right) ,
						\end{align*}%
						and similarly,
						\begin{equation*}
						\bar{P}_{\theta}^{\nu} \left( S_{i}|X_{-m+1}^{k-1}\right) =\sum_{s \in \mathbb{S}}\bar{P}_{\theta}^{\nu}
						\left( S_{i}|S_{k-1}=s,X_{-m+1}^{k-1}\right) \bar{P}_{\theta}^{\nu} \left(
						S_{k-1}=s|X_{-m+1}^{k-1}\right) .
						\end{equation*}	
						Thus, by Lemma  B.2.2 in \citesupp{stachurski2009},
						\begin{align*}
							& \left \Vert \bar{P}_{\theta}^{\nu} \left( S_{j} = \cdot |X_{-m+1}^{k}\right) -\bar{P}_{\theta}^{\nu} \left(
							S_{j} = \cdot |X_{-m+1}^{k-1}\right) \right \Vert_{1} \\
							\leq & \max_{a,b}  \left \Vert \bar{P}_{\theta}^{\nu} \left( S_{j} = \cdot |S_{k-1}=a,X_{-m+1}^{k-1}\right)
								-\bar{P}_{\theta}^{\nu} \left( S_{j} = \cdot |S_{k-1}=b,X_{-m+1}^{k-1}\right) \right\Vert_{1} \\
							\leq &  L \prod_{l=j}^{k-1}(1-q(X_{l})),
						\end{align*}%						
				        where the second line follows by Lemma \ref{lem:Ker-rev} with $n=k-1$. Thus, the desired
						result follows.
						
						\bigskip
						
						\textbf{Part (ii).} The proof is analogous to that of Lemma 5 (third part)
						in \citesupp{bickel98}. By analogous calculations to those in part (i),
						\begin{align*}
						& \left \Vert \bar{P}_{\theta}^{\nu} (S_{j}=\cdot \mid X_{-m}^{k})-\bar{P}_{\theta}^{\nu} (S_{j}=\cdot\mid
						X_{-n}^{k})\right\Vert_{1} \\
						\leq &  \max_{b,b^{\prime }} \left\Vert \bar{P}_{\theta}^{\nu} (S_{j}=\cdot \mid
						S_{-m}=b,X_{-m}^{k})-\bar{P}_{\theta}^{\nu} (S_{j}=\cdot \mid S_{-m}=b^{\prime
						},X_{-n}^{k})\right\Vert_{1} \\
						=&  \max_{b,b^{\prime }} \left\Vert \bar{P}_{\theta}^{\nu} (S_{j}=\cdot \mid
						S_{-m}=b,X_{-m}^{k})-\bar{P}_{\theta}^{\nu} (S_{j}=\cdot \mid S_{-m}=b^{\prime
						},X_{-m}^{k}) \right\Vert_{1},
						\end{align*}%
						where the last line follows from the fact that, given $S_{-m}$, it is the
						same to condition on $X_{-m}^{k}$ and on $X_{-n}^{k}$. The results thus follows from following the same steps as those in the proof of Theorem \ref{thm:Q-ergo}.
					\end{proof}

\subsection{Proof of Lemma \ref{lem:suff.HAC.Rate}}

\label{SM:exa.Canon.Inference copy(1)}\label{SM:exa.Canon.Inference}

\begin{proof}[Proof of Lemma \ref{lem:suff.HAC.Rate}]
	To simplify the exposition, we present the proof for the case where $\Delta_{t,-\infty} (\theta_{\ast}) (X^{t}_{-\infty}) $ is a scalar; since the dimension of this quantity is finite, the vector case follows readily from the results below.

	Part (a) follows easily from Lemma \ref{lem:score-rep} in the Supplemental Material \ref{SM:LAR1}.
	
	For part (b), it follows from part (a) that $\Delta_{t+j,-\infty} (\theta_{\ast}) (X^{t+j}_{-\infty})  \Delta_{t,-\infty} (\theta_{\ast}) (X^{t}_{-\infty}) $ depends only on $(X_{t-\bar{L}},....,X_{t+j})$. By Lemma \ref{lem:sta-ergo}, $(X_{t})_{t=-\infty}^{\infty}$ is $\beta$-mixing. Since, for any fixed $j$, the $\sigma$-algebra generated by $(\Delta_{s+j,-\infty} (\theta_{\ast})  \Delta_{s,-\infty} (\theta_{\ast}) )_{s \leq t}$ is contained in the $\sigma$-algebra generated by $(X_{s+j})_{s \leq t}$, and the $\sigma$-algebra generated by $(\Delta_{s+j,-\infty} (\theta_{\ast})  \Delta_{s,-\infty} (\theta_{\ast}) )_{s \geq t}$ is contained in the $\sigma$-algebra generated by $(X_{s})_{s \geq t - \bar{L}}$, it follows that, for each $j$, $(\Delta_{t+j,-\infty} (\theta_{\ast})  \Delta_{t,-\infty} (\theta_{\ast}) )_{t=-\infty}^{\infty}$ is also $\beta$-mixing with mixing coefficients that decay at rate $O( \gamma^{n - 2 \bar{L} })$ as $n$ diverges through the positive integers; as $\bar{L}$ is taken to be fixed, the decay rate is $O( \gamma^{n})$. As is well known, this result implies that the corresponding $\alpha$-mixing coefficients $(\alpha _{n})_{n}$ decay at the same rate; i.e., $\alpha_{n} = O( \gamma^{n})$ for all $n$.
	
	Henceforth, let $\Omega_{t+j,t,-\infty}(\theta_{\ast}) \equiv \Delta_{t+j,-\infty} (\theta_{\ast}) (X^{t+j}_{-\infty})  \Delta_{t,-\infty} (\theta_{\ast}) (X^{t}_{-\infty})$ and $\bar{\Omega}_{t+j,t,-\infty}(\theta_{\ast}) \equiv \Omega_{t+j,t,-\infty}(\theta_{\ast}) - E_{\bar
		{P}_{\ast}^{\nu}}[\Omega_{t+j,t,-\infty}(\theta_{\ast})]$, for any $t,j$.
	Observe that, for any $j$,
	%\begin{align*}
	%	& E\left[ \left( T^{-1} \sum_{t=1}^{T} \bar{\Omega}_{t+j,t,-\infty}(\theta^{\ast})   \right) \left( T^{-1} \sum_{t=1}^{T} \bar{\Omega}_{t+j,t,-\infty}(\theta^{\ast})   \right)   \right] \\
	%	= &  T^{-1} \sum_{t=1}^{T}  E\left[ \left( \bar{\Omega}_{t+j,t,-\infty}(\theta^{\ast})   \right) \left(  \bar{\Omega}_{t+j,t,-\infty}(\theta^{\ast})  \right) \right] \\
	%	& + 2 T^{-2} \sum_{t=1}^{T} \sum_{s=0}^{t-1} E\left[ \left( \bar{\Omega}_{t+j,t,-\infty}(\theta^{\ast})  \right) \left( \bar{\Omega}_{s+j,s,-\infty}(\theta^{\ast})   \right)  \right] \\
	%	& + T^{-2} \sum_{t=1}^{T} \sum_{s=0}^{t-1} E\left[ \left( \bar{\Omega}_{s+j,s,-\infty}(\theta^{\ast})  \right) \left( \bar{\Omega}_{t+j,t,-\infty}(\theta^{\ast})   \right)   \right] \\
	%	= & T^{-1} E\left[ \left( \bar{\Omega}_{j,0,-\infty}(\theta^{\ast})   \right) \left(  \bar{\Omega}_{j,0,-\infty}(\theta^{\ast})  \right)  \right] \\
	%	& + 2 \sum_{t=0}^{T-1} (1-t/T) E\left[ \left( \bar{\Omega}_{t+j,t,-\infty}(\theta^{\ast})  \right) \left( \bar{\Omega}_{j,0,-\infty}(\theta^{\ast})   \right)   \right]
	%\end{align*}
	\begin{align*}
		E_{\bar
			{P}_{\ast}^{\nu}}\left[ \left( T^{-1} \sum_{t=1}^{T} \bar{\Omega}_{t+j,t,-\infty}(\theta_{\ast})   \right)^{2}  \right] = &  T^{-1} \sum_{t=1}^{T}  E_{\bar
			{P}_{\ast}^{\nu}}\left[ \left( \bar{\Omega}_{t+j,t,-\infty}(\theta_{\ast})   \right)^{2}  \right] \\
		& + 2 T^{-2} \sum_{t=1}^{T} \sum_{\tau=0}^{t-1} E_{\bar
			{P}_{\ast}^{\nu}}\left[ \left( \bar{\Omega}_{t+j,t,-\infty}(\theta_{\ast})  \right) \left( \bar{\Omega}_{s+j,\tau,-\infty}(\theta_{\ast})   \right)  \right] \\
		%	& + T^{-2} \sum_{t=1}^{T} \sum_{s=0}^{t-1} E\left[ \left( \bar{\Omega}_{s+j,s,-\infty}(\theta^{\ast})  \right) \left( \bar{\Omega}_{t+j,t,-\infty}(\theta^{\ast})   \right)   \right] \\
		= & T^{-1} E_{\bar
			{P}_{\ast}^{\nu}}\left[ \left( \bar{\Omega}_{j,0,-\infty}(\theta_{\ast})   \right) ^{2}  \right] \\
		&+ 2 T^{-1}  \sum_{t=0}^{T-1} (1-t/T) E_{\bar
			{P}_{\ast}^{\nu}}\left[ \left( \bar{\Omega}_{t+j,t,-\infty}(\theta_{\ast})  \right) \left( \bar{\Omega}_{j,0,-\infty}(\theta_{\ast})   \right)   \right],
	\end{align*}
	where the last equality follows by stationarity.
	By Corollary 6.17 in \citesupp{white2001}, for any $m \in \mathbb{N}$, $$|E_{\bar
		{P}_{\ast}^{\nu}}[\bar{\Omega}_{j,0,-\infty}(\theta_{\ast}) \bar{\Omega}_{j+m,m,-\infty}(\theta_{\ast}) ]| \precsim (\alpha_{m} )^{\frac{2}{2+2 \delta}} \sqrt{ E_{\bar
			{P}_{\ast}^{\nu}}[ ( \bar{\Omega}_{j,0,-\infty}(\theta_{\ast}))^{2}  ] } ( E_{\bar
		{P}_{\ast}^{\nu}}[ | \bar{\Omega}_{j+m,m,-\infty}(\theta_{\ast})|^{2+2\delta}  ] )^{\frac{1}{2+2\delta}},$$ for any $j$ (the implicit constant in the display does not depend on $j$). Thus, for any $j$,
	%\begin{align*}
	%	& E\left[ \left( T^{-1} \sum_{t=1}^{T} \bar{\Omega}_{t+j,t,-\infty}(\theta_{\ast})   \right) \left( T^{-1} \sum_{t=1}^{T} \bar{\Omega}_{t+j,t,-\infty}(\theta_{\ast})   \right)   \right] 	\precsim  T^{-1} E[ \left( \bar{\Omega}_{j,0,-\infty}(\theta_{\ast})   \right)^{2} ] \\
	%	& + 2 \sum_{t=0}^{T-1} (1-t/T) (\alpha_{t})^{2/(2+2 \delta)} \sqrt{ E[ ( \bar{\Omega}_{j,0,-\infty}(\theta_{\ast}))^{2}  ] } ( E[ | \bar{\Omega}_{j+m,m,-\infty}(\theta_{\ast})|^{2+2\delta}  ] )^{1/(2+2\delta)}.
	%\end{align*}
	\begin{align*}
		& E_{\bar
			{P}_{\ast}^{\nu}}\left[ \left( T^{-1} \sum_{t=1}^{T} \bar{\Omega}_{t+j,t,-\infty}(\theta_{\ast})   \right) ^{2}  \right] 	\precsim  T^{-1} E_{\bar
			{P}_{\ast}^{\nu}}[ \left( \bar{\Omega}_{j,0,-\infty}(\theta_{\ast})   \right)^{2} ] \\
		& + 2 T^{-1} \sum_{t=0}^{T-1} (1-t/T) (\alpha_{t})^{2/(2+2 \delta)} \sqrt{ E_{\bar
				{P}_{\ast}^{\nu}}[ ( \bar{\Omega}_{j,0,-\infty}(\theta_{\ast}))^{2}  ] } ( E_{\bar
			{P}_{\ast}^{\nu}}[ | \bar{\Omega}_{j+m,m,-\infty}(\theta_{\ast})|^{2+2\delta}  ] )^{1/(2+2\delta)}.
	\end{align*}
	By the Cauchy-Schwarz inequality and sationarity, for any $j$,
	\begin{align*}
		E_{\bar
			{P}_{\ast}^{\nu}}[ \left( \bar{\Omega}_{j,0,-\infty}(\theta^{\ast})   \right)^{2} ] \precsim E_{\bar
			{P}_{\ast}^{\nu}}[ \Delta_{0,-\infty} (\theta_{\ast}) (X^{0}_{-\infty}) ^{4}   ],
	\end{align*}
	which is bounded by assumption. In addition, by similar calculations, it follows that, for any $j$,
	\begin{align*}
		E_{\bar
			{P}_{\ast}^{\nu}}[ \left( \bar{\Omega}_{j,0,-\infty}(\theta^{\ast})   \right)^{2+2\delta} ] \precsim E_{\bar
			{P}_{\ast}^{\nu}}[ \Delta_{0,-\infty} (\theta_{\ast}) (X^{0}_{-\infty}) ^{4+4\delta}   ],
	\end{align*}
	which is bounded by assumption.
	Therefore, there exists a finite constant $C$ (which does not depend on $j$) such that
	%	\begin{align*}
	%		 E\left[ \left( T^{-1} \sum_{t=1}^{T} \bar{\Omega}_{t+j,t,-\infty}(\theta^{\ast})   \right) \left( T^{-1} \sum_{t=1}^{T} \bar{\Omega}_{t+j,t,-\infty}(\theta^{\ast})   \right)   \right] 	\leq  C ( T^{-1} + 2 \sum_{t=0}^{T-1} (1-t/T) (\alpha_{t} )^{2/(2+2 \delta)} ).
	%	\end{align*}
	\begin{align*}
		E_{\bar
			{P}_{\ast}^{\nu}}\left[ \left( T^{-1} \sum_{t=1}^{T} \bar{\Omega}_{t+j,t,-\infty}(\theta^{\ast})   \right) ^{2}  \right] 	\leq  C ( T^{-1} + 2 T^{-1} \sum_{t=0}^{T-1} (1-t/T) (\alpha_{t} )^{2/(2+2 \delta)} ).
	\end{align*}
	As $\alpha_{t} = O(\gamma^{t})$, it follows that $\sum_{t=0}^{T-1} (1-t/T) (\alpha_{t} )^{2/(2+2 \delta)} = O \left(  \sum_{t=0}^{T-1} (1-t/T) (\gamma^{2/(2+2 \delta)})^{t}   \right)$. Since $\gamma < 1$, it follows that $\sum_{t=0}^{T-1} (1-t/T) (\alpha_{t} )^{2/(2+2 \delta)} = O(1)$, which in turn implies that $E_{\bar
		{P}_{\ast}^{\nu}}\left[ \left( T^{-1} \sum_{t=1}^{T} \bar{\Omega}_{t+j,t,-\infty}(\theta^{\ast})   \right)^{2}    \right] \leq C T^{-1}$.
	Hence, by the Markov inequality, for any $a>0$,
	{\footnotesize{	
			\begin{align*}
				& \bar{P}_{\ast}^{\nu} \left( \max_{j \in \{0,...,L\}} 	|| T^{-1} \sum_{t=1}^{T} \Delta_{t+j,-\infty} (\theta_{\ast}) (X^{t+j}_{-\infty})  \Delta_{t,-\infty} (\theta_{\ast}) (X^{t}_{-\infty})^{\intercal}  -  E_{\bar
					{P}_{\ast}^{\nu}}[ \Delta_{j,-\infty} (\theta_{\ast}) (X^{j}_{-\infty})  \Delta_{0,-\infty} (\theta_{\ast}) (X^{0}_{-\infty})^{\intercal}    ] || \geq a   \right) \\
				\leq & \sum_{j=0}^{L} \bar{P}_{\ast}^{\nu} \left( || T^{-1} \sum_{t=1}^{T} \Delta_{t+j,-\infty} (\theta_{\ast}) (X^{t+j}_{-\infty})  \Delta_{t,-\infty} (\theta_{\ast}) (X^{t}_{-\infty})^{\intercal}  -  E_{\bar
					{P}_{\ast}^{\nu}}[ \Delta_{j,-\infty} (\theta_{\ast}) (X^{j}_{-\infty})  \Delta_{0,-\infty} (\theta_{\ast}) (X^{0}_{-\infty})^{\intercal}    ] || \geq a   \right) \\
				\leq & C a^{-2} LT^{-1},
			\end{align*}
	}}
	which implies the desired result.
	
\end{proof}

\section{Proof of Theorem \ref{thm:StdErrors}}

\label{sm:StdErrors}

To prove Theorem \ref{thm:StdErrors}, we use the following results, whose
proofs are relegated to the end of this section. The following lemma shows
that we can \textquotedblleft quantify\textquotedblright\ convergence in probability.

\begin{lemma}
\label{lem:Quant.op} Suppose a random sequence $(X_{T})_{T=0}^{\infty}$
converges to zero in probability. Then, there exists a deterministic positive
sequence $(r_{T})_{T=0}^{\infty}$ such that $r_{T}=o(1)$ and, for any
$\epsilon>0$, there exists $T_{\epsilon}$ such that
\[
\Pr(|X_{T}|\geq r_{T})\leq\epsilon,
\]
for all $T\geq T_{\epsilon}$. In particular, $|X_{T}|=O_{\mathrm{Pr}}(r_{T})$.
\end{lemma}

The next lemma presents some useful properties for the \textquotedblleft score
process\textquotedblright.

\begin{lemma}
\label{lem:Charac.ScoreProcess} Under the Assumptions of Theorem
\ref{thm:StdErrors}, the following are true:

\begin{itemize}
\item[1. ] $||\sup_{\theta\in B(\delta,\theta_{\ast})}||\Delta_{0,-\infty
}(\theta)||||_{L^{2}(\bar{P}_{\ast}^{\nu})}<\infty$ ($\delta>0$ is as in
Assumption \ref{ass:deriva-bdd}).

\item[2. ] $\Delta_{t,-\infty}$ and $\Delta_{t,-\infty}\Delta_{0,-\infty
}^{\intercal}$ are continuous in $L^{1}(\bar{P}_{\ast}^{\nu})$-norm, i.e., for
any $\epsilon>0$, there exists $\delta>0$ such that
\[
\max_{t}\left\Vert \sup_{||\theta-\theta_{0}||<\delta}||\Delta_{t,-\infty
}(\theta)\Delta_{t,-\infty}(\theta)^{\intercal}-\Delta_{t,-\infty}(\theta
_{0})\Delta_{t,-\infty}(\theta_{0})^{\intercal}||\right\Vert _{L^{1}(\bar
{P}_{\ast}^{\nu})}\leq\epsilon
\]
and
\[
\ddot{\varpi}(\delta)\equiv\max_{t}\left\Vert \sup_{||\theta-\theta
_{0}||<\delta}||\Delta_{t,-\infty}(\theta)\Delta_{0,-\infty}(\theta
)^{\intercal}-\Delta_{t,-\infty}(\theta_{0})\Delta_{0,-\infty}(\theta
_{0})^{\intercal}||\right\Vert _{L^{1}(\bar{P}_{\ast}^{\nu})}\leq\epsilon.
\]

\item[3. ] There exists a constant $C<\infty$ such that, for any $t$ and $M$,
\[
\left\Vert \sup_{\theta\in B(\delta,\theta_{\ast})}\Delta_{t,-\infty}%
(\theta)-\Delta_{t,t-M}(\theta)\right\Vert _{L^{2}(\bar{P}_{\ast}^{\nu})}\leq
CM^{1-1/p}.
\]
Moreover, by Assumption \ref{ass:q-sum-p}, $p\in(0,2/3)$ and thus the RHS
vanishes as $M$ diverges.
\end{itemize}
\end{lemma}

\begin{proof}[Proof of Theorem \ref{thm:StdErrors}]
	The proof has several parts and steps.
	
	\bigskip
	
	\textsc{Part (a).} We show that
	\begin{align*}
		||T^{-1}\sum_{t=1}^{T}\nabla_{\theta}^{2}\log p_{t}^{\nu}( X_{t} | X^{t-1}_{0} ,\hat{\theta}_{\nu,T} )- E_{ \bar{P}_{\ast}^{\nu} } [  \xi_{1}(\theta_{\ast}) ] || = o_{\bar{P}_{\ast}^{\nu}}(1).
	\end{align*}
	We do this by using the triangle inequality and showing that the following expressions hold:
	\begin{align*}
		\lim_{T\rightarrow\infty} \left \Vert  T^{-1}\sum_{t=1}^{T}\{\nabla_{\theta}^{2}\log p_{t}^{\nu}( \cdot | \cdot , \hat{\theta}_{\nu,T}  )-\xi_{t}( \hat{\theta}_{\nu,T} )\} \right\Vert _{L^{1}(\bar{P}_{\ast}^{\nu})}   = 0
	\end{align*}
	(which implies convergence in probability),
	\begin{align*}
		|| T^{-1}\sum_{t=1}^{T} \{ \xi_{t}( \hat{\theta}_{\nu,T} ) -  \xi_{t}( \theta_{\ast} ) \}  ||    = o_{\bar{P}_{\ast}^{\nu}}(1),
	\end{align*}
	and
	\begin{align*}
		||T^{-1}\sum_{t=1}^{T} \xi_{t}(\theta_{\ast})  - E_{ \bar{P}_{\ast}^{\nu} } [  \xi_{1}(\theta_{\ast}) ]  || = o_{\bar{P}_{\ast}^{\nu}}(1).
	\end{align*}
	
	The first expression holds true because, by Theorem \ref{thm:consistent}, for any $\delta' \leq \delta$, $\hat{\theta}_{\nu,T}  \in B(\delta',\theta_{\ast}) $ w.p.a.1 and, hence, by Lemma \ref{lem:Hess-approx}, the desired result follows.
	
	Regarding the second expression,  again by Theorem \ref{thm:consistent}, $\hat{\theta}_{\nu,T}  \in B(\delta',\theta_{\ast}) $ w.p.a.1. Thus, it follows that, for any $\epsilon>0$, there exists $T(\epsilon)$ such that, for any $ t \geq T(\epsilon)$,
	\begin{align*}
		\bar{P}_{\ast}^{\nu}  \left(  	 || T^{-1}\sum_{t=1}^{T} \{  \xi_{t}( \hat{\theta}_{\nu,T} ) -  \xi_{t}( \theta_{\ast} ) \}  ||    \geq \epsilon \right) \leq &  \bar{P}_{\ast}^{\nu}  \left(  	 || T^{-1}\sum_{t=1}^{T}  \{ \xi_{t}( \hat{\theta}_{\nu,T} ) -  \xi_{t}( \theta_{\ast} ) \} ||    \geq \epsilon, \hat{\theta}_{\nu,T}  \in B(\delta',\theta_{\ast})   \right)  + \epsilon \\
		\leq &  \epsilon^{-1} E_{ \bar{P}_{\ast}^{\nu} } \left[  \sup_{\theta \in B(\delta',\theta_{\ast}) }	 ||  \xi_{1}( \theta ) -  \xi_{1}( \theta_{\ast} )  ||      \right] + 0.5 \epsilon,
	\end{align*}
	where the second line follows from the Markov inequality and stationarity. By Lemma \ref{lem:Hess-approx}, $\xi_{1}$ is continuous -- and thus uniformly continuous over compact sets. Since $\delta'>0$ can be chosen to be any number less than $\delta$ ($\delta$ as in Assumption \ref{ass:deriva-bdd}), we can choose it so that the first term in the RHS is less than  $0.5 \epsilon$. Hence, the desired follows.
	
	Finally, ergodicity of $X_{-\infty}^{\infty}$ (Lemma~\ref{lem:sta-ergo}) implies ergodicity of $(\xi _{t}(\theta_{\ast}))_{t=-\infty}^{\infty} $; therefore, by Lemma~\ref{lem:Hess-approx} and Birkhoff's ergodic theorem, $  ||T^{-1}\sum_{t=1}^{T} \xi_{t}(\theta_{\ast})  - E_{ \bar{P}_{\ast}^{\nu} } [  \xi_{1}(\theta_{\ast}) ]  || = o_{\bar{P}_{\ast}^{\nu}}(1)$. Hence,
	\begin{align*}
		||T^{-1}\sum_{t=1}^{T} \nabla_{\theta}^{2}\log p_{t}^{\nu}( X_{t} | X^{t-1}_{0}  ,\hat{\theta}_{\nu,T} )-  E_{ \bar{P}_{\ast}^{\nu} } [  \xi_{1}(\theta_{\ast}) ]  || = o_{\bar{P}_{\ast}^{\nu}}(1).
	\end{align*}	
	
	With this result and the Fisher information equality (established in the proof of Corollary \ref{cor:anormal}), the result of part (a) of the theorem follows.
	
	\bigskip
	
	\textsc{Part (b). Step 1} To prove part (b), it suffices to show that
	\begin{align*}
		|| E_{\bar{P}_{\ast}^{\nu}}[\xi_{1}(\theta_{\ast})]- H_{T}(\hat{\theta}_{\nu,T}  ) || = o_{\bar{P}_{\ast}^{\nu}}(1)
	\end{align*}
	and
	\begin{align*}
		||  \Sigma_{T}(\theta_{\ast}) - J_{T}(\hat{\theta}_{\nu,T}  )  || = o_{\bar{P}_{\ast}^{\nu}}(1).
	\end{align*}
	
	The first expression was established in Part (A). Regarding the second expression, we introduce some notation. For any $t,l,M \in \mathbb{N}$ and any $\theta \in \Theta$, let
	\begin{align*}
		\Omega_{t,l,M}( \theta ) \equiv \Delta_{t,M}(\theta)(X^{t}_{M}) \Delta_{l,M}(\theta) (X^{l}_{M})^{\intercal},
	\end{align*}
	where it is left implicit that this quantity depends on $X^{\max\{t,l\}}_{-M}$. Also, for any $\tau \in \{1,$\ldots $,L_{T}\}$,
	\begin{align*}
		\theta \mapsto 	\hat{\gamma}_{T,\tau,0}(\theta)  \equiv (T-\tau)^{-1} \sum_{t=1}^{T-\tau} \Omega_{t+\tau,t,0}(\theta).
	\end{align*}
	Recall that $ \Delta_{t,0} (\theta) (X^{t}_{0}) = \nabla_{\theta} p^{\nu}_{t}(X_{t} \mid X^{t-1}_{0},\theta)$, so  $	\hat{\gamma}_{T,\tau,0}(\theta) $ is the sample covariance of $\Omega_{t+\tau,t,0}(\theta)$.

	Given this notation, observe that
		\begin{align*}
		\Sigma_{T}(\theta_{\ast})  = & T^{-1} \sum_{t=1}^{T}  E_{\bar{P}_{\ast}^{\nu}} [  \Omega_{t,t,-\infty}(\theta_{\ast}) ] + T^{-1} \sum_{t=1}^{T} \sum_{l=0}^{t-1} \{ E_{\bar{P}_{\ast}^{\nu}} [ \Omega_{t,l,-\infty}(\theta_{\ast})  ] + E_{\bar{P}_{\ast}^{\nu}} [ \Omega_{t,l,-\infty}(\theta_{\ast})^{\intercal}]\} .
	\end{align*}
	The aim is to show that each of the summands above is well-approximated by its counterpart in $J_{T}$. For the first summand, we show in Step 2 below that
\begin{align*}
	||	T^{-1} \sum_{t=1}^{T} \{  \Omega_{t,t,0} (\hat{\theta}_{\nu,T})   -  E_{\bar{P}_{\ast}^{\nu}} [\Omega_{t,t,-\infty}(\theta_{\ast})]  \}  || = o_{  \bar{P}_{\ast}^{\nu}  }(1)  .
\end{align*}
	Regarding the second summand, we observe that, for any $t \geq l$,
\begin{align*}
	E_{\bar{P}_{\ast}^{\nu} } [ \Omega_{t,l,-\infty}(\theta_{\ast}) ] = E_{\bar{P}_{\ast}^{\nu} } [ \Omega_{t-l,0,-\infty}(\theta_{\ast}) ] \equiv \gamma_{t-l}(\theta_{\ast})
\end{align*}
	(the first equality, which follows from stationarity, can be established by analogous arguments to those presented at the beginning of Step 2). Hence,
\begin{align*}
	T^{-1}  \sum_{t=1}^{T} \sum_{l=0}^{t-1}  	E_{\bar{P}_{\ast}^{\nu} } [ \Omega_{t,l,-\infty}(\theta_{\ast}) ] =  & T^{-1} ( \gamma_{1}(\theta_{\ast}) + ( \gamma_{2}(\theta_{\ast}) + \gamma_{1} (\theta_{\ast})  ) + \cdots + (\gamma_{T} (\theta_{\ast}) + \cdots + \gamma_{1} (\theta_{\ast}) )) \\
	= &  \sum_{j=0}^{T-1} (1-j/T) \gamma_{j+1}(\theta_{\ast}).
\end{align*}	
Thus, it suffices to show that
	\begin{align*}
		||	 \sum_{j=0}^{T-1} (1-j/T) \gamma_{j+1}(\theta_{\ast}) - \sum_{j=0}^{L_{T}-1} \omega(j,L)  \hat{\gamma}_{T,j+1,0}(\hat{\theta}_{T,\nu})   || = o_{  \bar{P}_{\ast}^{\nu}  }(1),
	\end{align*}
the proof of which is in Step 3 below.

	\bigskip

	\textsc{Step 2} We now show that
	\begin{align*}
		||	T^{-1} \sum_{t=1}^{T} \{  \Delta_{t,0}(\hat{\theta}_{\nu,T})(X^{t}_{0})  \Delta_{t,0}(\hat{\theta}_{\nu,T})(X^{t})^{\intercal}   -  E_{\bar{P}_{\ast}^{\nu}} [\Delta_{t,-\infty}(\theta_{\ast})(X^{t}_{-\infty}) \Delta_{t,-\infty}(\theta_{\ast}) (X^{t}_{-\infty})^{\intercal}]   \}  || = o_{  \bar{P}_{\ast}^{\nu}  }(1),
	\end{align*}
	where, by the definition of $\Delta_{t,M}$, $\Delta_{t,0}(\theta)(X^{t}) = \nabla_{\theta} \log p_{t}(X_{t} \mid X^{t-1}_{0} , \theta)$.
	
	First, observe that $$E_{\bar{P}_{\ast}^{\nu}} [\Delta_{t,-\infty}(\theta_{\ast})(X^{t}_{-\infty}) \Delta_{t,-\infty}(\theta_{\ast}) (X^{t}_{-\infty})^{\intercal}] = E_{\bar{P}_{\ast}^{\nu}} [\Delta_{0,-\infty}(\theta_{\ast})(X^{0}_{-\infty}) \Delta_{0,-\infty}(\theta_{\ast}) (X^{0}_{-\infty})^{\intercal}].$$ This follows from stationarity (see Lemma \ref{lem:sta-ergo}) and the fact that $\Delta_{t,-\infty}(\theta_{\ast}) $ can be approximated (uniformly in $t$) by $\Delta_{t,t-M}(\theta_{\ast}) $ (see Lemma \ref{lem:score_approx}). Hence, it suffices to show that \begin{align*}
		||	T^{-1} \sum_{t=1}^{T}  \Delta_{t,0}(\hat{\theta}_{\nu,T})(X^{t})  \Delta_{t,0}(\hat{\theta}_{\nu,T})(X^{t})^{\intercal}   -   E_{\bar{P}_{\ast}^{\nu}} [\Delta_{0,-\infty}(\theta_{\ast})(X^{0}_{-\infty}) \Delta_{0,-\infty}(\theta_{\ast}) (X^{0}_{-\infty})^{\intercal}]    || = o_{  \bar{P}_{\ast}^{\nu}  }(1).
	\end{align*}
	
	By Lemma \ref{lem:sta-ergo},  ergodicity of $\Delta_{t,t-M}(\theta_{\ast}) $ for any $M$ follows. This, Birkhoff's ergodic theorem, and Lemma \ref{lem:score_approx} imply that
	\begin{align*}
		||	T^{-1} \sum_{t=1}^{T} \bar{\Delta}_{\infty} (\theta_{\ast})(X^{t}_{-\infty})  - E_{\bar{P}_{\ast}^{\nu}}[  \bar{\Delta}_{\infty} (\theta_{\ast})(X^{0}_{-\infty})   ] || = o_{  \bar{P}_{\ast}^{\nu}  }(1),
	\end{align*}
	where, for any $M \in \mathbb{Z} \cup \{ \infty \}$, $\bar{\Delta}_{M}(\theta)(X^{t}_{t-M})  \equiv \Delta_{0,-M}(\theta)(X^{t}_{t-M})   \Delta_{0,-M}(\theta)(X^{t}_{t-M})^{\intercal}   $.

	Hence, in order to obtain the desired result it suffices to show that
	\begin{align*}
		||	T^{-1} \sum_{t=1}^{T} \{ \bar{\Delta}_{t}(\hat{\theta}_{\nu,T} )(X^{t}_{0})     -  \bar{\Delta}_{\infty}(\theta_{\ast})(X^{t}_{-\infty}) \}   || = o_{  \bar{P}_{\ast}^{\nu}  }(1).
	\end{align*}
	In order to do so, by the triangle inequality, it is sufficient to show that
	\begin{align}\label{eqn:HAC-step2-1}
		||	T^{-1} \sum_{t=1}^{T} \{ \bar{\Delta}_{\infty} (\theta_{\ast})(X^{t}_{-\infty}) - \bar{\Delta}_{\infty} (\hat{\theta}_{\nu,T} )(X^{t}_{-\infty})  \} || = o_{  \bar{P}_{\ast}^{\nu}  }(1)
	\end{align}
	and
	\begin{align}\label{eqn:HAC-step2-2}
		||	T^{-1} \sum_{t=1}^{T} \{ \bar{\Delta}_{t}(\hat{\theta}_{\nu,T} )(X^{t}_{0})     -  \bar{\Delta}_{\infty} (\hat{\theta}_{\nu,T} )(X^{t}_{-\infty})  \}   || = o_{  \bar{P}_{\ast}^{\nu}  }(1).
	\end{align}
	
	Expression (\ref{eqn:HAC-step2-1}) holds by Lemma \ref{lem:Charac.ScoreProcess}, the fact that, for any $\delta>0$, $\hat{\theta}_{T,\nu} \in B(\delta,\theta_{\ast})$ w.p.a.1-$\bar{P}^{\nu}_{\ast}$ (by Theorem \ref{thm:consistent}) and the Markov inequality. Regarding expression (\ref{eqn:HAC-step2-2}), by the Markov inequality and the fact that $\hat{\theta}_{T,\nu} \in B(\delta,\theta_{\ast})$ w.p.a.1-$\bar{P}^{\nu}_{\ast}$ (by Theorem \ref{thm:consistent} ), it is sufficient to show that
	\begin{align*}
		T^{-1} \sum_{t=1}^{T} 	E_{ \bar{P}_{\ast}^{\nu} } [ \sup_{  \theta \in B(\delta,\theta_{\ast})  } ||  \bar{\Delta}_{t}(\theta )(X^{t}_{0})     -  \bar{\Delta}_{\infty} (\theta_{\ast} )(X^{t}_{-\infty})   ||  ] = o(1).
	\end{align*}
	In turn, the LHS is bounded by
	\begin{align*}
		& T^{-1} \sum_{t=1}^{T}  || 	\sup_{  \theta \in B(\delta,\theta_{\ast})  } ||  \Delta_{t,0}(\theta )    -  \Delta_{t,-\infty} (\theta_{\ast} )||  ||_{L^{2}(\bar{P}_{\ast}^{\nu} )} \\
		& \times (  ||  \sup_{\theta \in B(\delta,\theta_{\ast})  } ||  \Delta_{t,0}(\theta ) ||  ||_{L^{2}(\bar{P}_{\ast}^{\nu}  )}  +  ||  \sup_{\theta \in B(\delta,\theta_{\ast})  } ||  \Delta_{t,-\infty}(\theta )|| ||_{L^{2}(\bar{P}_{\ast}^{\nu}  )}    ).
	\end{align*}
	By Lemma \ref{lem:Charac.ScoreProcess}(3), the first term in the RHS is bounded (up to constants) by $T^{-1} \sum_{t=1}^{T} t^{1-1/p}$. The second term in the RHS is bounded by Lemma \ref{lem:Charac.ScoreProcess}(1). Thus, under Assumption \ref{ass:q-sum-p}, the whole expression converges to zero and the desired result follows.

	\bigskip
	
	\textsc{Step 3.} We next show that, for any $L \equiv L_{T}$ 	such that $\lim_{T \rightarrow \infty} L_{T} = \infty$ and  $L_{T} (  \ddot{\varpi}(  T^{-1/2} \log \log T ) \log \log T +r_{T} + T^{-1/2} ) = o(1) $,
	\begin{align*}
		||	 \sum_{j=0}^{T-1} (1-j/T) \gamma_{j+1}(\theta_{\ast}) - \sum_{j=0}^{L_{T}-1} \omega(j,L)  \hat{\gamma}_{T,j+1,0}(\hat{\theta}_{T,\nu})   || = o_{  \bar{P}_{\ast}^{\nu}  }(1),
	\end{align*}
	where, for any $\tau \in \{1,\ldots,L_{T}\}$ and any $M \leq 1$,
	\begin{align*}
		\theta \mapsto 	\hat{\gamma}_{T,\tau,M}(\theta)  \equiv T^{-1} \sum_{t=1}^{T-\tau}  \Delta_{t+\tau,M} (\theta) (X^{t+\tau}_{M})  \Delta_{t,M} (\theta) (X^{t}_{M})^{\intercal}
	\end{align*}
	(recall that $ \Delta_{t,0} (\theta) (X^{t}_{0}) = \nabla_{\theta} p^{\nu}_{t}(X_{t} \mid X^{t-1}_{0},\theta)$).
	%and, for any $t \in \{0,...,T,...\}$ and $M \in \mathbb{N} \cup \{ \infty\}$,   $\tilde{\Delta}_{t,0,-M}(\theta) \equiv  \Delta_{t,-M}(\theta) \Delta_{0,t-M}(\theta) ' $.
	
	Putting $\hat{\gamma}_{T,\tau} \equiv \hat{\gamma}_{T,\tau,-\infty} $, we have, by the triangle inequality,
		\begin{align*}
		& ||	 \sum_{j=0}^{T-1} (1-j/T) \gamma_{j+1}(\theta_{\ast}) - \sum_{j=0}^{L_{T}-1} \omega(j,L)  \hat{\gamma}_{T,j+1,0}(\hat{\theta}_{T,\nu})   || \\
		\leq & ||	 \sum_{j=0}^{T-1} (1-j/T) \gamma_{j+1}(\theta_{\ast}) - \sum_{j=0}^{L_{T}-1} \omega(j,L) \gamma_{j+1}(\theta_{\ast})   ||\\
		& + ||	\sum_{j=0}^{L_{T}-1} \omega(j,L)  \{ \gamma_{j+1}(\theta_{\ast})  -   \hat{\gamma}_{T,j+1}(\theta_{\ast}) \}  || \\
		& + ||	\sum_{j=0}^{L_{T}-1} \omega(j,L)  \{  \hat{\gamma}_{T,j+1}(\theta_{\ast})  -  \hat{\gamma}_{T,j+1}(\hat{\theta}_{T,\nu})  \}  || \\
		& + ||	\sum_{j=0}^{L_{T}-1} \omega(j,L)  \{ \hat{\gamma}_{T,j+1,0}(\hat{\theta}_{T,\nu}) -   \hat{\gamma}_{T,j+1}(\hat{\theta}_{T,\nu}) \}  ||.	
	\end{align*}
We now bound each term in the RHS individually.

	By assumption, for any $l \geq 0$, $||\gamma_{l}(\theta_{\ast}) || \leq  \upsilon(l) $ and thus, for any $L$,
	\begin{align*}
		||  \sum_{j=0}^{T-1} (1-j/T) \gamma_{j+1}(\theta_{\ast}) -  \sum_{j=0}^{L-1} \omega(j,L) \gamma_{j+1}(\theta_{\ast}) || \leq &	||  \sum_{j=L}^{T-1}  \gamma_{j+1}(\theta_{\ast}) || \\
		& + 	||  \sum_{j=0}^{L-1} \{ (1-j/T) -  \omega(j,L) \} \gamma_{j+1}(\theta_{\ast}) || \\
		\leq & 	  \sum_{j=L}^{\infty}  \upsilon(j) \\
		& + 	||  \sum_{j=0}^{L-1} \{ (1-j/T) -  \omega(j,L) \} \gamma_{j+1}(\theta_{\ast}) || .
	\end{align*}
	Since, $\upsilon$ is integrable, the first term in the RHS converges to zero as $L$ diverges. Furthermore, since $\omega(\cdot,\cdot)$ is bounded, $||\gamma_{j+1}(\theta_{\ast})|| \leq \upsilon(j+1)$, which is integrable, and $(1-j/T) -  \omega(j,L)$ converges to zero pointwise in $j$ as $T$ (and thus $L=L_{T}$) diverges so, by the dominated convergence theorem, the second term also converges to zero as $T$ (and thus $L=L_{T}$) diverges. Therefore, for any $\epsilon>0$, there exists $T_{\epsilon}$ such that, for all $T \geq T_{\epsilon}$,
	\begin{align*}
		||  \sum_{j=0}^{T-1} (1-j/T) \gamma_{j+1}(\theta_{\ast}) -  \sum_{j=0}^{L_{T}-1} \omega(j,L_{T}) \gamma_{j+1}(\theta_{\ast}) || \leq \epsilon.
	\end{align*}
	
	We now consider
	\begin{align*}
		\bar{P}_{\ast}^{\nu} \left(  ||  \sum_{j=0}^{L-1} \omega(j,L) \{  T^{-1} \sum_{t=1}^{T-j}  \Delta_{t+j,-\infty} (\theta_{\ast}) (X^{t+j}_{-\infty})  \Delta_{t,-\infty} (\theta_{\ast}) (X^{t}_{-\infty})^{\intercal} - \gamma_{j}(\theta_{\ast})   \} ||    \geq  \delta   \right)
	\end{align*}
	for any $\delta>0$. Since $ \sum_{j=0}^{L-1} \omega(j,L) \leq L$, it follows that this expression is bounded above by
	\begin{align*}
		\bar{P}_{\ast}^{\nu} \left(  \max_{j \in \{0,...,L_{T}\}}  || T^{-1} \sum_{t=1}^{T-j}  \Delta_{t+j,-\infty} (\theta_{\ast}) (X^{t+j}_{-\infty})  \Delta_{t,-\infty} (\theta_{\ast}) (X^{t}_{-\infty})^{\intercal} - \gamma_{j}(\theta_{\ast}) ||   \geq  \delta /L   \right).
	\end{align*}
	By similar arguments to those presented in Step 2 and Birkhoff's ergodic theorem, it follows that, for each $L$,
	\begin{align*}
		\max_{j \in \{0,...,L\}}  || T^{-1} \sum_{t=1}^{T-j}  \Delta_{t+j,-\infty} (\theta_{\ast}) (X^{t+j}_{-\infty})  \Delta_{t,-\infty} (\theta_{\ast}) (X^{t}_{-\infty})^{\intercal} - \gamma_{j}(\theta_{\ast})   || = o_{\bar{P}_{\ast}^{\nu}}(1).
	\end{align*}
	By Lemma \ref{lem:Quant.op}, for each $L$, there exists a positive sequence $(r_{T})_{T}$ such that $r_{T} = o(1)$ and  $\bar{P}_{\ast}^{\nu} \left( \max_{j \in \{0,...,L\}} 	|| T^{-1} \sum_{t=1}^{T} \Delta_{t+j,-\infty} (\theta_{\ast}) (X^{t+j}_{-\infty})  \Delta_{t,-\infty} (\theta_{\ast}) (X^{t}_{-\infty})^{\intercal}  -  \gamma_{j}(\theta_{\ast}) ||  \geq  r_{T}      \right) = o(1)$.	Thus, by setting $\delta = 2 r_{T} L$ , it follows that, for any $\epsilon>0$, there exists $T_{\epsilon}$ such that, for all $T \geq T_{\epsilon}$,
	\begin{align*}
		\bar{P}_{\ast}^{\nu} \left(  ||  \sum_{j=0}^{L-1} \omega(j,L) \{  T^{-1} \sum_{t=1}^{T-j} \Delta_{t+j,-\infty} (\theta_{\ast}) (X^{t+j}_{-\infty})  \Delta_{t,-\infty} (\theta_{\ast}) (X^{t}_{-\infty})^{\intercal} - \gamma_{j}(\theta_{\ast})   \} ||    \geq 2 r_{T} L  \right) \leq \epsilon.
	\end{align*}
	
	%XXXXXXX
	%
	%
	%By similar arguments to those presented in Step 2, it follows that
	%\begin{align*}
	%\max_{s \in \{ 1,...,L  \}} 	|| T^{-1} \sum_{t=1}^{T}  \tilde{\Delta}_{s,-\infty} (\theta_{\ast})(X^{s}_{-\infty})  - \gamma_{s}(\theta_{\ast})  || = o_{\bar{P}_{\ast}^{\nu}}(1).
	%\end{align*}
	%By Lemma XXX, for any $L$, there exists a positive sequence $(r_{T})_{T}$ such that $r_{T} = o(1)$ and  $\lim_{T \rightarrow \infty} \bar{P}_{\ast}^{\nu} \left(  \max_{s \in \{ 1,...,L  \}} 	|| T^{-1} \sum_{t=1}^{T}  \tilde{\Delta}_{s,-\infty} (\theta_{\ast})(X^{s}_{-\infty})  -  \gamma (s) ||  \geq  r_{T}      \right) = 0$. Hence, by  letting $\hat{\gamma}_{T,s}(\theta)  \equiv 	T^{-1} \sum_{t=1}^{T}  \tilde{\Delta}_{s,-\infty} (\theta)(X^{s}_{-\infty})  $, it follows that
	%\begin{align*}
	%	\bar{P}_{\ast}^{\nu} \left(  ||  \sum_{j=0}^{L-1} \omega(j,L) \{  \hat{\gamma}_{T,j+1}( \theta_{\ast}  ) - 	   \gamma_{j+1}(\theta_{\ast})   \} ||    \geq  2 r_{T} L   \right) \leq  	\bar{P}_{\ast}^{\nu} \left(   r_{T} L  \geq 2 r_{T} L    \right) = 0
	%\end{align*}
	%where the last inequality follows because $ \sum_{j=0}^{L-1} \omega(j,L) \leq L$.

	By Theorem \ref{thm:anormal},  $\hat{\theta}_{\nu,T} \in B(T^{-1/2} \log \log T,\theta_{\ast})$ w.p.a.1-$\bar{P}^{\nu}_{\ast}$. Hence, for any $\epsilon>0$ there exists $T_{\epsilon}$ such that for, all $T \geq T_{\epsilon}$,
	\begin{align*}
		& \bar{P}_{\ast}^{\nu} \left(  ||  \sum_{j=0}^{L-1} \omega(j,L) \{  \hat{\gamma}_{T,j+1}( \theta_{\ast}  ) - \hat{\gamma}_{T,j+1}( \hat{\theta}_{\nu,T}  )   \} ||    \geq L \ddot{\varpi}(  T^{-1/2} \log \log T ) \log \log T   \right) \\
		\leq & \bar{P}_{\ast}^{\nu} \left( \sup_{\theta \in B(T^{-1/2} \log \log T,\theta_{\ast})}   ||  \sum_{j=0}^{L-1} \omega(j,L) \{  \hat{\gamma}_{T,j+1}( \theta_{\ast}  ) - 	   \hat{\gamma}_{T,j+1}( \theta )   \} ||    \geq L \ddot{\varpi}(  T^{-1/2} \log \log T ) \log \log T    \right)  \\
		&+  \bar{P}_{\ast}^{\nu}(\hat{\theta}_{\nu,T} \notin B(T^{-1/2} \log \log T,\theta_{\ast})) \\
		\leq & \bar{P}_{\ast}^{\nu} \left( \sup_{\theta \in B(T^{-1/2} \log \log T,\theta_{\ast})}   \sum_{j=0}^{L-1} \omega(j,L)  ||  \hat{\gamma}_{T,j+1}( \theta_{\ast}  ) - \hat{\gamma}_{T,j+1}( \theta )    ||    \geq L \ddot{\varpi}(  T^{-1/2} \log \log T ) \log \log T    \right)  + \epsilon.
	\end{align*}
	Moreover, by the Markov inequality,
	\begin{align*}
		& \bar{P}_{\ast}^{\nu} \left( \sup_{\theta \in B(T^{-1/2} \log \log T,\theta_{\ast})}   \sum_{j=0}^{L-1} \omega(j,L)  ||  \hat{\gamma}_{T,j+1}( \theta_{\ast}  ) - \hat{\gamma}_{T,j+1}( \theta )    ||    \geq L \ddot{\varpi}(  T^{-1/2} \log \log T ) \log \log T    \right) \\
		&	\leq  \frac{  1 }{L \ddot{\varpi}(  T^{-1/2} \log \log T ) \log \log T   } E_{\bar{P}_{\ast}^{\nu} } \left[ \sup_{\theta \in B(T^{-1/2} \log \log T,\theta_{\ast})}  \sum_{j=0}^{L-1} \omega(j,L)   ||  \hat{\gamma}_{T,j}( \theta_{\ast}  ) -   \hat{\gamma}_{T,j}( \theta )    ||     \right] \\
		&	\leq  \frac{ 1 }{ L \ddot{\varpi}(  T^{-1/2} \log \log T ) \log \log T  }  \sum_{j=0}^{L-1} \omega(j,L)  E_{\bar{P}_{\ast}^{\nu} } \left[   \sup_{\theta \in B(T^{-1/2} \log \log T,\theta_{\ast})} ||  \tilde{\Delta}_{j,-\infty} ( \theta_{\ast}  ) -   \tilde{\Delta}_{j,-\infty} ( \theta )    ||     \right] \\
		&	\leq  \frac{  \sum_{j=0}^{L-1} \omega(j,L)  \ddot{\varpi}(  T^{-1/2} \log \log T )    }{ L \ddot{\varpi}(  T^{-1/2} \log \log T ) \log \log T   },
	\end{align*}
	where the last line follows from Lemma \ref{lem:Charac.ScoreProcess} and $\tilde{\Delta}_{t,-M}(\theta) \equiv  \Delta_{t,-M}(\theta) \Delta_{0,-M}(\theta)^{\intercal}$ for any $t \in \{0,\ldots,T,\ldots\}$ and $M \in \mathbb{N} \cup \{ \infty\}$. Since $ \sum_{j=0}^{L-1} \omega(j,L)/L \leq 1$, the last expression is less than $\epsilon$ for sufficiently large $T$. Thus,
	\begin{align*}
		& \bar{P}_{\ast}^{\nu} \left(  ||  \sum_{j=0}^{L-1} \omega(j,L) \{  \hat{\gamma}_{T,j+1}( \theta_{\ast}  ) - \hat{\gamma}_{T,j+1}( \hat{\theta}_{\nu,T}  )   \} ||    \geq L \ddot{\varpi}(  T^{-1/2} \log \log T ) \log \log T   \right) \leq \epsilon.
	\end{align*}
	
	Finally, since, for any $\tau \in \{1,\ldots,L\}$, $\theta \mapsto \hat{\gamma}_{T,\tau,0}(  \theta   )  \equiv T^{-1} \sum_{t=1}^{T-\tau} \Delta_{t+\tau,0}(\theta)(X^{t+\tau}_{0})   \Delta_{t,0}(\theta)(X^{t}_{0})^{\intercal}  $, by Lemma \ref{lem:Charac.ScoreProcess}(3) (with $M=t$),
	\begin{align*}
		|| \sum_{j=0}^{L-1} \omega(j,L) \{  \hat{\gamma}_{T,j+1}(  \hat{\theta}_{\nu,T}   ) - \hat{\gamma}_{T,j+1,0}( \hat{\theta}_{\nu,T}  )   \} ||_{L^{2}(\bar{P}_{\ast}^{\nu})} \precsim \sum_{j=0}^{L-1} \omega(j,L) T^{-1} \sum_{t=1}^{T-j} t^{1-1/p}.
	\end{align*}
	By the proof of Lemma \ref{lem:score_approx},  $T^{-1/2} \sum_{t=1}^{T} t^{1-1/p}$ vanishes; thus, the RHS is of order $o(L T^{-1/2}) $.
	
	Therefore, we have shown that, for any $\epsilon>0$, there exists $T_{\epsilon}$ such that, for all $T \geq T_{\epsilon}$,
	{\footnotesize{
			\begin{align*}
				\bar{P}_{\ast}^{\nu} \left(   	||	 \sum_{j=0}^{T-1} (1-j/T) \gamma_{j+1}(\theta_{\ast}) - \sum_{j=0}^{L_{T}-1} \omega(j,L)  \hat{\gamma}_{T,j+1,0}(\hat{\theta}_{T,\nu})   ||  \geq \epsilon + L_{T} (  \ddot{\varpi}(  T^{-1/2} \log \log T ) \log \log T + 2r_{T} + T^{-1/2} )   \right) \leq \epsilon,
	\end{align*}}}
	where the $\epsilon$ inside the probability arises from bounding $	||	 \sum_{j=0}^{T-1} (1-j/T) \gamma_{j+1}(\theta_{\ast})  - 	 \sum_{j=0}^{L_{T}-1} \omega(j,L_{T}) \gamma_{j+1}(\theta_{\ast}) ||$ and  requires $L_{T}$ to diverge. Therefore, by taking $L \equiv L_{T}$ such that $\lim_{T \rightarrow \infty} L_{T} = \infty$ and  $L_{T} (  \ddot{\varpi}(  T^{-1/2} \log \log T ) \log \log T +r_{T} + T^{-1/2} ) = o(1) $, we establish the desired result.

\end{proof}

\subsection{Proofs of Supplementary Lemmas}

\begin{proof}[Proof of Lemma \ref{lem:Quant.op}]
	Since $(X_{T})_{T=0}^{\infty}$ converges to zero in probability, for any $\epsilon>0$, there exists $T_{\epsilon}$ such that
	\begin{align*}
		\Pr( |X_{T} | \geq \epsilon ) \leq \epsilon,
	\end{align*}
	for all $T \geq T_{\epsilon}$. Now consider the sequence $\epsilon_{l} = 1/2^{l}$ for $l \in \{0,1,\ldots\}$. For each $l$, let $T_{l}$ be the smallest $T \geq T_{l-1}$ ($T_{-1}$ is set to $0$) for which
	\begin{align*}
		\Pr( |X_{T} | \geq 1/2^{l} ) \leq 1/2^{l}.
	\end{align*}
	Such $T$ always exists by the definition of convergence in probability.
	
	Next, we construct $(r_{T})_{T=0}^{\infty}$. For any $T \in \{0,1,\ldots\}$, there exists an $l$ for which $T_{l} \leq T \leq  T_{l+1}$, so set $r_{T} = 1/2^{l}$. First, observe that for any $\delta>0$, there exists a sufficiently large $T_{\delta}$ for which the corresponding $l$ is such that $1/2^{l} < \delta$. Thus, $r_{T}$ converges to 0 as $T$ diverges. Second, take any $\epsilon>0$ and choose $l_{\epsilon}$ such that $1/2^{l_{\epsilon}} \leq \epsilon$, and finally, take $T_{\epsilon} \equiv T_{l_{\epsilon}}$. For any $T \geq T_{\epsilon}$, it follows that
	\begin{align*}
		\Pr( |X_{T} | \geq r_{T} )=	\Pr( |X_{T} | \geq 1/2^{l_{\epsilon}} ) \leq 1/2^{l_{\epsilon}} \leq \epsilon,
	\end{align*}
	as desired.
	
\end{proof}

\begin{proof}[Proof of Lemma \ref{lem:Charac.ScoreProcess}] We show that $|| \sup_{ \theta \in B(\delta,\theta_{\ast})  } || \Delta_{0,-\infty}(\theta)|| ||_{L^{2}(\bar{P}_{\ast}^{\nu})}$ is bounded ($\delta$ is as in Assumption \ref{ass:deriva-bdd})  and that  $\Delta_{t,-\infty}$ and $\Delta_{t,-\infty}\Delta_{0,-\infty} $ are continuous in $L^{1}(\bar{P}_{\ast}^{\nu})$ norm, i.e., for any $\epsilon>0$, there exists  $\delta>0$ such that
	\begin{align}\label{eqn:HAC-a1p}
		\omega(\delta) \equiv \max_{t}  ||  \sup_{||\theta  - \theta_{0}|| < \delta } \{  \Delta_{t,-\infty}(\theta)\Delta_{t,-\infty}(\theta)^{\intercal}  - \Delta_{t,-\infty}(\theta_{0}) \Delta_{t,-\infty}(\theta_{0})^{\intercal} \} ||_{L^{1}(\bar{P}_{\ast}^{\nu})} \leq \epsilon
	\end{align}
	and
	\begin{align}\label{eqn:HAC-a2p}
		\ddot{\varpi}(\delta) \equiv \max_{t} 	|| \sup_{||\theta  - \theta_{0}|| < \delta } \{  \Delta_{t,-\infty}(\theta) \Delta_{0,-\infty}(\theta)^{\intercal}   - \Delta_{t,-\infty}(\theta_{0}) \Delta_{0,-\infty}(\theta_{0})^{\intercal} \} ||_{L^{1}(\bar{P}_{\ast}^{\nu})} \leq \epsilon.
	\end{align}
We also show that there exist a constant $C < \infty$ such that, for any $t$ and $M$,
	\begin{align}\label{eqn:HAC-b1p}
		|| \sup_{  \theta \in B(\delta,\theta_{\ast})   } \Delta_{t,-\infty}(\theta)   - \Delta_{t,t-M}(\theta) ||_{L^{2}(\bar{P}_{\ast}^{\nu})} \leq C M^{1-1/p}.
	\end{align}
	Moreover, by Assumption \ref{ass:q-sum-p}, $p \in (0,2/3)$ and thus the RHS vanishes as $M$ diverges.

	We first establish (\ref{eqn:HAC-b1p}). To do so, observe that, by inspection of the proof of Lemma \ref{lem:score_approx}, it follows that the conclusion of that lemma holds uniformly in $\theta$ (and also in $t$), i.e.,
	\begin{align}\label{eqn:HAC-b}
		|| \sup_{  \theta \in B(\delta,\theta_{\ast})   } \Delta_{t,-\infty}(\theta)   - \Delta_{t,t-M}(\theta) ||_{L^{2}(\bar{P}_{\ast}^{\nu})}  \precsim \sum_{j=[t-M/2]}^{t-1} \varrho(j,t-M) + \sum_{j=[t-M]}^{[t-M/2]-1} \varrho(t-1,j).
	\end{align}
By the definition of $\varrho$ and stationarity, we have that, for any $j \geq k$, $\varrho(j,k) = \varrho(j-k,0)$ and thus $\sum_{j=[t-M/2]}^{t-1} \varrho(j,t-M) \leq \sum_{j=[t-M/2]}^{t-1} 1/( j - (t-M)  )^{1/p} \leq \int_{[M/2]+1}^{M} 1/( x )^{1/p} dx \leq \frac{p}{1-p} (M/2)^{1-1/p} $, and $\sum_{j=[t-M]}^{[t-M/2]-1} \varrho(t-1-j,0) \leq \sum_{j=[M/2]-1}^{M-1} 1/( x  )^{1/p} \leq  \frac{p}{1-p} (M/2)^{1-1/p} $. Thus,
	\begin{align*}
		|| \sup_{  \theta \in B(\delta,\theta_{\ast})   } \Delta_{t,-\infty}(\theta)   - \Delta_{t,t-M}(\theta) ||_{L^{2}(\bar{P}_{\ast}^{\nu})}  \precsim M^{1-1/p},
	\end{align*}
	as desired.

	Since, under Assumption \ref{ass:deriva-bdd}, $|| \sup_{  \theta \in B(\delta,\theta_{\ast})   } || \Delta_{0,-M}(\theta)|| ||_{L^{2}(\bar{P}_{\ast}^{\nu})} < \infty$ for any finite $M$, (\ref{eqn:HAC-b1p}) implies that $|| \sup_{  \theta \in B(\delta,\theta_{\ast})  } || \Delta_{0,-\infty}(\theta)|| ||_{L^{2}(\bar{P}_{\ast}^{\nu})}$ is bounded.

	We show next that (\ref{eqn:HAC-a2p}) holds (the proof of (\ref{eqn:HAC-a1p}) is completely analogous and is, therefore, omitted). To this end, observe that
	\begin{align*}
		|| \sup_{||\theta  - \theta_{0}|| < \delta  }  \tilde{\Delta}_{t,-\infty}(\theta) - \tilde{\Delta}_{t,-\infty}(\theta_{0}) ||_{L^{1}(\bar{P}_{\ast}^{\nu})} \leq &	|| \sup_{\theta \in B(\delta,\theta_{\ast}) } \tilde{\Delta}_{t,-\infty}(\theta) - \tilde{\Delta}_{t,t-M}(\theta) ||_{L^{1}(\bar{P}_{\ast}^{\nu})} \\
		& + 	||  \sup_{ \{||\theta  - \theta_{0}|| < \delta\} \cap B(\delta,\theta_{\ast})  }   \tilde{\Delta}_{t,t-M}(\theta)  - \tilde{\Delta}_{t,t-M}(\theta_{0}) ||_{L^{1}(\bar{P}_{\ast}^{\nu})} \\
		& + 	||  \sup_{ \theta_{0}  \in B(\delta,\theta_{\ast})  }   \tilde{\Delta}_{t,-\infty}(\theta_{0})  - \tilde{\Delta}_{t,t-M}(\theta_{0}) ||_{L^{1}(\bar{P}_{\ast}^{\nu})} \\
		\equiv & Term_{1,t,M}  + Term_{2,t,M}   + Term_{3,t,M},
	\end{align*}
where, for any $t \in \{0,\ldots,T,\ldots\}$ and $M \in \mathbb{N} \cup \{ \infty\}$,   $\tilde{\Delta}_{t,-M}(\theta) \equiv  \Delta_{t,-M}(\theta) \Delta_{0,-M}(\theta)^{\intercal}$. We now bound each of these terms to obtain the desired result.
	
	Regarding terms 1 and 3, by simple algebra and the fact that $|| \sup_{\theta \in B(\delta,\theta_{\ast})  } || \Delta_{0,-M}(\theta)|| ||_{L^{2}(\bar{P}_{\ast}^{\nu})} < \infty$ for any $M \in \mathbb{N} \cup \{ \infty\}$, it follows that, for some finite $C < \infty$,
	\begin{align*}
		Term_{1,t,M}  + Term_{3,t,M}  \leq C	|| \sup_{\theta \in B(\delta,\theta_{\ast})} || \Delta_{t,-\infty}(\theta) - \Delta_{t,t-M}(\theta) ||   ||_{L^{2}(\bar{P}_{\ast}^{\nu})},
	\end{align*}
	and, by (\ref{eqn:HAC-b1p}), the RHS is bounded by $O(M^{1-1/p})$. Therefore, under Assumption \ref{ass:q-sum-p},  for any $\epsilon>0$, there exists an $M$ such that, uniformly over $t$, $Term_{1,t,M}   + Term_{3,t,M}  \leq \epsilon$. Henceforth, fix this $M$.
	
	Regarding term 2, observe that $M < \infty$ and that $\tilde{\Delta}_{t,t-M}$ is the product of two functions that are comprised of $M$-term-sums of products of $\theta \mapsto \log p_{\theta} (x,s,x')$, $\theta \mapsto \log Q_{\theta} (x,s,x')$  and their derivatives, all of which are continuous functions by Assumption \ref{ass:Theta-int}. Thus, it can be shown that $\Delta_{t,t-M}$ is continuous, thereby implying that, for any $\epsilon>0$, there exists some $\delta_{M,\epsilon}$ (which could always be chosen to be smaller than $\delta>0$) such that $Term_{2,t,M} < \epsilon$. This completes the proof.
\end{proof}

%\begin{proof}
%	Fix any $L$, then for any $j \in \{1,...,L\}$, by Birkhoff Ergodic theorem,
%\begin{align*}
%	|| T^{-1} \sum_{t=1}^{T-j}  \Delta_{t+j,-\infty} (\theta_{\ast}) (X^{t+j}_{-\infty})  \Delta_{t,-\infty} (\theta_{\ast}) (X^{t}_{-\infty}) ' - \gamma_{j}(\theta_{\ast})   ||  = o_{ \bar{P}_{\ast}^{\nu}}(1)
%\end{align*}	
%	
%%	By the Markov inequality,
%%		\begin{align*}
%%		& \bar{P}_{\ast}^{\nu} \left(   \sum_{j=0}^{L-1} \omega(j,L)  || T^{-1} \sum_{t=1}^{T-j}  \Delta_{t+j,-\infty} (\theta_{\ast}) (X^{t+j}_{-\infty})  \Delta_{t,-\infty} (\theta_{\ast}) (X^{t}_{-\infty}) ' - \gamma_{j}(\theta_{\ast})   ||   \geq  \delta   \right) \\
%%		\leq & \delta^{-1}  \sum_{j=0}^{L-1} \omega(j,L) E_{\bar{P}_{\ast}^{\nu}}  \left[   || T^{-1} \sum_{t=1}^{T-j}  \Delta_{t+j,-\infty} (\theta_{\ast}) (X^{t+j}_{-\infty})  \Delta_{t,-\infty} (\theta_{\ast}) (X^{t}_{-\infty}) ' - \gamma_{j}(\theta_{\ast})   ||      \right]
%%	\end{align*}
%\end{proof}

\section{Additional Empirical Example}

\label{SM:empiricalf}

In this empirical example, we examine the predictive ability of an index of
leading indicators for regime changes in U.S. output growth. The model, based
on \citesupp{filardo94}, allows for shifts in the mean output growth, and is
given by
\begin{align}
Y_{t}  &  =\mu_{0}+\mu_{1}S_{t}+\sum_{i=1}^{4}\phi_{i}(Y_{t-i}-\mu_{0}-\mu
_{1}S_{t-i})+\sigma_{1}U_{1,t},\label{eqyf}\\
Z_{t}  &  =\mu_{2}+\sum_{i=1}^{4}\psi_{i}Z_{t-i}+\sigma_{2}U_{2,t},
\label{eqzf}%
\end{align}
with the hidden, two-state Markov chain $(S_{t})_{t}$ being governed by the
transition probabilities
\begin{equation}
\Pr(S_{t}=s\mid Z_{t-1}=z,S_{t-1}=s)=[1+\exp(-\alpha_{s}-\beta_{s}%
z)]^{-1},\quad s\in\{0,1\}, \label{eqtpf}%
\end{equation}
and $(U_{1,t},U_{2,t})_{t}$ postulated to be i.i.d. $\mathcal{N}\left(
\left[
\begin{array}
[c]{c}%
0\\
0
\end{array}
\right]  ,\left[
\begin{array}
[c]{cc}%
1 & \rho\\
\rho & 1
\end{array}
\right]  \right)  $ and independent of $(S_{t})_{t}$.
%\footnote{It should be noted that our asymptotic results do not allow for a structure like (\ref{eqyf}) due to the dependence of the conditional distribution of $Y_{t}$ on $S_{t-i}$, $i>0$. However, extending our theory to this case is straightforward.}
Here, $Y_{t}$ stands for the growth rate of total industrial production and
$Z_{t}$ stands for the growth rate in the Composite Index of Eleven Leading
Indicators (CLI). The data, taken from \citesupp[Ch. 4]{kim99}, are monthly,
span the period 1948:1--1991:4, and are transformed in the way detailed in
\citesupp{filardo94}. Note that the analysis in \cite{filardo94} is based on a
slightly longer sample period ending in 1992:8. The likelihood ratio test of
\citesupp{Hansen92} rejects the hypothesis that $\mu_{0}=\mu_{1}$ in
(\ref{eqyf}).

{\footnotesize \ \begin{table}[ptb]
\caption{ML Estimates (Real Output, Index of Leading Indicators)}%
\label{tab:app-3}%
{\footnotesize  \ \ \  }
\par
\begin{center}
{\footnotesize \ \ \ \
%\multicolumn{1}{l}{%
\begin{tabular}
[c]{lllllllll}\hline\hline
\multicolumn{2}{c}{Partial ML} &  &  &  & \multicolumn{4}{c}{Joint ML}\\\hline
$\mu_{0}$ & -0.8567 (0.3257 ) &  &  &  & $\mu_{0}$ & -0.7892 (0.2786) &  & \\
$\mu_{1}$ & 1.3762 (0.3257) &  &  &  & $\mu_{1}$ & 1.3099 (0.2179) & $\mu_{2}$
& 0.0090 (0.0678)\\
$\phi_{1}$ & 0.1949 (0.0784) &  &  &  & $\phi_{1}$ & 0.1952 (0.0694) &
$\psi_{1}$ & 0.3769 (0.0490)\\
$\phi_{2}$ & 0.0802 (0.0540) &  &  &  & $\phi_{2}$ & 0.1278 (0.0493) &
$\psi_{2}$ & 0.1891 (0.0538 )\\
$\phi_{3}$ & 0.1094 (0.0602) &  &  &  & $\phi_{3}$ & 0.1161 (0.0453) &
$\psi_{3}$ & 0.0049 (0.0486 )\\
$\phi_{4}$ & 0.1229 (0.0611) &  &  &  & $\phi_{4}$ & 0.1169 (0.0484) &
$\psi_{4}$ & 0.0570 (0.0515 )\\
$\alpha_{0}$ & 1.6311 (0.5351) &  &  &  & $\alpha_{0}$ & 1.7158 (0.6145) &
$\sigma_{2}$ & 0.7729 (0.0437)\\
$\beta_{0}$ & -1.0428 (0.7048) &  &  &  & $\beta_{0}$ & -0.8408 (0.4131) &
$\rho$ & 0.4498 (0.0345)\\
$\alpha_{1}$ & 4.3743 (0.6891 ) &  &  &  & $\alpha_{1}$ & 4.1188 (0.7618) &  &
\\
$\beta_{1}$ & 1.7669 (0.4805) &  &  &  & $\beta_{1}$ & 1.6020 (0.8209) &  & \\
$\sigma_{1}$ & 0.6969 (0.0386) &  &  &  & $\sigma_{1}$ & 0.6928 (0.0393) &  &
\\\hline\hline
\end{tabular}
}
\end{center}
\end{table}}

We compute two sets of ML estimates: partial ML estimates based on
(\ref{eqyf}) and (\ref{eqtpf}) (as in \citesupp{filardo94}) and joint ML
estimates based on the system (\ref{eqyf})--(\ref{eqtpf}). Results are
presented in Table \ref{tab:app-3}, with estimated standard errors given in
parentheses. The latter are obtained from the \textquotedblleft
sandwich\textquotedblright\ estimator defined in Theorem~\ref{thm:StdErrors}%
(b); the weights $\omega(\cdot,L_{T})$ are generated from the Parzen kernel
and the truncation parameter $L_{T}$ is determined using the automatic plug-in
method of \citesupp{andrews91}. The differences between partial and joint ML
estimates are not very substantial, a finding which is not perhaps surprising
in view of the fact that the estimated conditional correlation $\rho$ is
moderate (0.4498). Conventional $t$-type tests based on joint ML estimates
reject the hypotheses that $\beta_{0}$ and $\beta_{1}$ are zero, suggesting
that the lagged CLI has significant information about the probability of
switching between the two Markov regimes. Needless to say, inference based on
the partial ML estimator must be viewed with caution unless $\rho=0$.
%XXXXIn addition, it should be noted that, without appropriate generalization, our asymptotic results do not allow for a structure like (\ref{eqyf}) due to the dependence of the conditional distribution of $Y_{t}$ on $S_{t-i}$, $i>0$. XXXX \emph{TRUE}?

\begin{figure}[h]
\caption{Difference between the inferred probability of $S_{t}=1$ based on
partial ML and joint ML estimates.}%
\label{probs}
\centering
\includegraphics[width=\linewidth,height=7cm]{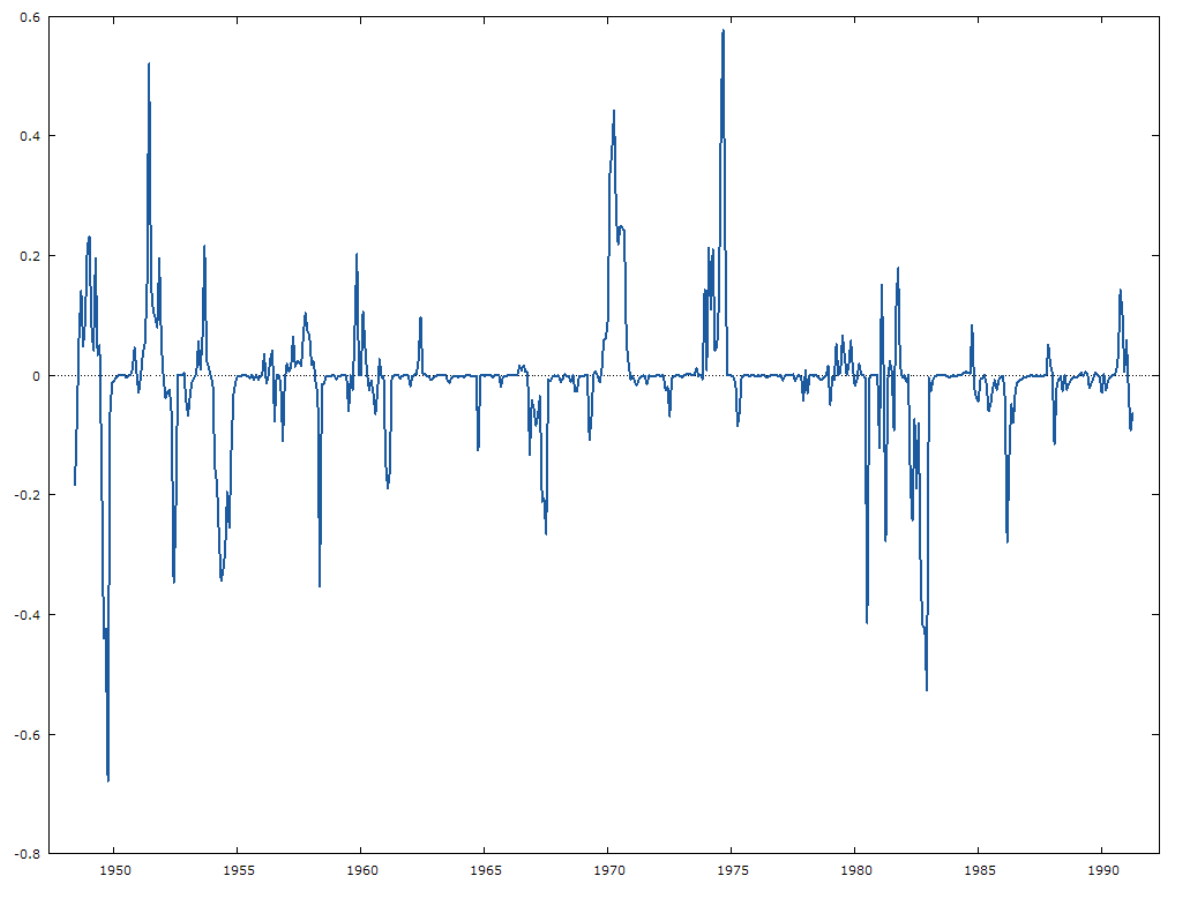}\end{figure}

It is worth noting that, although partial ML and joint ML estimates of the
parameters are not substantially different, the inferred probabilities that
$S_{t}=1$, conditional on sample information available at time $t$, based on
these estimates are not always close. This can be seen in Figure~\ref{probs},
which plots the difference between the inferred probabilities based on the two
sets of estimates. The estimated covariate-dependent transition probabilities
based on the two sets of ML estimates (not shown) are, on the other hand,
fairly similar.

\newpage

\bibliographystylesupp{plainnat} {\small {
\bibliographystylesupp{plain}
\bibliographysupp{markov}
}}

\end{document}